\documentclass[11pt]{article}
\usepackage[english]{babel}
\usepackage{times}
\usepackage[T1]{fontenc}
\usepackage{amsmath, amsthm, amsfonts}
\usepackage{amssymb}
\usepackage{amscd}
\usepackage{xspace}
\usepackage{verbatim}
\usepackage{latexsym}
\usepackage{amsfonts}
\usepackage{graphicx}
\usepackage{float}
\usepackage{color}

\newcommand{\mysection}[1]{
\section{#1}\setcounter{equation}{0}}
\title{\bf Boundary singularities of  positive solutions of quasilinear Hamilton-Jacobi equations\footnote{To appear in Calc. Var. Part. Diff. Equ.}}
\author{{\bf Marie-Fran\c{c}oise Bidaut-V\'eron} \\{\bf Marta Garcia-Huidobro}\\
 {\bf Laurent V\'eron}\\[2mm]
}

\date{}
\begin{document}
 \maketitle


\newcommand{\txt}[1]{\;\text{ #1 }\;}
\newcommand{\tbf}{\textbf}
\newcommand{\tit}{\textit}
\newcommand{\tsc}{\textsc}
\newcommand{\trm}{\textrm}
\newcommand{\mbf}{\mathbf}
\newcommand{\mrm}{\mathrm}
\newcommand{\bsym}{\boldsymbol}
\newcommand{\scs}{\scriptstyle}
\newcommand{\sss}{\scriptscriptstyle}
\newcommand{\txts}{\textstyle}
\newcommand{\dsps}{\displaystyle}
\newcommand{\fnz}{\footnotesize}
\newcommand{\scz}{\scriptsize}
\newcommand{\be}{\begin{equation}}
\newcommand{\bel}[1]{\begin{equation}\label{#1}}
\newcommand{\ee}{\end{equation}}
\newcommand{\eqnl}[2]{\begin{equation}\label{#1}{#2}\end{equation}}
\newcommand{\barr}{\begin{eqnarray}}
\newcommand{\earr}{\end{eqnarray}}
\newcommand{\bars}{\begin{eqnarray*}}
\newcommand{\ears}{\end{eqnarray*}}
\newcommand{\nnu}{\nonumber \\}
\newtheorem{subn}{\name}
\renewcommand{\thesubn}{}
\newcommand{\bsn}[1]{\def\name{#1}\begin{subn}}
\newcommand{\esn}{\end{subn}}
\newtheorem{sub}{\name}[section]
\newcommand{\dn}[1]{\def\name{#1}}   
\newcommand{\bs}{\begin{sub}}
\newcommand{\es}{\end{sub}}
\newcommand{\bsl}[1]{\begin{sub}\label{#1}}
\newcommand{\bth}[1]{\def\name{Theorem}
\begin{sub}\label{t:#1}}
\newcommand{\blemma}[1]{\def\name{Lemma}
\begin{sub}\label{l:#1}}
\newcommand{\bcor}[1]{\def\name{Corollary}
\begin{sub}\label{c:#1}}
\newcommand{\bdef}[1]{\def\name{Definition}
\begin{sub}\label{d:#1}}
\newcommand{\bprop}[1]{\def\name{Proposition}
\begin{sub}\label{p:#1}}
\newcommand{\R}{\eqref}
\newcommand{\rth}[1]{Theorem~\ref{t:#1}}
\newcommand{\rlemma}[1]{Lemma~\ref{l:#1}}
\newcommand{\rcor}[1]{Corollary~\ref{c:#1}}
\newcommand{\rdef}[1]{Definition~\ref{d:#1}}
\newcommand{\rprop}[1]{Proposition~\ref{p:#1}}
\newcommand{\BA}{\begin{array}}
\newcommand{\EA}{\end{array}}
\newcommand{\BAN}{\renewcommand{\arraystretch}{1.2}
\setlength{\arraycolsep}{2pt}\begin{array}}
\newcommand{\BAV}[2]{\renewcommand{\arraystretch}{#1}
\setlength{\arraycolsep}{#2}\begin{array}}
\newcommand{\BSA}{\begin{subarray}}
\newcommand{\ESA}{\end{subarray}}
\newcommand{\BAL}{\begin{aligned}}
\newcommand{\EAL}{\end{aligned}}
\newcommand{\BALG}{\begin{alignat}}
\newcommand{\EALG}{\end{alignat}}
\newcommand{\BALGN}{\begin{alignat*}}
\newcommand{\EALGN}{\end{alignat*}}
\newcommand{\note}[1]{\textit{#1.}\hspace{2mm}}
\newcommand{\Proof}{\note{Proof}}
\newcommand{\qeda}{\hspace{10mm}\hfill $\square$}
\newcommand{\Remark}{\note{Remark}}
\newcommand{\modin}{$\,$\\[-4mm] \indent}
\newcommand{\forevery}{\quad \forall}
\newcommand{\set}[1]{\{#1\}}
\newcommand{\setdef}[2]{\{\,#1:\,#2\,\}}
\newcommand{\setm}[2]{\{\,#1\mid #2\,\}}
\newcommand{\mt}{\mapsto}
\newcommand{\lra}{\longrightarrow}
\newcommand{\lla}{\longleftarrow}
\newcommand{\llra}{\longleftrightarrow}
\newcommand{\Lra}{\Longrightarrow}
\newcommand{\Lla}{\Longleftarrow}
\newcommand{\Llra}{\Longleftrightarrow}
\newcommand{\warrow}{\rightharpoonup}
\newcommand{
\paran}[1]{\left (#1 \right )}
\newcommand{\sqbr}[1]{\left [#1 \right ]}
\newcommand{\curlybr}[1]{\left \{#1 \right \}}
\newcommand{\abs}[1]{\left |#1\right |}
\newcommand{\norm}[1]{\left \|#1\right \|}
\newcommand{
\paranb}[1]{\big (#1 \big )}
\newcommand{\lsqbrb}[1]{\big [#1 \big ]}
\newcommand{\lcurlybrb}[1]{\big \{#1 \big \}}
\newcommand{\absb}[1]{\big |#1\big |}
\newcommand{\normb}[1]{\big \|#1\big \|}
\newcommand{
\paranB}[1]{\Big (#1 \Big )}
\newcommand{\absB}[1]{\Big |#1\Big |}
\newcommand{\normB}[1]{\Big \|#1\Big \|}
\newcommand{\produal}[1]{\langle #1 \rangle}

\newcommand{\thkl}{\rule[-.5mm]{.3mm}{3mm}}
\newcommand{\thknorm}[1]{\thkl #1 \thkl\,}
\newcommand{\trinorm}[1]{|\!|\!| #1 |\!|\!|\,}
\newcommand{\bang}[1]{\langle #1 \rangle}
\def\angb<#1>{\langle #1 \rangle}
\newcommand{\vstrut}[1]{\rule{0mm}{#1}}
\newcommand{\rec}[1]{\frac{1}{#1}}
\newcommand{\opname}[1]{\mbox{\rm #1}\,}
\newcommand{\supp}{\opname{supp}}
\newcommand{\dist}{\opname{dist}}
\newcommand{\myfrac}[2]{{\displaystyle \frac{#1}{#2} }}
\newcommand{\myint}[2]{{\displaystyle \int_{#1}^{#2}}}
\newcommand{\mysum}[2]{{\displaystyle \sum_{#1}^{#2}}}
\newcommand {\dint}{{\displaystyle \myint\!\!\myint}}
\newcommand{\q}{\quad}
\newcommand{\qq}{\qquad}
\newcommand{\hsp}[1]{\hspace{#1mm}}
\newcommand{\vsp}[1]{\vspace{#1mm}}
\newcommand{\ity}{\infty}
\newcommand{\prt}{\partial}
\newcommand{\sms}{\setminus}
\newcommand{\ems}{\emptyset}
\newcommand{\ti}{\times}
\newcommand{\pr}{^\prime}
\newcommand{\ppr}{^{\prime\prime}}
\newcommand{\tl}{\tilde}
\newcommand{\sbs}{\subset}
\newcommand{\sbeq}{\subseteq}
\newcommand{\nind}{\noindent}
\newcommand{\ind}{\indent}
\newcommand{\ovl}{\overline}
\newcommand{\unl}{\underline}
\newcommand{\nin}{\not\in}
\newcommand{\pfrac}[2]{\genfrac{(}{)}{}{}{#1}{#2}}

\def\ga{\alpha}     \def\gb{\beta}       \def\gg{\gamma}
\def\gc{\chi}       \def\gd{\delta}      \def\ge{\epsilon}
\def\gth{\theta}                         \def\vge{\varepsilon}
\def\gf{\phi}       \def\vgf{\varphi}    \def\gh{\eta}
\def\gi{\iota}      \def\gk{\kappa}      \def\gl{\lambda}
\def\gm{\mu}        \def\gn{\nu}         \def\gp{\pi}
\def\vgp{\varpi}    \def\gr{\rho}        \def\vgr{\varrho}
\def\gs{\sigma}     \def\vgs{\varsigma}  \def\gt{\tau}
\def\gu{\upsilon}   \def\gv{\vartheta}   \def\gw{\omega}
\def\gx{\xi}        \def\gy{\psi}        \def\gz{\zeta}
\def\Gg{\Gamma}     \def\Gd{\Delta}      \def\Gf{\Phi}
\def\Gth{\Theta}
\def\Gl{\Lambda}    \def\Gs{\Sigma}      \def\Gp{\Pi}
\def\Gw{\Omega}     \def\Gx{\Xi}         \def\Gy{\Psi}

\def\CS{{\mathcal S}}   \def\CM{{\mathcal M}}   \def\CN{{\mathcal N}}
\def\CR{{\mathcal R}}   \def\CO{{\mathcal O}}   \def\CP{{\mathcal P}}
\def\CA{{\mathcal A}}   \def\CB{{\mathcal B}}   \def\CC{{\mathcal C}}
\def\CD{{\mathcal D}}   \def\CE{{\mathcal E}}   \def\CF{{\mathcal F}}
\def\CG{{\mathcal G}}   \def\CH{{\mathcal H}}   \def\CI{{\mathcal I}}
\def\CJ{{\mathcal J}}   \def\CK{{\mathcal K}}   \def\CL{{\mathcal L}}
\def\CT{{\mathcal T}}   \def\CU{{\mathcal U}}   \def\CV{{\mathcal V}}
\def\CZ{{\mathcal Z}}   \def\CX{{\mathcal X}}   \def\CY{{\mathcal Y}}
\def\CW{{\mathcal W}} \def\CQ{{\mathcal Q}}
\def\BBA {\mathbb A}   \def\BBb {\mathbb B}    \def\BBC {\mathbb C}
\def\BBD {\mathbb D}   \def\BBE {\mathbb E}    \def\BBF {\mathbb F}
\def\BBG {\mathbb G}   \def\BBH {\mathbb H}    \def\BBI {\mathbb I}
\def\BBJ {\mathbb J}   \def\BBK {\mathbb K}    \def\BBL {\mathbb L}
\def\BBM {\mathbb M}   \def\BBN {\mathbb N}    \def\BBO {\mathbb O}
\def\BBP {\mathbb P}   \def\BBR {\mathbb R}    \def\BBS {\mathbb S}
\def\BBT {\mathbb T}   \def\BBU {\mathbb U}    \def\BBV {\mathbb V}
\def\BBW {\mathbb W}   \def\BBX {\mathbb X}    \def\BBY {\mathbb Y}
\def\BBZ {\mathbb Z}

\def\GTA {\mathfrak A}   \def\GTB {\mathfrak B}    \def\GTC {\mathfrak C}
\def\GTD {\mathfrak D}   \def\GTE {\mathfrak E}    \def\GTF {\mathfrak F}
\def\GTG {\mathfrak G}   \def\GTH {\mathfrak H}    \def\GTI {\mathfrak I}
\def\GTJ {\mathfrak J}   \def\GTK {\mathfrak K}    \def\GTL {\mathfrak L}
\def\GTM {\mathfrak M}   \def\GTN {\mathfrak N}    \def\GTO {\mathfrak O}
\def\GTP {\mathfrak P}   \def\GTR {\mathfrak R}    \def\GTS {\mathfrak S}
\def\GTT {\mathfrak T}   \def\GTU {\mathfrak U}    \def\GTV {\mathfrak V}
\def\GTW {\mathfrak W}   \def\GTX {\mathfrak X}    \def\GTY {\mathfrak Y}
\def\GTZ {\mathfrak Z}   \def\GTQ {\mathfrak Q}

\font\Sym= msam10 
\def\SYM#1{\hbox{\Sym #1}}
\newcommand{\bdw}{\prt\Gw\xspace}
\date{}
\maketitle\medskip

\noindent{\small {\bf Abstract}
We study the boundary behaviour of the solutions of (E) $\;-\Gd_p u+|\nabla u|^q=0$ in a domain $\Gw \sbs \BBR^N$, when $N\geq p> q>p-1$. We show the existence of a critical exponent $q_*<p$ such that if $p-1<q<q_*$ there exist positive solutions of (E) with an isolated singularity on $\prt\Gw$ and that these solutions belong to two different classes of singular solutions. If $q_*\leq q<p$ no such solution exists and actually any boundary isolated singularity of a positive solution of (E) is removable. We prove that all the singular positive solutions are classified according the two types of singular solutions that we have constructed.
}\smallskip

\tableofcontents
\vspace{1mm}
\hspace{.05in}
\medskip
\mysection{Introduction}
Let  $N\geq p>1$, $q>p-1$ and $\Gw \sbs \BBR^N$ ($N>1$) be a $C^2$ bounded domain such that $0\in\prt\Gw$. In this article we study the boundary behavior at $0$ of nonnegative functions $u\in C^1(\Gw)\cap C(\overline\Gw\setminus\{0\})$ which satisfy
\be\label{A1}\left.\BA{ll}
-\Gd_p u+|\nabla u|^q=0\qquad\text{in }\Gw,
\EA\right.\ee
where $\Gd_pu:={\rm div}\left(|\nabla u|^{p-2}\nabla u\right)$.
The two main questions we consider are as follows:\smallskip

 \noindent {\bf Q-1}- Existence of positive solutions of  $(\ref{A1})$. \smallskip

 \noindent {\bf Q-2}- Description of positive solutions with an isolated boundary singularity at $0$.\medskip

 When $p=2$ a fairly complete description of positive solutions of
 \be\label{A1-0}\left.\BA{ll}
-\Gd u+|\nabla u|^q=0\\
\EA\right.\ee
in $\Gw$ is provided by Nguyen-Phuoc and V\'eron \cite{NpV}. In particular they prove the following series of results in the range of values $1<q<2$.\smallskip

\noindent {\bf 1}- Any signed solution of $(\ref{A1-1})$ verifies the estimates

 \be\label{A1-1}\left.\BA{ll}
|\nabla u(x)|\leq c_{N,q}\left(d(x)\right)^{-\frac{1}{q-1}}\qquad\forall x\in\Gw,\\
\EA\right.\ee
where $d(x)=\dist(x,\prt\Gw)$. As a consequence, if $u\in C(\overline\Gw\setminus\{0\})$ is a solution which vanishes on $\prt \Gw\setminus\{0\}$, it satisfies
 \be\label{A1-2}\left.\BA{ll}
|u(x)|\leq c_{q,\Gw}d(x)|x|^{-\frac{1}{q-1}}\qquad\forall x\in\Gw.\\
\EA\right.\ee

\noindent {\bf 2}- If $\frac{N+1}{N}\leq q<2$ any positive solution of $(\ref{A1-1})$ in $\Gw$  which vanishes on $\prt\Gw\setminus\{0\}$ is identically $0$.  {\it An isolated boundary point is a removable singularity for $(\ref{A1-0})$}. \smallskip

\noindent {\bf 3}- If $1<q<\frac{N+1}{N}$  and $k>0$ there exists a unique positive solution $u:=u_{k}$ of $(\ref{A1-0})$ in $\Gw$  which vanishes on $\prt\Gw\setminus\{0\}$ and satisfies $u(x)\sim c_N kP^\Gw(x,0)$ as $x\to 0$, where $P^\Gw$ is the Poisson kernel in $\Gw\ti\prt\Gw$.
\smallskip

\noindent {\bf 4}- If $1<q<\frac{N+1}{N}$ there exists a unique positive solution $u$ of  $(\ref{A1-0})$ in the half-space $\BBR^N_+:=\{x=(x',x_N):x'\in\BBR^{N-1}, x_N>0\}$ under the form
$u(x)=|x|^{-\frac{2-q}{q-1}}\gw(|x|^{-1}x)$ which vanishes on $ \prt\BBR^N_+\setminus\{0\}$. The function $\gw$ is the unique positive solution of
\be\label{A1-4}\left.\BA{ll}
-\Gd' \gw+\left((\frac{2-q}{q-1})^2\gw^2+|\nabla'\gw|^2\right)^{\frac q 2}-\gl_{N,q}\gw=0\qquad\text{in }S^{N-1}_+,\\
\phantom{-\Gd' u+((q-1)^2\gw^2+|\nabla'\gw|^2)^{\frac q 2}-\gl_{N,q}}
\gw=0\qquad\text{in }\prt S^{N-1}_+,
\EA\right.\ee
where $S^{N-1}$ is the unit sphere  of $\BBR^N$, $\prt S^{N-1}_+=\prt\BBR^N_+\cap S^{N-1}$, $\Gd'$ the Laplace-Beltrami operator  and $\gl_{N,q}>0$ an explicit constant.

\smallskip

\noindent {\bf 5}- If $1<q<\frac{N+1}{N}$ and $u$ is a positive solution of $(\ref{A1-1})$ in $\Gw$, which is continuous in $\overline \Gw\setminus\{0\}$ and vanishes on $\prt\Gw\setminus\{0\}$ the following dichotomy occurs: \smallskip

\noindent {\it (i)} either $u(x)\sim |x|^{-\frac{2-q}{q-1}}\gw(|x|^{-1}x)$ as $x\to 0$,\smallskip

\noindent {\it (ii)} or $u(x)\sim kc_NP^{\Gw}(x,0)$ as $x\to 0$ for some $k\geq 0$.\medskip

The aim of this article is to extend to the quasilinear case $1<p\leq N$ the above mentioned results. The following pointwise gradient estimate valid for any signed solution $u$ of $(\ref{A1})$ has been proved in \cite{BVGHV1}: if $0<p-1<q$ there exists a constant  $c_{N,p,q}>0$ such that
 \bel{I-1}
 \abs{\nabla u(x)}\leq c_{N,p,q}(d(x))^{-\frac{1}{q+1-p}}\qquad\forall x\in \Gw.
 \ee

As a consequence, any solution $u\in C^1(\overline \Gw\setminus\{0\}$ satisfies

 \bel{I-2}
 \abs{ u(x)}\leq c_{p,q,\Gw}d(x)\abs x^{-\frac{1}{q+1-p}}\qquad\forall x\in \Gw.
 \ee

Concerning boundary singularities, the situation is much more complicated than in the case $p=2$ and the threshold of critical exponent less explicit. We first consider the problem in  $\BBR^N_+$. Assuming $p-1<q\leq p$, separable solutions of $(\ref{A1})$ in $\BBR^N_+$ vanishing on $\partial\BBR^N_+\setminus\{0\}$ can be looked for in spherical coordinates $(r,\gs)\in\BBR^*_+\ti S^{N-1}$ (we denote $\BBR^*_+=(0,\infty)$) under the form
   \bel{I-7}
u(x)=u(r,\gs)=r^{- \gb_q}\gw(\gs),\quad r>0\,,\;\gs\in S^{N-1}_+:=\{S^{N-1}\cap \BBR^N_+\}.
\ee
Then $\gw$ is solution of the following problem
   \bel{I-8}\BA {ll}
-div'\left(\left(\gb_q^2\gw^2+|\nabla' \gw|^2\right)^{\frac{p-2}2}\nabla' \gw\right)-\gb_q\Gl_{\gb_q}\left(\gb_q^2\gw^2+|\nabla' \gw|^2\right)^{\frac{p-2}2}\gw
\\[2mm]
\phantom{----------------}
+\left(\gb_q^2\gw^2+|\nabla' \gw|^2\right)^{\frac{q}2}=0\quad \mbox{in } S_+^{N-1}
\\[2mm]
\phantom{------------------------,}
\gw=0\quad \mbox{on } \prt S_+^{N-1},
\EA\ee
where
   \bel{I-8*}
   \gb_q=\frac{p-q}{q+1-p}\,\text{ and }\;\Gl_{\gb_q}=\gb_q(p-1)+p-N,
   \ee
and $\nabla'$ is the covariant derivative on $S^{N-1}$ identified to the tangential gradient thanks to the canonical isometrical imbedding of $S^{N-1}$ into $\BBR^N$, and $div'$ the divergence operator acting on vector fields on $S^{N-1}$.

The existence of a positive solution to this problem cannot be separated from the problem of existence of  {\it separable $p$-harmonic functions} which are
 $p$-harmonic  in $\BBR^N_+$ which vanish on $\prt\BBR^N_+\setminus\{0\}$ and have the form
 $\Psi(x)=\Psi(r,\gs)=r^{-\gb}\psi(\gs)$ for some real number $\gb$. Necessarily such a $\psi$ must satisfy
    \bel{I-9}\BA {ll}
-div'\left(\left(\gb^2\psi^2+|\nabla' \psi|^2\right)^{\frac{p-2}2}\nabla' \psi\right)-\gb\Gl_{\gb}\left(\gb^2\psi^2+|\nabla' \psi|^2\right)^{\frac{p-2}2}\psi=0\quad& \mbox{in } S_+^{N-1}
\\[2mm]
\phantom{-div'\left(\left(\gb^2\psi^2+|\nabla' \psi|^2\right)^{\frac{p-2}2}\nabla' \psi\right)-\gb\Gl_{\gb}\left(\gb^2\psi^2+|\nabla' \psi|^2\right)^{\frac{p-2}2}}
\psi=0\quad &\mbox{on } \prt S_+^{N-1},
\EA\ee
where $\Gl_{\gb}=\gb(p-1)+p-N$.
We will refer to $(\ref{I-9})$ as {\it the spherical $p$-harmonic eigenvalue problem}. The study of this problem has been initiated in the 2-dim case by Krol \cite{Kr} ($\gb<0$) and Kichenassamy and V\'eron \cite {KV} ($\gb>0$). In this case $\gw$ satisfies a completely integrable second order differential equation. In the case where $S^{N-1}_+$ is replaced by a smooth domain $S\subset S^{N-1}$ with $N\geq 3$, Tolksdorf \cite{To} proved the existence of a unique couple $(\tilde \gb_s,\tilde\psi_s)$ where $\tilde \gb_s<0$ and $\tilde\psi_s$ has constant sign and is defined up to an homothety. Recently Porretta and V\'eron \cite {PV1} gave a simpler and more general proof of the existence of two couples
 $(\tilde \gb_s,\tilde\psi_s)$ and $(\gb_{*\,s},\psi_{*\,s})$ where $\gb_{*\,s}>0$ and $\tilde\psi_s$ and $\psi_{*\,s}$ are positive solutions of $\eqref{I-9}$ with $\gb=\tilde\gb_s$ and $\gb=\gb_{*\,s}$ respectively and are unique up to a multiplication by a real number. When $p=2$ this problem is an eigenvalue problem for the Laplace-Beltrami operator on a subdomain of $S^{N-1}$. If $S=S^{N-1}_+$, $\tilde \gb_s$ and  $\gb_{*\,s}$ are respectively denoted by $\tilde \gb$ and  $\gb_*$ and accordingly $\tilde\psi_s$ and $\psi_{*\,s}$ by $\tilde\psi$ and $\psi_*$. Since $x\mapsto x_N$ is $p$-harmonic,  $\tilde \gb=-1$. Except in the cases $N=2$ where it is the positive root of some algebraic equation of degree 2, $p=2$ where it is $N-1$ and $p=N$ where it is $1$, the value of $\gb_*$ is unknown besides the straightforward estimate $\gb_*\geq \max\{1,\frac{N-p}{p-1}\}$. Using the fact that $\psi_*$ depends only on the azimuthal variable and satisfies a differential equation,  we prove in Appendix II  the following new estimate:
\medskip

\noindent {\bf Theorem A }{\it Let $1<p\leq N$.  \smallskip

\noindent (i) If $2\leq p\leq N$, then $\gb_*\leq \frac{N-1}{p-1}$ with equality only if $p=2$ or $N$.
\smallskip

\noindent (ii) If $1\leq p<2$, then  $\gb_*> \frac{N-1}{p-1}$.}\medskip

The $p$-harmonic function $\Psi_*(x)=\Psi_*(r,\gs)=r^{-\gb_*}\psi_*(\gs)$ endows the role of a Poisson kernel. To this exponent $\gb_*$ is associated the critical value $q_*$ of $q$ defined by $\gb_*=\gb_q$, or equivalently
  \bel{I-10}
  q_*:=\frac{\gb_*(p-1)+p}{\gb_*+1}=p-\frac{\gb_*}{\gb_*+1}.
\ee

The following result characterizes strong singularities.\medskip

\noindent {\bf Theorem B }{\it Let $0<p-1\leq N$, then\smallskip

\noindent (i) If $p-1<q<q_*$ problem $(\ref{I-8})$ admits a unique positive solution $\gw_*$.
\smallskip

\noindent (ii) If $q_*\leq q <p$ problem $(\ref{I-8})$ admits no positive solution.}\medskip

This critical exponent corresponds to the threshold of criticality for boundary isolated singularities. \medskip

\noindent {\bf Theorem C }{\it Assume $q_*\leq q <p\leq N$. If $u\in C^1(\overline\Gw\setminus\{0\})$ is a nonnegative solution of $(\ref{A1})$ in $\Gw$ which vanishes on $\prt\Gw\setminus\{0\}$, it is identical zero.}\medskip

As in the case $p=2$, there exist positive solutions $(\ref{A1})$ in $\Gw$ with weak boundary singularities which are characterized by their blow-up near the singularity. By opposition to the case $p=2$ where existence is obtained by use of a weak formulation of the boundary value problem, combined with uniform integrability of the absorption term thanks to Poisson kernel estimates (see \cite{NpV}), this approach cannot be performed in the case $p\neq 2$; the obtention of solutions with weak singularities necessitates a very  long and delicate construction of subsolutions and supersolutions. Furthermore, when $p\neq N$,  the construction  is done only if $\Gw$ is locally an hyperplane near $0$.\medskip

In the sequel we denote by $B_R(a)$ the open ball of center $a$ and  radius $R>0$ and $B_R=B_R(0)$. We also set $B^+_R(a):=\BBR^N_+\cap B_R(a)$, $B^+_R:=\BBR^N_+\cap B_R$, $B^-_R(a):=\BBR^N_-\cap B_R(a)$ and $B^-_R:=\BBR^N_-\cap B_R$, where
$\BBR^N_-:=\{x=(x',x_N):x'\in\BBR^{N-1}, x_N<0\}$. If
$\Gw$ is an open domain and $R>0$, we put $\Gw_R=\Gw\cap B_R$ .\medskip

\noindent {\bf Theorem D }{\it Let $\Gw\subset\BBR^N_+ $ be a bounded domain such that $0\in\prt\Gw$. Assume there exists $\gd>0$  such that
$\Gw_\gd=B^+_\gd$ and $0<p-1<q<q_* <p\leq N$. Then for any $k>0$ there exists a unique $u:=u_k\in C^1(\overline\Gw\setminus\{0\})$, solution of  $(\ref{A1})$ in $\Gw$, vanishing on
$\prt\Gw\setminus\{0\}$ and such that
  \bel{I-11}
\lim_{\tiny\BA{cl}x\to 0\\\frac {x}{|x|}\to\gs\in S^{N-1}_+
\EA}\!\!\!\!|x|^{\gb_*}u_k(x)=k\psi_*(\gs).
\ee
Furthermore $\lim_{k\to \infty}u_{k}=u_\infty$ and
  \bel{I-12}
\lim_{\tiny\BA{cl}x\to 0\\\frac {x}{|x|}\to\gs\in S^{N-1}_+
\EA}\!\!\!\!|x|^{\gb_q}u_\infty(x)=\psi_*(\gs).
\ee
}\medskip

When $p=N$, then $q_*=N-\frac{1}{2}$; in such a range of values we use the conformal invariance of  $\Gd_N$ and prove that the previous result holds if $\Gw$ is {\it any} $C^2$ domain. Finally, the isolated singularities of positive solutions of $(\ref{A1})$ are completely described by the two types of singular solutions obtained in the previous theorem and we prove:\medskip

\noindent {\bf Theorem E }{\it Let $\Gw $ be a bounded domain such that $0\in\prt\Gw$. Assume there exists $\gd>0$  such that
$\Gw_\gd=B^+_\gd$ and $0<p-1<q<q_* <p\leq N$. If $u\in C^1(\overline\Gw\setminus\{0\})$ is a positive solution of  $(\ref{A1})$ in $\Gw$ which vanishes on
$\prt\Gw\setminus\{0\}$, then\smallskip

\noindent (i) either there exists $k\geq 0$ such that
  \bel{I-13}
\lim_{\tiny\BA{cl}x\to 0\\\frac {x}{|x|}\to\gs\in S^{N-1}_+
\EA}\!\!\!\!|x|^{\gb_*}u(x)=k\psi_*(\gs);
\ee

\noindent (ii) or
  \bel{I-14}
\lim_{\tiny\BA{cl}x\to 0\\\frac {x}{|x|}\to\gs\in S^{N-1}_+
\EA}\!\!\!\!|x|^{\gb_q}u(x)=\psi_*(\gs).
\ee
}\medskip

\noindent{\bf Acknowledgements} This article has been prepared with the support of the MathAmsud collaboration program 13MATH-03 QUESP. The first two authors were supported by Fondecyt grant N$^{\circ}$1110268. The authors are grateful to the referee for a careful reading of the manuscript.

\mysection{A priori estimates}
\subsection{The gradient estimates and its applications}
We recall the following estimate and its consequences which are proved in \cite{BVGHV1}.

\bprop{lets0} Assume $q>p-1$ and $u$ is a $C^1$ solution of $(\ref{A1})$ in a domain $\Gw$. Then
\bel{est0}
\abs{\nabla u(x)}\leq c_{N,p,q}(d(x))^{-\frac{1}{q+1-p}}\qquad\forall x\in \Gw.
\ee
\es


The first application is a pointwise upper bound for solutions with isolated singularities.

\bcor {ptws}Assume $q>p-1>0$, $R^*>0$ and $\Gw$ is a domain containing $0$ such that $d(0)\geq 2R^*$. Then for any $x\in B_{R^*}\setminus\{0\}$, and $0<R\leq R^*$, any  $u\in C^1(\Gw\setminus\{0\})$ solution of $(\ref{A1})$ in $\Gw\setminus\{0\})$ satisfies
\bel{estu0-1}
\abs{u(x)}\leq c_{N,p,q}\abs {|x|^{\frac{q-p}{q+1-p}}-R^{\frac{q-p}{q+1-p}}}+\max\{\abs{u(z)}:\abs z=R\},
\ee
if $p\neq q$, and
\bel{estu0-2}
\abs{u(x)}\leq c_{N,p}\left(\ln R-\ln \abs x\right)+\max\{\abs{u(z)}:\abs z=R\},
\ee
if $p=q$.
\es

The second application corresponds to solutions with boundary blow-up. For $\gd>0$ small enough we set $\Gw_{\gd}:=\{z\in\Gw:d(z)<\gd\}$.

\bcor{unif}Assume $q>p-1>0$, $\Gw$ is a bounded domain with a $C^2$ boundary. Then there exists $\gd_1>0$ which depends only on $\Gw$ such that any   $u\in C^1(\Gw) $ solution of $(\ref{A1})$ in $\Gw$ satisfies
\bel{estu1}
\abs{u(x)}\leq c_{N,p,q}\abs {(d(x))^{\frac{q-p}{q+1-p}}-\gd_1^{{\frac{q-p}{q+1-p}}}}+\max\{\abs{u(z)}:d(z)=\gd_1\}\quad\forall x\in \Gw_{\gd_1}
\ee
 if $p\neq q$, and
\bel{estu2}
\abs{u(x)}\leq c_{N,p,q}\left(\ln \gd_1-\ln d(x)\right)+\max\{\abs{u(z)}:d(z)=\gd_1\}\quad\forall x\in \Gw_{\gd_1}
\ee
if $p=q$.
\es


\noindent\Remark As a consequence of $(\ref{estu1})$ there holds for $p>q>p-1$
\bel{estu3}
u(x)\leq \left(c_{N,p,q}+K\max\{\abs{u(z)}:d(z)\geq\gd_1\}\right)(d(x))^{\frac{q-p}{q+1-p}}\quad\forall x\in \Gw
\ee
where $K=({\rm diam} (\Gw))^{\frac{p-q}{q+1-p}}$, with the standard modification if $p=q$.\medskip

As a variant of \rcor{unif} the following upper estimate of solutions in an exterior domain will be used in the sequel.
\bcor{unif2}Assume $q>p-1>0$, $R>0$ and   $u\in C^1(B_{R_0}^c)$ is any solution of $(\ref{A1})$ in $B_{R_0}^c$. Then for any $R>R_0$ there holds
\bel{estu1-1}
\abs{u(x)}\leq c_{N,p,q}\abs{(\abs x-R_0)^{\frac{q-p}{q+1-p}}-(R-R_0)^{\frac{q-p}{q+1-p}}}+\max\{\abs{u(z)}:\abs z=R\}\quad\forall x\in B_{R}^c
\ee
 if $p\neq q$ and
\bel{estu2-2}
\abs{u(x)}\leq c_{N,p,q}\left(\ln (\abs x-R_0)-\ln (R-R_0)\right)+\max\{\abs{u(z)}:\abs z=R\}\quad\forall x\in B_{R}^c
\ee
if $p=q$.
\es
\noindent\Proof The proof is a consequence of the identity
$$u(x)=u(z)+\myint{0}{1}\myfrac{d}{dt}u(tx+(1-t)z)dt=\myint{0}{1}\langle\nabla u(tx+(1-t)z),x-z\rangle dt
$$
where $z=\frac{R}{\abs x}x$. Since by  $(\ref{est0})$
$$\abs{\nabla u(tx+(1-t)z)}\leq C_{N,p,q}(t\abs{x}+(1-t)R-R_0)^{-\frac{1}{q+1-p}},$$
$(\ref{estu1-1})$ and \eqref{estu2-2} follow by integration. \qeda\medskip



\subsection{Boundary a priori estimates}

The next result is the extension to isolated boundary singularities of a previous regularity estimate  dealing with singularity in a domain proved in \cite[Lemma 3.10]{BVGHV1}.
\blemma{lest2} Assume $p-1<q<p$, $\Gw$ is a bounded $C^2$ domain such that $0\in\prt\Gw$. Let $u\in C^1(\overline\Gw\setminus\{0\})$ be a solution of $(\ref{A1})$ in $\Gw$ which vanishes on $\prt\Gw\setminus\{0\}$ and satisfies
\bel{est2*}
\abs{u(x)}\leq \gf(\abs x)\qquad\forall x\in \Gw,
\ee
where $\gf: \BBR^*_+\mapsto \BBR_+$ is continuous, nonincreasing and satisfies
\bel{est2°}\gf(rs)\leq \gamma\gf(r)\gf(s)\quad\text{and }\;r^{\frac{p-q}{q+1-p}}\gf(r)\leq c,
\ee
 for some $\gg,c>0$ and any $r,s>0$. There exist $\ga\in (0,1)$ and $ c_1=c_1(p,q,\Gw)>0$ such that
\bel{est2**}\BA {lll}
\!(i)&\abs{\nabla u(x)}\leq c_1\gf(\abs x)\abs x^{-1}&\quad\forall x\in \Gw,\\[2mm]
\!(ii)&\abs{\nabla u(x)-\nabla u(y)}\leq c_1\gf(\abs x)\abs x^{-1-\ga}\abs{x-y}^\ga&\quad\forall x,y\in \Gw,\;\abs x\leq\abs y.
\EA\ee
Furthermore
\bel{est3}
 u(x)\leq c_1\gf(\abs x)\frac{d(x)}{\abs x}\qquad\forall x\in \Gw.
\ee
\es
\noindent\Proof For $\ell>0$, we set $\Gw^\ell:=\frac{1}{\ell}\Gw$.
If $\ell\in (0,1]$ the curvature of $\prt\Gw^\ell$ remains uniformly bounded. As in \cite[p 622]{BoVe}, there exists $0<\gd_0\leq 1$ and an involutive diffeomorphism $\psi$ from $\overline B_{\gd_0}\cap\overline \Gw^{\gd_0}$ into $\overline B_{\gd_0}\cap(\Gw^{\gd_0})^c$ which is the identity on $\overline B_{\gd_0}\cap\prt \Gw^{\gd_0}$ and such that $D\psi(\xi)$ is the symmetry with respect to the tangent plane $T_\xi\prt\Gw$ for any
$\xi\in \prt\Gw\cap \overline B_{\gd_0}$. We extend any function $v$ defined in  $\overline B_{\gd_0}\cap\overline \Gw^{\gd_0}$ and vanishing on $\overline B_{\gd_0}\cap\prt \Gw^{\gd_0}$ into a function $\tilde v$ defined in $\overline B_{\gd_0}$
by
\bel{est4-1}\tilde v(x)=\left\{\BA {ll}
v(x)\qquad&\text {if }x\in \overline B_{\gd_0}\cap\overline \Gw^{\gd_0}\\
-v\circ\psi(x)&\text {if }x\in \overline B_{\gd_0}\cap(\Gw^{\gd_0})^c,
\EA\right.\ee
If $v\in C^1(\overline B_{\gd_0}\cap\overline \Gw^{\gd_0})$ is a solution of $(\ref{A1})$ in $B_{\gd_0}\cap\Gw^{\gd_0}$ which vanishes on $\prt\Gw^{\gd_0}\cap \overline B_{\gd_0}$, $\tilde v$ satisfies
\bel{est4-2}
-\sum_j\myfrac{\prt}{\prt x_j}\tilde A_j(x,\nabla \tilde v)+B(x,\nabla \tilde v)=0\quad \text{in }B_{\gd_0}.
\ee
As in \cite[(2.37)]{BoVe} the $A_j$ and $B$ satisfy the following estimates

\bel{est4-3}\left.\BA {llll}
&(i) \quad &\tilde A_j(x,0)=0\\[2mm]

&(ii)\quad &\displaystyle\sum_{i,j}\myfrac{\prt}{\prt \eta_i}\tilde A_j(x,\eta)\xi_i\xi_j\geq C_1\abs\eta^{p-1}\abs\xi^2
\\[4mm]\displaystyle

&(iii) \quad &\displaystyle\sum_{i,j}\abs{\myfrac{\prt}{\prt \eta_j}\tilde A_j(x,\eta)}\leq C_2\abs\eta^{p-2},
\EA\right.\ee
and
\bel{est4-4}
\abs{B(x,\eta)}\leq C_3(1+\abs\eta)^p,
\ee
where the $C_j$ are positive constants. These estimates are the ones needed to apply Tolksdorf's result \cite[Th 1,2]{To1}. There exists a constant $C$, such that for any ball $\overline B_{3R}\subset \overline B_{\gd_0}$, there holds
\bel{est4-5}
\norm{\nabla \tilde v}_{L^\infty(B_R)}\leq C,
\ee
where $C$ depends on the constants $C_k$ ($k=1,2,3$), $N$, $p$ and $\norm{\tilde v}_{L^\infty(B_{3R})}$. We define
\bel{est4-6}
 \Phi_\ell[u](y):=u_\ell=\frac1{\gf(\ell)} u(\ell y)\qquad\forall y\in \Gw^\ell.
\ee
Then
\bel{est4-7}
 |u_\ell(y)|\leq \frac{\gf(\ell\abs y)}{\gf(\ell)}\leq \gg \gf(\abs y)\qquad\forall y\in \Gw^\ell
\ee
and
\bel{est4-8}
-\Gd_pu_\ell +(\ell^{\gb_q}\gf(\ell))^{q+1-p}\abs{\nabla u_\ell}^q=0\qquad\text{in }\Gw^\ell.
\ee
Using formula $(\ref{est4-1})$ we extend $u_\ell$  into a function $\tilde u_\ell$ which satisfies
\bel{est4-9}
-\sum_j\myfrac{\prt}{\prt y_j}\tilde A_j(y,\nabla \tilde u_\ell)+(\ell^{\gb_q}\gf(\ell))^{q+1-p}B(y,\nabla \tilde u_\ell)=0\quad \text{in }B_{\gd_0}.
\ee
For $0<\abs x<\gd_0$ there exists $\ell\in (0,2)$ such that $\frac{\gd_0\ell}{2}\leq\abs x\leq \gd_0\ell$. Then $y\mapsto \tilde u_\ell(y)$ with $y=\frac{x}{\ell}$ satisfies $(\ref{est4-9})$ in $B_{\gd_0}$ and $|\tilde u_\ell(y)|\leq \gamma_*\gf(\abs y)$ since
$\psi$ is a diffeomorphism and $D\psi(\xi)\in O(N)$ for any $\xi\in \prt\Gw\cap B_{\gd_0}$. The function $\tilde u_\ell$ remains bounded on any ball $B_{3R}(z)\subset \Gamma:=\{y\in\BBR^N:\frac{\gd_0}{2}<\abs y<\gd_0\}$, therefore   $\abs{\nabla \tilde u_\ell(y)}\leq c$ for any $y\in B_{R}(z)$, for some constant $c>0$. This implies
\bel{est4-10}\BA{ll}
\abs{\nabla u(x)}\leq c\gg_*\gd_0\gf(\frac{2}{\gd_0})\gf(|x|)|x|^{-1}\quad\forall x\in \Gw\cap B_{\gd_0},
\EA\ee
which is $(\ref{est2**})$-(i). Moreover, by standard regularity estimates \cite{Lieb}, there exists $\ga\in (0,1)$ such that $\abs{\nabla \tilde u_\ell(y)-\nabla \tilde u_\ell(y')}\leq c\abs{y-y'}^\ga$ for all $y$ and $y'$ belonging to  $B_{R}(z)$. This implies $(\ref{est2**})$-(ii).

Next we prove $(\ref{est3})$. Let $0<\gd_1\leq\gd_0$ such that at any boundary point $z$ there exist two closed balls of radius $\gd_1$ tangent to $\prt\Gw$ at $z$ and which are included in $\Gw\cup\{z\}$ and in $\overline\Gw^c\cup\{z\}$ respectively ($\gd_1$ corresponds to the maximal radius of the interior and exterior sphere condition). Let  $x\in\Gw$ such that $d(x)\leq \gd_1$
(this is not a loss of generality) and $z_x$ be the projection of $x$ on $\prt\Gw$.  We first assume that $x$ does not belong to the cone $\Gs_{\frac{\gp}{4}}$ with vertex $0$, axis $-{\bf n}_0$, where ${\bf n}_0$ is the normal outward unit vector at $0$, and angle $\frac{\gp}{4}$. Consider the path $\zeta$ from $z_x$ to $x$ defined by
$\zeta(t)=tx+(1-t)z_x$ with $0\leq t\leq 1$. Then
\bel{est4-10'}\BA{ll}u(x)=\myint{0}{1}\myfrac{d}{dt}u\circ\zeta(t)dt=\myint{0}{1}\langle\nabla u\circ\zeta(t),x-z_x\rangle dt
\EA\ee
Thus, by the Cauchy-Schwarz inequality, using $(\ref{est2**})$,
\bel{est4-10''}\BA{ll}
\abs{u(x)}\leq c_1d(x)\myint{0}{1}\myfrac{\phi(|\zeta(t)|)}{\abs{\zeta(t)}}dt.
\EA\ee
Since $x\notin\Gs_{\frac{\gp}{4}}$, $\zeta(t)\notin\Gs_{\frac{\gp}{4}}$ and there exists $c_2>0$ depending on $\Gw$ such that
$c^{-1}_2\abs x\leq \abs{\zeta(t)}\leq c_2\abs x$ for all $0\leq t\leq 1$. Therefore $\phi(|\zeta(t)|)\leq \phi(c_2\abs x)\leq \gamma\phi(c_2)\phi(\abs x)$ by $(\ref{est2°})$. This implies
\bel{est4-10'''}\BA{ll}
\abs{u(x)}\leq \gamma c_1c_2\phi(c_2)\myfrac{d(x)\phi(\abs x)}{\abs{x}}
\EA\ee
by $(\ref{est3})$ whenever $x\notin \Gs_{\frac{\gp}{4}}$. When $x\in \Gs_{\frac{\gp}{4}}$ then $d(x)\leq \abs x\leq c_3d(x)$
where $c_3>0$ depends on the curvature of $\prt\Gw$. Then $(\ref{est2*})$ combined with $(\ref{est2°})$ implies the claim. \qeda


\blemma{lest1} Assume $p-1<q\leq p$, $\Gw$ is a bounded $C^2$ domain such that $0\in\prt\Gw$ and $R_0=\max\{ \abs z:z\in\Gw\}$. If $u\in C(\overline\Gw\setminus\{0\})\cap C^1(\Gw)$ is a positive solution of
$(\ref{A1})$ which vanishes on $\prt\Gw\setminus\{0\}$, it satisfies
\bel{est1}
u(x)\leq\left\{\BA {ll} c_2\left(\abs x^{\frac{q-p}{q+1-p}}-R_0^{\frac{q-p}{q+1-p}}\right)\qquad&\text{if }q< p\\[2mm]
 (p-1)\ln\left(\frac{R_0}{\abs x}\right)\qquad&\text{if }q= p
 \EA\right.
\ee
for all $ x\in \Gw$, where $c_2=c_2(p,q)>0$.
\es
\noindent\Proof For $\ge>0$ we denote by $P_\ge:\BBR\mapsto\BBR_+$ the function defined  by
\bel{est2-0}
P_\ge(r)=\left\{\BA {ll}0\qquad&\text{if }0\leq r\leq \ge\\
-\frac{r^4}{2\ge^3}+\frac{3r^3}{\ge^2}-\frac{6r^2}{\ge}+5r-\frac{3\ge}{2}\qquad&\text{if }\ge<r<2 \ge\\
r-\frac{3\ge}{2}\qquad&\text{if }r\geq 2\ge,
\EA\right.\ee
and by $u_\ge$ the extension of $P_\ge(u)$ by zero outside $\Gw$. There exists $R_0$ such that $\Gw \subset B_{R_0}$ .
Since $0\leq P_\ge(r)\leq |r|$ and $P_\ge$  is convex, $u_\ge \in C(\BBR^N\setminus\{0\})\cap W^{1,p}_{loc}(\BBR^N\setminus\{0\})$ and
$$-\Gd_pu_\ge+\abs{\nabla u_\ge}^q\leq 0\qquad\text{in }\BBR^N.
$$
Let $R>R_0$. If $p-1<q<p$
\bel{est2-1}U_{\ge,R}(\abs x)=c_2 \left((|x|-\ge)^{\frac{q-p}{q+1-p}}-(R-\ge)^{\frac{q-p}{q+1-p}}\right)\qquad\text{in }B_R\setminus B_\ge,
\ee
with  $c_2=(p-q)^{-1}(q+p-1)^{\frac{q-p}{q+1-p}}$. Then $-\Gd_pU_{\ge,R}+\abs{\nabla U_{\ge,R}}^q\geq 0$. Since $u_\ge$ vanishes on $\prt B_R$ and is finite on $\prt B_\ge$, it follows $u_\ge\leq U_{\ge,R}$. Letting successively $\ge\to  0$ and $R\to R_0$ yields to
$(\ref{est1})$. If $q=p$ we take
\bel{est2-1-1}U_{\ge,R}(\abs x)=(p-1) \ln\left(\frac{R-\ge}{\abs x-\ge}\right)\qquad\text{in }B_R\setminus B_\ge,
\ee
which turns out to be a supersolution of $(\ref{A1})$; the end of the proof is similar.
\qeda\medskip

As a consequence of \rlemma{lest2} and \rlemma{lest1}, we obtain.

\bcor{lest3} Let $p,q$ $\Gw$ and $u$ be as in \rlemma{lest1}. Then there exists a constant $c_3=c_3(p,q,\Gw)>0$ such that
\bel{est4-11}
\abs{\nabla u(x)}\leq c_3\abs{x}^{-\frac{1}{q+1-p}}\qquad\forall x\in\Gw
\ee
and
\bel{est4-12}
u(x)\leq c_3d(x)\abs{x}^{-\frac{1}{q+1-p}}\qquad\forall x\in\overline\Gw\setminus\{0\}.
\ee
\es

\noindent\Remark  If $\Gw$ is locally flat near $0$, then estimates $(\ref{est4-11})$ and $(\ref{est4-12})$ are valid without any sign assumption on $u$. More precisely, if $\prt\Gw\cap B_{\gd_0}=T_0\prt\Gw\cap B_{\gd_0}$ we can perform the reflection of
$u$ through the tangent plane $T_0\prt\Gw$ to $\prt\Gw$ at $0$ and the new function $\tilde u$ is a solution of $(\ref{A1})$ in $B_{\gd_0}\setminus\{0\}$. By \rprop{lets0}, it satisfies
\bel{est4-12-1}
\abs{\nabla \tilde u(x)}\leq c_{N,p,q}\abs{x}^{-\frac{1}{q+1-p}}\qquad\forall x\in B_{\frac{\gd_0}{2}}\setminus\{0\}.
\ee
Integrating this relation as in \cite{BVGHV1}, we derive that for any $x\in B_{\frac{\gd_0}{2}}\cap\Gw$, there holds
\bel{est4-12-2}
\abs{u(x)}\leq \left\{\BA {ll}
c_{N,p,q}\left(|x|^{-\gb_q}-(\frac{\gd_0}{2})^{-\gb_q}\right)+\max\{\abs{u(z)}:\abs z=\frac{\gd_0}{2}\}\quad&\text{if } p\neq q\\
c_{N,p}\ln \left(\frac{\gd_0}{2\abs x}\right)+\max\{\abs{u(z)}:\abs z=\frac{\gd_0}{2}\}\quad&\text{if } p= q.
\EA\right.\ee\medskip

In the next result we allow the boundary singular set to be a compact set.

\bprop{lest4} Let $p-1<q<p$ and $\gd_1$ as above. There exist $ r^*\in (0, \gd_1]$ and $c_4=c_4(N,p,q)>0$ such that for any nonempty compact set $K\subset\prt\Gw$, $K\neq\prt\Gw$ and any positive solution $u\in C(\overline\Gw\setminus K)\cap C^1(\Gw)$ of $(\ref{A1})$ which vanishes on $\prt\Gw\setminus K$, there holds
\bel{est6}
u(x)\leq c_4d(x)(d_K(x))^{-\frac{1}{q+1-p}}\qquad\forall x\in\prt\Gw\text{ s.t. }d(x)\leq r^*,
\ee
where $d_K(x)=\dist (x,K)$.
\es
\noindent\Proof {\it Step 1: Tangential estimates}. Let $x\in\Gw$ such that $d(x)\leq \gd_1$. We denote by $\gs(x)$ the projection of $x$ onto $\prt\Gw$, unique since $d(x)\leq \gd_1$. Let $r\,,r',\gt>0$ such that $\frac{3}{4}r<r'<\frac{7}{8}r$ and $0<\gt\leq \frac{r'}{2}$ and put $\gw_{\gt,x}=\gs(x)+\gt{\bf n}_{\gs(x)}$. Since $\prt\Gw$ is $C^2$, there exists $0< r^*\leq \gd_1$ depending on $\Gw$ such that $d_K(\gw_{\gt,x})>\frac{7}{8}r$ whenever $d(x)\leq r^*$. Let $a>0$ and $b>0$ to be specified later on; we define $\tilde v(s)=a(r'-s)^{\frac{q-p}{q+1-p}}-b$ and $v(y)=\tilde v(\abs{y-\gw_{\gt,x}})$ in $[0,r')$ and $B_{r'}(\gw_{\gt,x})$ respectively. Then
$$\BA {l}\abs {\tilde v'}^{p-2}\left(\abs{\tilde v'}^{q+2-p}-(p-1)\tilde v''-\myfrac{N-1}{s}\tilde v'\right)
\!=a^{p-1}\left(\myfrac{p-q}{q+1-p}\right)^{p-1}\!\!\!\!(r'-s)^{-\frac{q}{q+1-p}}X(s)
\EA$$
where
$$X(s)=\left(a\myfrac{p-q}{q+1-p}\right)^{q+1-p}-\myfrac{p-1}{q+1-p}-\myfrac{(N-1)(r'-s)}{s}.
$$
For any $\gt\in (0,r')$ there exists $a>0$ such that
$$\left(a\myfrac{p-q}{q+1-p}\right)^{q+1-p}\geq \myfrac{p-1}{q+1-p}+\myfrac{(N-1)(r'-s)}{s}\qquad\forall \gt\leq s\leq r'.
$$
This implies
\bel{est7}
-\Gd_p v+\abs{\nabla v}^q\geq 0\qquad\text{in }B_{r'}(\gw_{\gt,x})\setminus B_{\gt}(\gw_{\gt,x}).
\ee
Next we take $b=a(r'-\gt)^{\frac{q-p}{q+1-p}}$, thus $v=0$ on $\prt B_{\gt}(\gw_{\gt,x})$. Clearly $B_{\gt}(\gw_{\gt,x})\subset\overline\Gw^c$ since $\gt<\gd_1$. Therefore $v\geq 0=u$ on $\prt\Gw\cap B_{r'}(\gw_{\gt,x})$  and
$u\leq v=\infty$ on $\Gw\cap\prt B_{r'}(\gw_{\gt,x})$. By the comparison principle, $v\geq u$ in $\Gw\cap B_{r'}(\gw_{\gt,x}).$ In particular
$$u(x)\leq v(x)\leq a(r'-\gt-d(x))^{\frac{q-p}{q+1-p}}-a(r'-\gt)^{\frac{q-p}{q+1-p}}.
$$
We take now $\gt=\frac{r'}{2}$ and $d(x)\leq \frac{r}{4}$ and we derive by the mean value theorem
\bel{est8}
u(x)\leq c'_4r'^{-\frac{1}{q+1-p}}d(x)=c'_4d(x)(d_K(x))^{-\frac{1}{q+1-p}},
\ee
with $c'_4=c'_4(p,q)>0$ Letting $r'\to \frac{7}{8}r$, we get $(\ref{est3})$.\smallskip

\noindent{\it Step 2: Global estimates}.  If $d(x)\geq \frac{1}{4}d_K(x)$, there holds
$$d(x)(d_K(x))^{-\frac{1}{q+1-p}}\geq 2^{-\frac{2}{q+1-p}}(d(x))^{\frac{q-p}{q+1-p}}.$$
Combining this inequality with $(\ref{estu3})$ and obtain $(\ref{est6})$. \qeda\medskip

\noindent\Remark  Under the assumption of \rprop{lest4}, it follows from the maximum principle that $u$ is upper bounded  in the set $\Gw'_{r^*}:=\{x\in\Gw:d(x)>r^*\}=\Gw\setminus\overline\Gw_{r^*}$ by the solution $w$ of
\bel{est9}\BA {ll}
-\Gd_pw+|\nabla w|^q=0\qquad&\text {in }\Gw_{r^*}\\\phantom{-\Gd_p+|\nabla v|^q}
w=c_4d(x)(d_K(x))^{-\frac{1}{q+1-p}}\qquad&\text {in }\prt\Gw_{r^*},
\EA\ee
and $w$ itself is bounded by $d^*=\max\{cd(x)(d_K(x))^{-\frac{1}{q+1-p}}:d(x)=r^*\}$.\medskip


Next we prove a boundary Harnack inequality. We recall that $\gd_1$ has been introduced at \rcor{unif}, and that the interior and exterior sphere conditions hold in the set $\{x\in\BBR^N:d (x)\leq\gd_1\}$.

\bth{bharnack} Let $q>p-1$ and $0\in\prt\Gw$. Then there exists $c_5=c_5(N,p,q,\Gw)>0$ such that for any
positive solution $u\in C(\Gw\cup((\prt\Gw\setminus\{0\})\cap B_{2\gd_1})\cap C^1(\Gw)$ of $(\ref{A1})$ in $\Gw$, vanishing on $\prt\Gw\setminus\{0\})\cap B_{2\gd_1}$, there holds
\bel{harn3}
\myfrac{u(y)}{c_5d(y)}\leq \myfrac{u(x)}{d(x)}\leq c_5\myfrac{u(y)}{d(y)}
\ee
for all $x,y\in B_{\frac{2\gd_1}{3}}\cap\Gw$ such that $\frac{1}{2}\abs x\leq \abs y\leq 2\abs x$.
\es

For proving \rth{bharnack} we need some intermediate lemmas. First we recall the following result from \cite {Bau}.

\blemma{chain}Assume that $a \in \prt \Gw$, $0<r<\gd_1$ and $h>1$ is an integer. There exists an integer $N_0$, depending only on $\gd_1$, such that for any points $x$ and $y$ in $\Gw \cap B_{\frac{3r}{2}}(a)$ verifying $\min\{d(x),d(y)\} \geq r/2^h$, there exists a connected chain of balls $B_1,...,B_j$ with $j\leq N_0h$ such that
	\bel{geo} \BA{ll} x \in B_1, y \in B_j, \q B_i\cap B_{i+1} \neq \ems  \text{ for } 1\leq i \leq j-1 \\
        \text{and } 2B_i \sbs B_{2r}(Q) \cap \Gw \text{ for } 1\leq i \leq j.
	\EA \ee
\es
The next result is a standard Harnack inequality.
\blemma{Har1} Assume $a \in (\prt \Gw\sms\{0\}) \cap B_{\frac{2\gd_1}{3}}$ and $0<r\leq \abs{a}/4$. Let $u\in C(\Gw\cup ((\prt\Gw\sms\{0\})\cap B_{2\gd_1}))\cap C^1(\Gw)$ be a positive solution of $(\ref{A1})$ vanishing on $(\prt\Gw\sms\{0\})\cap B_{2\gd_1}$. Then there exists a positive constant $c_6>1$ depending on $N,$ $p,$ $q$ and $\gd_1$ such that
 \bel{Har1-1} u(x) \leq c_6^h u(y),\ee
for every $x,y \in B_{\frac{3r}{2}}(a)\cap \Gw$ such that $\min\{d(x),d(y)\} \geq r/2^h$ for some $h \in \BBN$.
\es
\noindent\Proof For $\ell>0$, we define $T_{\ell}[u]$ by
\bel{scal2}
T_\ell[u](x)=\ell^{\frac{p-q}{q+1-p}}u(\ell x),
\ee
and we notice that if $u$ satisfies $(\ref{A1})$ in $\Gw$, then $T_\ell[u]$ satisfies the same equation in $\Gw^{\ell}:=\ell^{-1}\Gw$. If we take in particular $\ell=|a|$, we can assume $|a|=1$, thus the curvature of the domain $\Gw^{|a|}$ remains bounded. By \rprop{lest4}
	\bel{Har1-2} u(x) \leq c'_6\forevery x \in B_{2r}(a)\cap \Gw \ee
where $c'_6$ depends on $N$, $q$, $\gd_1$. Then we proceed as in \cite {NpV}, using \rlemma{chain} and internal Harnack inequality as quoted in \cite[Corollary 10]{Tru}. \qeda\medskip

Since the solutions are H\"older continuous, the following statement holds as in \cite[Theorem 4.2] {Tru}:

\blemma{Har2} Let the assumptions on $a$ and $u$ of \rlemma{Har1} be fulfilled. If $b \in \prt \Gw \cap B_r(a)$ and $0<s\leq 2^{-1}r$, there exist two positive constants $\gd$ and $c_7$ depending on $N$, $p$, $q$ and $\Gw$ such that
	\bel{Har2-1} u(x) \leq c_7\myfrac{\abs{x-b}^\gd}{s^\gd}\max\{u(z): z\in B_r(b)\cap \Gw\}\ee
for every $x \in B_s(b) \cap \Gw$.
\es
As a consequence we derive the following Carleson type estimate.

\blemma{Har3} Assume $a \in (\prt \Gw\sms\{0\}) \cap B_{\frac{2\gd_1}{3}}$ and $0<r\leq \abs{a}/8$. Let $u\in C(\Gw\cup ((\prt\Gw\sms\{0\})\cap B_{2\gd_1}))\cap C^2(\Gw)$ be a positive solution of $(\ref{A1})$ vanishing on $(\prt\Gw\sms\{0\})\cap B_{2\gd_1}$. Then there exists a constant $c_8$ depending only on $N$, $p$ and $q$ such that
	\bel{Har3-1}\BA {ll} u(x) \leq c_8 u(a-\frac{r}{2}{\bf n}_{_a}) \forevery x \in B_r(a) \cap \Gw.
	\EA\ee
\es
\noindent\Proof By \rlemma{Har1} it is clear that for any integer $h$ and $x \in B_r(a) \cap \Gw$ such that $d(x)\geq 2^{-h}r$, there holds
\bel{Har3-2} \BA {ll}
u(x) \leq c_6^h u(a-\frac{r}{2}{\bf n}_{_a}).
\EA\ee
Therefore $u$ satisfies inequality $(\ref{Har2-1})$ as any H\"older continuous function does. The proof that the constant is independent of $r$ and $u$ is more delicate. It is done in  \cite[Lemma 2.4]{Bau} for linear equations, but it is based only on \rlemma{Har2} and a geometric construction, thus it is also valid in our case. \qeda

\blemma{LUE1} Assume $a \in (\prt \Gw\sms\{0\}) \cap B_{\frac{2\gd_1}{3}}$ and $0<r\leq \abs{a}/8$. Let $u\in C(\Gw\cup ((\prt\Gw\sms\{0\})\cap B_{2\gd_1}))\cap C^2(\Gw)$ be a positive solution of $(\ref{A1})$ vanishing on $(\prt\Gw\sms\{0\})\cap B_{2\gd_1}$. Then there exist $\ga \in (0,1/2)$ and $c_9>0$ depending on $N$, $p$ and $q$ such that
	\bel{LUE1} \BA {ll}
	\myfrac{1}{c_9}\myfrac{t}{r} \leq \myfrac{u(b-t{\bf n}_{_b})}{u(a-\frac{r}{2}{\bf n}_{_a})} \leq c_9\myfrac{t}{r}\EA
	\ee
for any $b \in B_r(a) \cap \prt \Gw$ and $0 \leq t < \frac{\ga}{2}r$.
\es
\noindent\Proof It is similar to the one of \cite[Lemma 3.15]{NpV}. \qeda\medskip

\noindent {\it Proof of \rth{bharnack}}. Assume $x \in B_{\frac{2\gd_1}{3}} \cap \Gw$ and set $r=\frac{\abs{x}}{8}$. \medskip

\noindent{\it Step 1: Tangential estimate: we  suppose  $d(x) < \frac{\ga}{2}r$}. Let $a \in \prt \Gw \sms \{0\}$ such that $\abs{a}=\abs{x}$ and $x \in B_r(a)$. By \rlemma{LUE1},
	\bel{Hark1-2} \myfrac{8}{c_9}\myfrac{u(a-\frac{r}{2}{\bf n}_{_a})}{\abs{x}} \leq \myfrac{u(x)}{d(x)} \leq 8c_9\myfrac{u(a-\frac{r}{2}{\bf n}_{_a})}{\abs{x}}. \ee
We can connect $a-\frac{r}{2}{\bf n}_{_a}$ with $-2r{\bf n}_{_0}$ by $m_1$ (depending only on $N$) connected balls $B_i=B_\frac{r}{4}(x_i)$ with $x_i \in \Gw$ and $d(x_i) \geq \frac{r}{2}$ for every $1 \leq i \leq m_1$. It follows from $(\ref{Har3-1})$ that
	$$\BA {ll} c^{-m_1}_6u(-2r{\bf n}_{_0}) \leq u(a-\frac{r}{2}{\bf n}_{_a}) \leq c^{m_1}_6u(-2r{\bf n}_{_0}),
	\EA$$
which, together with $(\ref{Hark1-2})$ leads to
	\bel{Hark1-3} \myfrac{1}{c_{10}}\myfrac{u(-2r{\bf n}_{_0})}{\abs{x}} \leq \myfrac{u(x)}{d(x)} \leq c_{10}\myfrac{u(-2r{\bf n}_{_0})}{\abs{x}}, \ee
	with $c_{10}=8c_9c^{m_1}_6$.\smallskip
	
\noindent{\it Step 2: Internal estimate: we suppose $d(x) \geq \frac{\ga}{2}r$}. We can connect $-2r{\bf n}_{_0}$ with $x$ by $m_2$ (depending only on $N$) connected balls $B'_i=B_\frac{\ga r}{4}(x'_i)$ with $x'_i \in \Gw$ and $d(x'_i) \geq \frac{\ga}{2}r$ for every $1 \leq i \leq m_2$. By Harnack and Carleson inequalities $(\ref{Har1-1})$ and $(\ref{Har3-1})$ and since $\frac{\ga}{4}\abs{x}<d(x)\leq \abs{x}$, we get
	\bel{Hark1-4}  \myfrac{\ga}{4c_{6}'^{m_2}}\myfrac{u(-2r{\bf n}_{_0})}{\abs{x}} \leq \myfrac{u(x)}{d(x)} \leq \myfrac{4c_{6}'^{m_2}}{\ga}\myfrac{u(-2r{\bf n}_{_0})}{\abs{x}}. \ee
\noindent{\it Step 3: End of proof.} Suppose $\frac{\abs{x}}{2}\leq s \leq 2\abs{x}$, we can connect $-2r{\bf n}_{_Q}$ with $-s{\bf n}_{_Q}$ by $m_3$ (depending only on $N$) connected balls $B''_i=B_\frac{r}{2}(x''_i)$ with $x''_i \in \Gw$ and $d(x''_i) \geq r$ for every $1 \leq i \leq m_3$. This fact, jointly with $(\ref{Hark1-3})$ and $(\ref{Hark1-4})$, yields to
	\bel{Hark1-5}  \myfrac{1}{c_{11}}\myfrac{u(-s{\bf n}_{_0})}{\abs{x}} \leq \myfrac{u(x)}{d(x)} \leq c_{11}\myfrac{u(-s{\bf n}_{_0})}{\abs{x}} \ee
where $c_{11}=c_{11}(N,q,\Gw)$. Finally, if $y \in B_{\frac{2r_0}{3}} \cap \Gw$ satisfies $\frac{\abs{x}}{2} \leq \abs{y} \leq 2\abs{x}$, then by applying twice $(\ref{Hark1-5})$ we get $(\ref{harn3})$ with $c_{5}=c_{11}^2$. \qeda\medskip

The following inequality is a consequence of \rth{bharnack}.
\bcor{BdryHarn}
Assume $q>p-1$ and $0\in\prt\Gw$. Then there exists $c_{12}>0$ depending on $p$, $q$ and $\Gw$ such that for any
positive solutions $u_1,\,u_2\in C(\Gw\cup((\prt\Gw\setminus\{0\})\cap B_{2\gd_1}))\cap C^1(\Gw)$ of $(\ref{A1})$ in $\Gw$, vanishing on $(\prt\Gw\setminus\{0\})\cap B_{2\gd_1}$, there holds
\bel{harn3'}
\sup\left\{\myfrac{u_1(y)}{u_2(y)}:y\in B_r\setminus B_{\frac{r}{2}}\right\}\leq c_{12}\inf\left\{\myfrac{u_1(y)}{u_2(y)}:y\in B_r\setminus B_{\frac{r}{2}}\right\}.
\ee
\es


\mysection{Boundary singularities}
\subsection{Strongly singular  solutions}
In this section we consider the equation $(\ref{A1})$ in $\BBR_+^N$. We denote by $(r,\gs) \in \BBR_+ \ti S^{N-1}$ the spherical coordinates in $\BBR^N $ and
$$S_+^{N-1}=\left\{(\sin\gf\gs',\cos\gf):\gs'\in S^{N-2},\gf\in [0,\frac{\gp}{2})\right\}.
$$
If $v(x)=r^{-\gb}\gw(\gs)$ satisfies $(\ref{A1})$ in $\BBR_+^N$ and vanishes on $\prt{\BBR^N_+}\setminus\{0\}$, then $\gb=\gb_q$ and $\gw$ is a solution of
\begin{equation}\label{AS1}\left.
\BA {ll}
-div'\left(\left(\gb_q^2\gw^2+|\nabla' \gw|^2\right)^{\frac{p-2}2}\nabla' \gw\right)-\gb_q\Gl_{\gb_q}\left(\gb_q^2\gw^2+|\nabla' \gw|^2\right)^{\frac{p-2}2}\gw
\\[2mm]
\phantom{----------------}
+\left(\gb_q^2\gw^2+|\nabla' \gw|^2\right)^{\frac{q}2}=0\quad \mbox{in } S_+^{N-1}
\\[2mm]
\phantom{------------------------,}
\gw=0\quad \mbox{on } \prt S_+^{N-1}.
\EA\right.\end{equation}
where $\gb_q$ and $\Gl_{\gb_q}$ have been defined in $(\ref{I-8*})$. We denote by $(\gb_*,\psi_*)\in\BBR_+^*\ti C^2(\overline S_+^{N-1})$  the unique couple such $\max\psi_*=1$ with the property that the function
$(r,\gs)\mapsto r^{-\gb_*}\psi_*(\gs)$
is positive, $p$-harmonic in $\BBR_+^N$ and vanishes on $\partial{\BBR^N_+}\setminus\{0\}$. Then
 $\psi_*=\psi$ satisfies
 \begin{equation}\label{AS2}\left.
\BA {ll}
-div'\left(\left(\gb_*^2\psi^2+|\nabla' \psi|^2\right)^{\frac{p-2}2}\nabla' \psi\right)-\gb_*\Gl_{\gb_*}\left(\gb_*^2\psi^2+|\nabla' \psi|^2\right)^{\frac{p-2}2}\psi
=0\quad \mbox{in } S_+^{N-1}
\\[2mm]
\phantom{-div'\left(\left(\gb_*^2\psi^2+|\nabla' \psi|^2\right)^{\frac{p-2}2}\nabla' \psi\right)-\gb_*\Gl_{\gb_*}\left(\gb_*^2\psi^2+|\nabla' \psi|^2\right)^{\frac{p-2}2}}
\psi=0\quad \mbox{on } \prt S_+^{N-1}.
\EA\right.\end{equation}
Since the function $\psi_*$ is unique it depends only on the azimuthal variable $\gth_{N-1}=\cos^{-1}(\frac{x_N}{|x|})$ (see Appendix II). Our first result is the following
\bth{nonexist} If $q\geq q_*$, or equivalently $\gb_q\leq \gb_*$, there exists no positive solution to problem $(\ref{AS1})$.
\es
\noindent\Proof Suppose such a solution $\gw$ exists and put $\gth=\gb_q/\gb_*$, then $0<\gth\leq 1$. Set $\eta=\psi^\gth$, where $\psi$ is a positive solution of $(\ref{AS2})$, and define the operator $\CT$ by
\be\label{AS3}\BA {l}
\CT(\eta)=-div'\left(\left(\gb_q^2\eta^2+|\nabla' \eta|^2\right)^{\frac{p-2}2}\nabla' \eta\right)
-\gb_q\Gl_{\gb_q}\left(\gb_q^2\eta^2+|\nabla' \eta|^2\right)^{\frac{p-2}2}\eta
\\\phantom{-----------------------}+\left(\gb_q^2\eta^2+|\nabla' \eta|^2\right)^{\frac{q}2}.
\EA\ee
Since $\nabla\eta=\gth\psi^{\gth-1}\nabla\psi$,
$$\left(\gb_q^2\eta^2+|\nabla' \eta|^2\right)^{\frac{p-2}2}=\gth^{p-2}\psi^{(\gth -1)(p-2)}\left(\gb_*^2\psi^2+|\nabla' \psi|^2\right)^{\frac{p-2}2},$$
$$\left(\gb_q^2\eta^2+|\nabla' \eta|^2\right)^{\frac{p-2}2}\nabla' \eta=\gth^{p-1}\psi^{(\gth -1)(p-1)}\left(\gb_*^2\psi^2+|\nabla' \psi|^2\right)^{\frac{p-2}2}\nabla'\psi,$$
therefore
$$\BA {l}
\CT(\eta)=-\gth^{p-1}\psi^{(\gth -1)(p-1)}{\rm div'}\left(\left(\gb_*^2\psi^2+|\nabla' \psi|^2\right)^{\frac{p-2}2}\nabla'\psi\right)\\[2mm]
\phantom{\CT(\eta)}-
\gth^{p-1}(\gth -1)(p-1)\psi^{(\gth -1)(p-1)-1}\left(\gb_*^2\psi^2+|\nabla' \psi|^2\right)^{\frac{p-2}2}|\nabla'\psi|^2\\[2mm]
\phantom{\CT(\eta)}
-\gb_q\Gl_{\gb_q}\gth^{p-2}\psi^{(\gth -1)(p-1)}\left(\gb_*^2\psi^2+|\nabla' \psi|^2\right)^{\frac{p-2}2}\psi+\gth^{q}\psi^{(\gth -1)q}\left(\gb_*^2\psi^2+|\nabla' \psi|^2\right)^{\frac{q}2}.
\EA$$
But $\gb_q\Gl_{\gb_q}\gth^{p-2}=\gb_*\Gl_{\gb_q}\gth^{p-1}\leq \gb_*\Gl_{\gb_*}\gth^{p-1}$ since $\gb_q\leq \gb_*$. Using $(\ref{AS2})$, we see that $\CT(\eta)\geq 0$. Because Hopf Lemma is valid, there holds $\prt_{\bf n}\psi <0$ on $\prt S^{N-1}_+$. Since $\gw$ is $C^1$ in  $\overline{S^{N-1}_+}$ and $\psi$ is defined up to an homothety, there exists a smallest function $\psi$ such that $\eta\geq\gw$, and the graphs of $\eta$ and $\gw$ over $\overline{S^{N-1}_+}$
are tangent, either at some $\ga\in S^{N-1}_+$, or only at a point $\ga\in\prt S^{N-1}_+$. We put $w=\eta-\gw$. Then
\be\label{AS4}
\CT(\eta)=\CT(\eta)-\CT(\gw)=\Gf(1)-\Gf(0),
\ee
where $\Phi(t)=\CT(\gw_t)$ with $\gw_t=\gw+tw$.

We use local coordinates $(\gs_1,...,\gs_{N-1})$ on $S^{N-1}$ near $\alpha$. We denote by $g=(g_{ij})$  the
metric tensor on $S^{N-1}$ and by $g^{jk}$ its contravariant components. Then, for any $\varphi\in C^1(S^{N-1})$,
$${\abs {\nabla\varphi }^{2}}=
\sum_{j,k}g^{jk}\frac {\partial \varphi}{\partial \sigma_{j}}\frac {\partial
\varphi}{\partial \sigma_{k}}=\langle\nabla\varphi,\nabla\varphi\rangle_g.
$$
If $X=(X^1,...X^d)\in C^1(TS^{N-1})$ is a vector field, we lower indices by setting
$\displaystyle {X^\ell=\sum_{i}g^{\ell i}X_{i}}$ and define the divergence of $X$ by
$$div'_gX=\frac {1}{\sqrt {\abs g}} \sum_{\ell}
\frac {\partial}{\partial \sigma_{\ell}}\left(\sqrt {\abs g}X^\ell\right)
=\frac {1}{\sqrt {\abs g}} \sum_{\ell,i}\frac {\partial}{\partial \sigma_{\ell}}\left(\sqrt {\abs g}g^{\ell i}X_{i}\right).$$
We write $\Phi(t)=\Phi_1(t)+\Phi_2(t)+\Phi_3(t)$
where
$$\Phi_1(t)=-\gb_q\Gl_{\gb_q}\left(\gb_q^2\gw_t^2+|\nabla' \gw_t|^2\right)^{\frac{p-2}2}\gw_t
\;,\quad\Phi_2(t)=\left(\gb_q^2\gw_t^2+|\nabla' \gw_t|^2\right)^{\frac{q}2}
$$
and
$$\Phi_3(t)=-\rm{div'}\left(\left(\gb_q^2\gw_t^2+|\nabla' \gw_t|^2\right)^{\frac{p-2}2}\nabla' \gw_t\right).
$$
Then

$$\Phi_1(1)-\Phi_1(0)=-\sum_{j}a_{j}\frac {\partial w}{\partial
\sigma_{j}}-bw\,\text{ and }\;\Phi_2(1)-\Phi_2(0)=\sum_{j}c_{j}\frac {\partial w}{\partial
\sigma_{j}}+dw,
$$
where
$$b=\gb_q\Gl_{\gb_q}\left(\beta_q^{2}{\gw_t}^2+\abs {\nabla \gw_t}^2\right)^{\frac{p}{2}-2}
\left((p-1)\gb_q^2\gw_t^2+\abs{\nabla \gw_t}^2\right),
$$

$$a_{j}=(p-2)\gb_q\Gl_{\gb_q}\left(\beta_q^{2}{\gw_t}^2+\abs
{\nabla \gw_t}^2\right)^{\frac{p}{2}-2}\gw_t\sum_{k}g^{jk} \myfrac{\prt\gw_t}{\prt\gs_k},$$

$$d=q\beta_q^{2}\left(\beta^{2}{\gw_t}^2+\abs {\nabla \gw_t}^2\right)^{\frac{q}{2}-1}\gw_t,
$$
and
$$c_{j}=q\left(\beta_q^{2}{\gw_t}^2+\abs
{\nabla \gw_t}^2\right)^{\frac{q}{2}-1}\sum_{k}g^{jk} \myfrac{\prt\gw_t}{\prt\gs_k}.$$
Furthermore
$$\BA {ll}\Phi_3(1)-\Phi_3(0)
=-(p-2){\rm div'}\left(\left(\gb_q^2\gw_t^2+|\nabla' \gw_t|^2\right)^{\frac{p-4}2}\left(\gb_q^2\gw_t w+
\langle\nabla'\gw_t,\nabla' w\rangle_g\right)\nabla'\gw_t\right)
\\\phantom{\Phi_3(1)-\Phi_3(0)}
-{\rm div'}\left(\left(\gb_q^2\gw_t^2+|\nabla' \gw_t|^2\right)^{\frac{p-2}2}\nabla' w\right).
\EA$$
Therefore we can write $\Gf(1)-\Gf(0)$ under the form
\be\label{AS5}\Phi(1)-\Phi(0)=-{\rm div'}(A\nabla'w)+\langle B,\nabla' w\rangle_g+Cw:=\CL w
\ee
where
\be\label{AS6}\BA {ll}\langle AX,X\rangle_g=\left(\gb_q^2\gw_t^2+|\nabla' \gw_t|^2\right)^{\frac{p-4}2}
\left(p-2)\langle \nabla'\gw_t,X\rangle_g^2+|\nabla' \gw_t|^2|X|^2\right)\\\phantom{\langle AX,X\rangle_g}
\geq \left(\gb_q^2\gw_t^2+|\nabla' \gw_t|^2\right)^{\frac{p-4}2}\min\{1,p-1\}|\nabla' \gw_t|^2|X|^2.
\EA\ee
and $B$ and $C$ can be computed from the previous expressions. It is important to notice that $\gb_q^2\gw_t^2+|\nabla' \gw_t|^2$ is bounded between two positive constants $m_1$ and $m_2$ in $\overline{S^{N-1}_+}$. Thus the operator $\CL$ is uniformly elliptic with bounded coefficients. Since $w$ is nonnegative and either  at some point $\alpha$,  $\nabla 'w(\alpha)=0$ and $w(\alpha)>0$, or at some boundary point $\ga$ where $w(\ga)=0$ and $\prt_{\bf n}w(\ga)<0$, it follows from the strong maximum principle or Hopf Lemma (see
\cite{GT}) that $w=0$, contradiction.\qeda

\bth{exist} Assume $q<q_*$ or equivalently  $\gb_q> \gb_*$. There exists a unique positive solution $\gw_*$ to problem $(\ref{AS1})$.
\es
\noindent\Proof {\it Existence}.  It will follow from \cite{BoMuPu}. Indeed problem $(\ref{AS1})$ can be written under the form
\be\label{AS3-0}\BA {lll}
{\bf A}(\gw):=-div'\,{\bf a}(\gw,\nabla'\gw)={\bf B}(\gw,\nabla'\gw)\qquad &\text{in }S^{N-1}_+\\
\phantom{{\bf A}(\gw):=-div'\,{\bf a}(,\nabla'\gw)}
\gw=0\qquad &\text{on }\prt S^{N-1}_+,
\EA\ee
where
\be\label{AS3-01}\BA {lll}
{\bf a}(r,\xi)=\left(\gb_q^2r^2+|\xi|^2\right)^{\frac{p-2}2}\xi,\\
{\bf B}(r,\xi)=\gb_q\Gl_{\gb_q}\left(\gb_q^2r^2+|\xi|^2\right)^{\frac{p-2}2}r
-\left(\gb_q^2r^2+|\xi|^2\right)^{\frac{q}2}.
\EA\ee
The operator ${\bf A}$ is a Leray-Lions operator which satisfies the assumptions (1.6)-(1.8) of \cite[Theorem 2.1]{BoMuPu}, and
the term ${\bf B}$ satisfies (1.9),(1.10) in the same article. Therefore the existence of a positive solution $\gw\in W^{1,p}_0(S^{N-1}_+)\cap L^\infty(S^{N-1}_+)$ is ensured whenever we can find a supersolution $\overline\gw\in W^{1,p}(S^{N-1}_+)\cap L^\infty(S^{N-1}_+)$
and a nontrivial subsolution $\underline\gw\in W^{1,p}(S^{N-1}_+)$ of $(\ref{AS3-0})$ such that
\be\label{AS3-1}\BA {lll}
0\leq \underline\gw\leq \overline\gw \qquad &\text{in }S^{N-1}_+.
\EA\ee
First we note that $\eta=\eta_0$ is a supersolution if the positive constant $\eta_0$ is large enough. In order to find a subsolution, we set again $\eta=\gy^\gth$ with $\gth=\gb_q/\gb_*$ and $\gy$ as in $(\ref{AS2})$. Now $\gth> 1$, thus $\eta\in W^{1,p}_0(S^{N-1}_+)$.
As above we have
$$\BA {l}
\CT(\eta)=-\gth^{p-1}\gy^{(\gth -1)(p-1)}{\rm div'}\left(\left(\gb_*^2\gy^2+|\nabla' \gy|^2\right)^{\frac{p-2}2}\nabla'\gy\right)\\[2mm]
\phantom{\CT(\eta)}-
\gth^{p-1}(\gth -1)(p-1)\gy^{(\gth -1)(p-1)-1}\left(\gb_*^2\gy^2+|\nabla' \gy|^2\right)^{\frac{p-2}2}|\nabla'\gy|^2\\[2mm]
\phantom{\CT(\eta)}
-\gb_q\Gl_{\gb_q}\gth^{p-2}\gy^{(\gth -1)(p-1)}\left(\gb_*^2\gy^2+|\nabla' \gy|^2\right)^{\frac{p-2}2}\gy+\gth^{q}\gy^{(\gth -1)q}\left(\gb_*^2\gy^2+|\nabla' \gy|^2\right)^{\frac{q}2}.
\EA$$
Now $\gb_q\Gl_{\gb_q}\gth^{p-2}=\gb_*\Gl_{\gb_q}\gth^{p-1}=\gb_*(\Gl_{\gb_q}-\Gl_{\gb_*})\gth^{p-1}+\gb_*\Gl_{\gb_*}\gth^{p-1}$ and
$\Gl_{\gb_q}-\Gl_{\gb_*}=(\gb_q-\gb_*)(p-1)=\gb_*(p-1)(\gth-1)$, hence
$$\BA {l}
\CT(\eta)=-\gth^{p-1}\gy^{(\gth -1)(p-1)}{\rm div'}\left(\left(\gb_*^2\gy^2+|\nabla' \gy|^2\right)^{\frac{p-2}2}\nabla'\gy\right)\\[2mm]
\phantom{\CT(\eta)}-
\gth^{p-1}(\gth -1)(p-1)\gy^{(\gth -1)(p-1)-1}\left(\gb_*^2\gy^2+|\nabla' \gy|^2\right)^{\frac{p-2}2}|\nabla'\gy|^2\\[2mm]
\phantom{\CT(\eta)}
-\gb_*(\Gl_{\gb_q}-\Gl_{\gb_*})\gth^{p-1}\gy^{(\gth -1)(p-1)}\left(\gb_*^2\gy^2+|\nabla' \gy|^2\right)^{\frac{p-2}2}\gy\\[2mm]
\phantom{\CT(\eta)}-\gb_*\Gl_{\gb_*}\gth^{p-1}\gy^{(\gth-1)(p-1)}\left(\gb_*^2\gy^2+|\nabla' \gy|^2\right)^{\frac{p-2}2}\gy+\gth^{q}\gy^{(\gth -1)q}\left(\gb_*^2\gy^2+|\nabla' \gy|^2\right)^{\frac{q}2}.
\EA$$
Using the equation satisfied by $\psi$ yields to the relation
$$\BA {l}
\CT(\eta)=
-\gth^{p-1}(\gth -1)(p-1)\gy^{(\gth -1)(p-1)-1}\left(\gb_*^2\gy^2+|\nabla' \gy|^2\right)^{\frac{p-2}2}|\nabla'\gy|^2\\[2mm]
\phantom{\CT(\eta)=}
-\gb_*^2(p-1)(\gth-1)\gth^{p-1}\gy^{(\gth -1)(p-1)-1}\left(\gb_*^2\gy^2+|\nabla' \gy|^2\right)^{\frac{p-2}2}\gy^2\\[2mm]
\phantom{\CT(\eta)=}+\gth^{q}\gy^{(\gth -1)q}\left(\gb_*^2\gy^2+|\nabla' \gy|^2\right)^{\frac{q}2}\\[2mm]
\phantom{\CT(\eta)}
=-\gth^{p-1}(\gth-1)(p-1)\gy^{(\gth -1)(p-1)-1}\left(\gb_*^2\gy^2+|\nabla' \gy|^2\right)^{\frac{p}2}\\
\phantom{\CT(\eta)=}+\gth^{q}\gy^{(\gth -1)q}\left(\gb_*^2\gy^2+|\nabla' \gy|^2\right)^{\frac{q}2}.
\EA$$
If we replace $\eta:=\eta_1=\gy^\gth$ by $\eta:=\eta_m=(m\gy)^\gth$ in the above computation, the inequality $\CT\eta_m)\le 0$ will be true provided
$$m^{\gth(q+1-p)}\gy^{(\gth-1)(q+1-p)+1}\le \gth^{p-1-q}(\gth-1)(p-1)\left(\gb_*^2\gy^2+|\nabla' \gy|^2\right)^{\frac{p-q}2},$$
which is satisfied if we choose $m$ small enough so that $(m\gy)^\gth\leq \eta_0$ and  satisfying
$$m^{\gth(q+1-p)}\le \gb_*^{(\gth-1)(q+1-p)+1}\gth^{p-1-q}(\gth-1)(p-1)\frac{\min_{x\in S_+^{N-1}} \left(\gb_*^2\gy^2+|\nabla' \gy|^2\right)^{\frac{p-q}2}}{\max _{x\in S_+^{N-1}} \gy^{(\gth-1)(q+1-p)+1}}.$$
Therefore $0<\eta_m\leq \eta_0$ and standard regularity implies that the solution $\gw$ is $C^{1}$ in $\overline S^{N-1}_+$. Actually $\gw$ is $C^{\infty}$ since the operator is not degenerate.
\smallskip

\noindent {\it Uniqueness}. We use the tangency method developed in the proof of \rth{nonexist}. Assume $\gw_1$ and $\gw_2$ are two positive solutions of $(\ref{AS2})$, then they are positive in $S^{N-1}_+$ and $\prt_{\bf n}\gw_i<0$ on
$\prt S^{N-1}_+$. Either the $\gw_i$ are ordered and $\gw_1\leq\gw_2$, or their graphs intersect. In any case we can define
$$\gt=\inf\{s>1:s\gw_1\geq \gw_2\}.
$$
We set $\gw^*=\gt\gw_1$. Then either the graphs of $\gw_2$ and $\gw^*$ are tangent at some interior point $\ga$, or they are not tangent in $S^{N-1}_+$, $\prt_{\bf n}\gw^*\leq \prt_{\bf n}\gw_2<0$ on $\prt S^{N-1}_+$ and there exists $\ga\in \prt S^{N-1}_+$ such that $\prt_{\bf n}\gw^*(\ga)= \prt_{\bf n}\gw_2(\ga)<0$. Furthermore $\CT(\gw^*)\geq 0$. If we set
$w=\gw^*-\gw_2$, then, as in \rth{nonexist},
$$-{\rm div'}(\tilde A\nabla'w)+\langle \tilde B,\nabla' w\rangle_g+\tilde Cw=\tilde\CL w\geq 0
$$
where
\be\label{AS16}\BA {ll}\langle \tilde AX,X\rangle_g=\left(\gb_q^2\gw_t^2+|\nabla' \gw_t|^2\right)^{\frac{p-4}2}
\left(p-2)\langle \nabla'\gw_t,X\rangle_g^2+|\nabla' \gw_t|^2|X|^2\right)\\\phantom{\langle AX,X\rangle_g}
\geq \left(\gb_q^2\gw_t^2+|\nabla' \gw_t|^2\right)^{\frac{p-4}2}\min\{1,p-1\}|\nabla' \gw_t|^2|X|^2,
\EA\ee
in which $\gw_t=\gw_2+t(\gw^*-\gw_2)$ and $t\in (0,1)$ is obtained by applying the mean value theorem and $\tilde B$ and $\tilde C$ are defined accordingly. Since $\tilde\CL$ is uniformly elliptic and has bounded coefficients, it follows from the strong maximum principle that $w=0$. Thus $\gw^*=\gt\gw_1=\gw_2$ and $\gt=1$ from the equation. This ends the proof.\qeda


\subsection{Removable boundary singularities}

The following is the basic result for removability of isolated singularities. It is valid in the general case, but with a local geometric constraint.

\bth {remov} Assume $q^*\leq q<p\leq N$, $\Gw$ is a $C^2$ bounded domain with $0\in\prt\Gw$, such that  $\Gw\cap B_\gd=B_\gd^+$ for some $\gd>0$. If $u\in C^1(\overline\Gw\setminus\{0\})$ is a nonnegative solution of $(\ref{A1})$ in $\Gw$ which vanishes on $\prt\Gw\setminus\{0\}$, then it is identically $0$.
\es
\noindent\Proof {\it Step 1: Assume $\Gw\subset\BBR^N_+$}.
For $\ge>0$, we set
$\Gw'_\ge=\Gw\cap \overline {B^c_\ge}$  and $H_\ge=\BBR^N_+\cap \overline {B^c_\ge}$. For $k,n\in\BBN_*$, $n\geq{\rm diam\,}(\Gw)$, we denote by  $v_{k,n,\ge}$ ($n\in\BBN_*$) the solution of the problem
\bel{brem1-n}
\left.\BA {ll}
-\Gd_pv+\abs{\nabla v}^q=0\qquad&\text{in }H_\ge\cap B_n\\
\phantom{-\Gd_p+\abs{\nabla v}^q}
v=k\chi_{_{\BBR^N_+\cap \prt B_\ge}}\qquad&\text{on }\prt (H_\ge\cap B_n).
\EA\right.
\ee
If $k>c_2\ge^{\frac{q-p}{q+1-p}}$ for a suitable $c_2=c_2(p,q)>0$ (see \rlemma{lest1}), then  $v_{k,n,\ge}\geq u$ in $\Gw'_\ge$. Moreover there holds $v_{k,n,\ge}\leq v_{k',n',\ge}$ for $n\leq n'$ and $k\leq k'$. Furthermore the function
$$U_{\ge,n}(x)=c_2\left((\abs{x}-\ge)^{\frac{q-p}{q+1-p}}-(n-\ge)^{\frac{q-p}{q+1-p}}\right)$$
 is a supersolution in $B_n\setminus B_\ge$, and there holds
$v_{k,n,\ge}\leq U_{\ge,n}$. By monotonicity and standard a priori estimate, we obtain that $v_{k,n,\ge}\to v_\ge$ when $n, k\to\infty$  and that the function $v=v_\ge$ is solution of
\bel{brem1}
\left.\BA {ll}
-\Gd_pv+\abs{\nabla v}^q=0\qquad&\text{in }H_\ge\\
\phantom{?}
\lim_{\abs x\to \ge}v(x)=\infty
\\
\phantom{-\Gd_p+\abs{\nabla v}^q}
v=0\qquad&\text{on }\prt\BBR^N_+\cap \overline {B^c_\ge}.
\EA\right.
\ee
Furthermore
\bel{brem1e}
u(x)\leq v_\ge(x)\leq c_2(\abs{x}-\ge)^{\frac{q-p}{q+1-p}}\quad\text{in }\Gw'_\ge.
 \ee
The function $v_\ge$ may not be unique, however it is the minimal solution of the above problem since the $v_{k,n,\ge}$ is unique, and monotonicity in $n$ and $k$ holds. Actually, $v_\ge\leq v_{\ge'}$ if $0\leq \ge\leq\ge'$.
For $\ell>0$, we recall that the transformation  $v\mapsto T_\ell[v]$ defined by $(\ref{scal2})$
leaves equation $(\ref{A1})$ invariant. As a consequence of the uniqueness of the approximations we have $T_\ell[v_{k,n,\ge}]=v_{ \ell^{\frac{p-q}{q+1-p}}k,\ell^{-1}n,\ell^{-1}\ge}$, which implies
\bel{brem3}
T_\ell[v_\ge]=v_{\ell^{-1}\ge}.
\ee
Letting $\ge\to 0$,  we derive from the monotonicity with respect to $\ge$ and standard  $C^{1,\ga}$ estimates, that the following identity holds:
\bel{brem4}
T_\ell[v_0]=v_0\qquad\forall \ell>0.
\ee
The function $v_0$ is a positive and separable solution of $(\ref{A1})$ in $\BBR^N_+$ which vanishes on $\prt\Gw\setminus\{0\}$. It follows from \rth{nonexist} that $v_0=0$, and so is $u$. \smallskip

\noindent{\it Step 2: The general case}. We assume that $\Gw\cap B_\gd\subset \BBR^N_+$ and we denote by $M$ the maximum of $u$ on $\prt B_\gd\cap\Gw$. Then the function $(u-M)_+$ is a subsolution of $(\ref{A1})$ in $\Gw\cap B_\gd$ which vanishes on $\prt\Gw\cap B_\gd\setminus\{0\}$. By Step 1, it is dominated by $v_0$, which ends the proof.\qeda
\medskip

\noindent\Remark The previous result is valid if $u$ is a subsolution with the same regularity. If $u$ is no longer assumed to be nonnegative, only $u^+$ vanishes. Furthermore, the regularity of the boundary has not been used, but only the fact that $\Gw$ is locally contained into a half space to the boundary of which $0$ belongs. \medskip

\noindent\Remark If no geometric assumption is made on $\prt\Gw$, we can prove that $u(x)=o (\abs x^{-\gb_q})$ near $0$. The next result shows that the removability holds if $q>q_*$.

\bth {remov1} Assume $q^*< q<p\leq N$ and $\Gw$ is a $C^2$ bounded domain with $0\in\prt\Gw$. If $u$ is a nonnegative solution of $(\ref{A1})$ in $\Gw$ which belongs to $C^1(\overline\Gw\setminus\{0\})$ and vanishes on $\prt\Gw\setminus\{0\}$,  it is identically $0$.
\es
\noindent\Proof As it is proved in \cite{PV1}, for any smooth subdomain $S\subset S^{N-1}$, there exists a unique $\gb_{*\,s}>0$ and $\psi^*_s>0$, unique up to an homothety, such that $x\mapsto \abs x^{-\gb_{*\,s}}\psi^*_s(\abs x^{-1}x)$ is $p$ harmonic in the cone $C_S=\{x\in \BBR^N\setminus\{0\}:\abs x^{-1}x\in S\}$ and $\psi^*_s$ satisfies
   \bel{I-9s}\BA {ll}
-div'\left(\left(\gb_{*\,s}^2\psi^{*\,2}_s+|\nabla' \psi^*_s|^2\right)^{\frac{p-2}2}\nabla' \psi^*_s\right)-\gb_{*\,s}\Gl_{\gb_{*\,s}}\left(\gb_{*\,s}^2\psi^{*\,2}+|\nabla' \psi^*_s|^2\right)^{\frac{p-2}2}\psi^*_s=0\quad& \mbox{in } S
\\[2mm]
\phantom{-div'\left(\left(\gb_{*\,s}^2\psi^{*\,2}+|\nabla' \psi^*_s|^2\right)^{\frac{p-2}2}\nabla' \psi^*_s\right)-\gb_{*\,s}\Gl_{\gb_{*\,s}}\left(\gb_{*\,s}^2\psi^{*\,2}+|\nabla' \psi^*_s|^2\right)^{\frac{p-2}2}}
\psi^*_s=0\quad &\!\!\mbox{on } \prt S,
\EA\ee
Furthermore  $S\subset \tilde S\subset S^{N-1}$ implies $\gb_{*\,\tilde s}\leq \gb_{*\,s}$. Using  the  system of spherical coordinates defined in $(\ref{beta-4})$ in Appendix II, for $\ge>0$ we denote by $S:=S_\ge$ the spherical shell with vertex  the north pole $N$ and latitude angle $\gth_{N-1}\in [0,\frac{\gp}{2}+\ge]$. Because of uniqueness of $\gb_{*\,s}$, $\gb_{*\,s_\ge}\uparrow\gb_*$ as $\ge\to 0$. Therefore, if $q>q_*$, or equivalently $\gb_q<\gb_*$, there exists $\gd,\ge>0$ such that $\Gw\cap B_\gd\subset C_{S_\ge}\cap B_\gd$ and $\gb_q<\gb_{*\,s_\ge}$. Since \rth{nonexist} is valid if $S^{N-1}_+$ is replaced by $S_\ge$ and $\gb_q<\gb_{*\,s_\ge}$ it follows that $u=0$ as in the proof of
\rth {remov}, Steps 1 and 2.\qeda
\medskip

The next result, valid in the case $p=N$, is based upon the conformal invariance of the N-Laplacian. In this case the exponent $\gb_*$ corresponding to the first spherical N-harmonic eigenvalue is equal to $1$ and the corresponding spherical N-harmonic eigenfunction in $S^{N-1}_+$ is $x_N/\abs x^{2}$.

\bth {remov2} Assume $N-\frac{1}{2}\leq q< N$, $\Gw$ is a bounded domain and $0\in\prt\Gw$ is such that there exists a ball $B\subset\Gw^c$ to the boundary of which $0$ belongs. If $u$ is a nonnegative solution of
\bel{brem5}
-\Gd_Nu+\abs{\nabla u}^q=0\qquad\text{in }\Gw,
\ee
which belongs to $C(\overline\Gw\setminus\{0\})\cap W^{1,N}_0(\Gw\setminus \overline B_\ge(0))$ for any $\ge>0$, it is identically $0$.
\es
\noindent\Proof We  assume that  the inward normal unit vector to $\prt\Gw$ at $0$ is ${\bf e}_N=(0,0,...,1)$  and that the ball $B=B_{\frac{1}{2}}(a)$ of center $a=-\frac{1}{2}{\bf e}_N$ and radius $\frac{1}{2}$ touches $\prt\Gw$ at $0$ and is exterior to $\Gw$ (this can be assumed up to a rotation and a dilation). This is the consequence of the exterior sphere condition at the point $0$. It is always valid if $\prt\Gw$ is $C^2$. We denote by $\CI_\gw$ the inversion of center $\gw=-{\bf e}_N$ and power 1, i.e.
$\CI_\gw(x)=\gw+\frac{x-\gw}{\abs{x-\gw}^2}$. Under this transformation, the
complement of the ball $B_{\frac{1}{2}}(a)$, which contains $\Gw$, is transformed into the half space $\BBR^N_-$ which contains the image $\tilde\Gw$ of $\Gw$. Since $u$ satisfies $(\ref{brem5})$, $\tilde u=u\circ\CI_\gw$ satisfies
\bel{brem6}
-\Gd_N\tilde u+\abs{x-\gw}^{2(q-N)}\abs{\nabla \tilde u}^q=0\qquad\text{in }\tilde \Gw.
\ee
Furthermore since $0=\CI_\gw(0)$ and $\CI_\gw$ is a diffeomorphism, $\tilde u\in C(\overline{\tilde \Gw}\setminus\{0\})\cap C^1(\tilde \Gw)$  and it vanishes on $\prt\tilde\Gw\setminus\{0\}$. Since $\abs{x-\gw}\leq 1$ and $q<N$, $\tilde u$ is a subsolution for $(\ref{brem5})$ in $\tilde G$. By \rth {remov1}, $\tilde u=0$.\qeda

\subsection{Weakly singular solutions}

The main result of this section is the following existence and uniqueness result concerning solutions of $(\ref{A1})$ with a boundary weak singularity. We recall that $\psi_*$ is unique positive solution of $(\ref{I-9})$ such that $\sup\psi_*=1$. Our first result is valid for any $1<p\leq N$ but it needs a geometric constraint on $\Gw$.

\bth {thweak1}Let $p-1<q<q_*<p\leq N$ and $\Gw\subset\BBR^N_+$ be a bounded $C^2$ domain such that $0\in\prt\Gw$. Assume that there exists $\gd>0$ such that $\Gw_\gd:=\Gw\cap B_\gd= B^+_\gd$. Then for any $k>0$ there exists a unique positive solution $u:=u_k$
of $(\ref{A1})$ in $\Gw$, which belongs to $C^1(\overline\Gw\setminus\{0\})$, vanishes on $\prt\Gw\setminus\{0\}$ and satisfies
\bel{wl1}
\lim_{x\to 0}\frac{u_k(x)}{\Psi_*(x)}=k
\ee
in the $C^1$-topology of $S^{N-1}_+$, where
$$\Psi_*(x)=\abs x^{-\gb_*}\psi_*(|x|^{-1}x).$$
\es\medskip

The proof of this theorem is long and difficult and requires a certain number of intermediate results.

\blemma {lweak1} Let the assumptions on $p$, $q$ and $\Gw$ of \rth {thweak1} be satisfied. There exists a unique positive $p$-harmonic function $\Phi_{*}$ in $\Gw$, which is continuous in  $\overline\Gw\setminus\{0\}$, vanishes on $\prt \Gw\setminus\{0\}$ and satisfies
\bel{wl2}
\lim_{x\to 0}\frac{\Phi_{*}(x)}{\Psi_*(x)}=1.
\ee
\es
\noindent\Proof
For $0<\ge<\gd$ let $v_\ge$ be the unique nonnegative $p$-harmonic function in $\Gw\setminus \overline {B^+_{\ge}}$ which is continuous in
$\overline \Gw\setminus {B^+_{\ge}}$, vanishes on $\prt \Gw\setminus {B_{\ge}}$ and  achieves the value ${\Psi_*}$ on $\prt B_\ge\cap \Gw$. Since $\Gw\subset\BBR^N_+$, $v_\ge\leq \Psi_*$ in $\Gw\setminus {B^+_{\ge}}$. Hence inequalities
$0<\ge<\ge'\leq \gd$ imply $v_\ge\leq v_{\ge'}$ in $\Gw\setminus \overline {B^+_{\ge'}}$. Because $\Psi_*\leq \delta^{-\beta_*}$, there holds
\bel{wl3}
v_\epsilon+\delta^{-\beta_*}\geq \Psi_*,
\ee
in $\Gw\setminus B_\delta^+ $. Since $v_\ge$ and $\Psi_*$ coincide on $\prt B_\ge^+$ and vanish on $\prt\BBR^N_+\cap(B_\delta^+ \setminus B_\epsilon^+)$, $(\ref{wl3})$ holds also in $B_\delta^+ \setminus B_\epsilon^+$.
Because $v_\ge\geq 0$ there holds
\bel{wl3+}
(\Psi_*-\delta^{-\beta_*})_+\leq v_\epsilon\leq \Psi_*\quad\mbox{ in }\Omega \setminus B_\epsilon^+.
\ee
By a standard regularity result $v_\ge$ converges to a function $\Phi_{*}$ continuous in $\overline\Gw\setminus\{0\}$, $p$-harmonic in $\Gw$ such that
$$ (\Psi_*-\delta^{-\beta_*})_+\leq \Phi_{*}\leq \Psi_*$$
in $ \Gw$. Therefore $(\ref{wl2})$ holds provided $\frac{x}{\abs x}$ remains in a compact subset of $S^{N-1}_+$.
Let us define a function $\tilde \phi_{*}$ by $\tilde\gf_*(x)=\abs x^{\gb_*}\Gf_*(x)$, then $\tilde\gf_*(r,\gs)\leq {\psi_*}(\gs)$ where  $r=\abs x$ and $\gs=\frac{x}{\abs x}\in S^{N-1}_+$. By standard $C^{1,\ga}$ estimates, $\tilde\gf_*(r,.)$ is relatively compact in the
$C(\overline{S^{N-1}_+})$-topology. Therefore the convergence of  $\frac{\Gf_*(x)}{{\Psi_*}(x)}$ to $1$ when $x$ to $0$ holds not only when $\frac{x}{\abs x}$ remains in a compact subset of $S^{N-1}_+$, but uniformly on $S^{N-1}_+$, which implies $(\ref{wl2})$. Uniqueness follows classically by $(\ref{wl2})$ and the maximum principle.\qeda

\blemma {lweak2} Let the assumptions on $p$, $q$ and $\Gw$ of \rth {thweak1} be satisfied. If for some $k>0$ there exists a solution $u_{k}$
of $(\ref{A1})$ in $\Gw$, which belongs to $C^1(\overline\Gw\setminus\{0\})$, vanishes on $\prt\Gw\setminus\{0\}$ and satisfies
$(\ref{wl1})$, then for any $k>0$ there exists such a solution.
\es
\noindent\Proof We notice that for any $c<1$ (resp $c>1$), $cu_{k}$ is a subsolution (resp. supersolution) of $(\ref{A1})$ in $\Gw$. Let $\Phi_*$ be as in Lemma \ref{l:lweak1}. If $c<1$, the function $ck\Phi_*$ is a supersolution of $(\ref{A1})$ which vanishes on $\prt\Gw\setminus\{0\}$. Furthermore
$$\lim_{x\to 0}\myfrac{cu_{k}(x)}{\Psi_*(x)}=ck=\lim_{x\to 0}\myfrac{ck\Phi_*(x)}{\Psi_*(x)}.
$$
Then there exists a solution $u_{ck}$ of $(\ref{A1})$ in $\Gw$ which satisfies $cu_{k}\leq u_{ck}\leq ck\Phi_*$. If $c>1$, we set
$u^*:=T_{c^\theta}[u_k]$, which means $ u^*(x)= c^{\gb_q\theta}u_k(c^\theta\,x)$ with $\gth=(\gb_q-\gb_*)^{-1}$. Then $u^*$ is a solution of
$(\ref{A1})$ in $\Gw^{c^\theta}=\frac{1}{c^\theta}\Gw$. In particular, $u^*$ satisfies the equation in $B^+_{\frac{\gd}{c^\theta}}(0)$. Since $c^\theta>1$, $B^+_{\frac{\gd}{c^\theta}}(0)\subset B^+_{\gd}(0)$. Put $m=\max\{u^*:x\in \prt B^+_{\frac{\gd}{c^\theta}}(0)\}$. The function $(u^*-m)_+$, extended by $0$ outside $B^+_{\frac{\gd}{c^\theta}}(0)$, is a subsolution of $(\ref{A1})$ in $\Gw$. Furthermore it satisfies
$$\lim_{x\to 0}\myfrac{(u^*-m)_+(x)}{\Psi_*(x)}=ck,
$$
uniformly on any compact subset of $S^{N-1}_+$. Therefore there exists a solution $u_{ck}$ of $(\ref{A1})$ in $\Gw$ which satisfies $(u^*-m)_+\leq u_{ck}\leq ck\Phi_*$, and in particular it vanishes on $\prt\Gw\setminus\{0\}$ and belongs to $C^1(\overline\Gw\setminus\{0\})$. By \cite {PSZ}, $u_{ck}$ is positive in $\Gw$. Thus $u_{ck}$ belongs to
$C^{1,\ga}(\overline{B^+_\gd}(0)\setminus\{0\})$ and satisfies
$$\abs x^{\gb_*}\abs{u_{ck}(x)}+\abs x^{1+\gb_*}\abs{\nabla u_{ck}(x)}
+\abs x^{1+\gb_*+\ga}\sup_{\scriptsize\BA{ll}|y|\leq |x|\\
\phantom{-}x\neq y\EA}\frac{\abs{\nabla u_{ck}(x)-\nabla u_{ck}(y)}}{|x-y|^\ga}\leq M
$$
by $(\ref{est2**})$. Therefore the set of functions $\{r^{\gb_*+1}\nabla u_{ck}(r,.)\}_{r>0}$ is uniformly relatively compact in the topology of uniform convergence on $\overline S_+^{N-1}$. Since it converges to $ck\nabla'\psi_* $ uniformly on compact subsets of $S_+^{N-1}$ as $r\to 0$, this convergence holds in $C(\overline {S_+^{N-1}})$. This implies
\bel{wl4}
\lim_{x\to 0}\frac{u_{ck}(x)}{\Psi_*(x)}=ck.
\ee
\qeda\medskip

The next Lemma is the keystone of our construction. Its proof is very delicate and needs several intermediate steps.

\blemma {lweak3} Under the assumptions of \rth{thweak1}  there exists a real number $R_0$ such that $0<R_0\leq \gd$ and a positive subsolution $\tilde u$ of $(\ref{A1})$ in $B^+_{R_0}$ which is Lipschitz continuous in $\overline{B^+_{R_0}}\setminus\{0\}$, vanishes on  $\overline {B^+_{R_0}}\cap\prt\BBR^N_+\setminus\{0\}$, is smaller than $\Psi_*$ and satisfies
\bel{sub1}
\lim_{x\to 0}\frac{\tilde u(x)}{\Psi_*(x)}=1.
\ee
\es
\noindent\Proof The construction of the function $\tilde u$. We  look for a subsolution under the form $\tilde u=\Psi_*-w$ for a suitable nonnegative function $w$. \smallskip

\noindent{\it Step 1: reduction of the problem}. We use spherical coordinates for a $C^1$ function $u:x\mapsto u(x)=u(r,\gs)$, $r=|x|$, $\gs=\frac{x}{|x|}$. Then $\nabla u=u_r{\bf e}+r^{-1}\nabla' u$ where ${\bf e}=\abs x^{-1}x$, $\abs{\nabla u}^2=
u^2_r+r^{-2}\abs{\nabla' u}^2$ and $\abs{\nabla u}^q=\left(u^2_r+r^{-2}\abs{\nabla' u}^2\right)^{\frac{q}{2}}
$. The expression of the p-Laplacian in spherical coordinates is
$$\BA {ll}-\Gd_p u=-\left(\left(u^2_r+r^{-2}\abs{\nabla' u}^2\right)^{\frac{p-2}{2}}u_r\right)_r
-\myfrac{N-1}{r}\left(u^2_r+r^{-2}\abs{\nabla' u}^2\right)^{\frac{p-2}{2}}u_r\\[4mm]
\phantom{-\Gd_p u---------}-\myfrac{1}{r^2}
div'\left(\left(u^2_r+r^{-2}\abs{\nabla' u}^2\right)^{\frac{p-2}{2}}\nabla' u\right).
\EA$$
Put $v(t,\sigma)=r^{\gb_*}u(r,\gs)$ with $t=\ln r\in (-\infty,\ln\gd]$, then $v$ satisfies
\bel{cons1}\BA {ll}\CQ[v]:=\\[4mm]
-\left(\left((v_t-\gb_*v)^2+\abs{\nabla' v}^2\right)^{\frac{p-2}{2}}(v_t-\gb_*v)\right)_t
-div'\left(\left((v_t-\gb_*v)^2+\abs{\nabla' v}^2\right)^{\frac{p-2}{2}}\nabla'v\right)\\[4mm]
\phantom{-\Gd_p u}
+\Gl_{\gb_*}\left((v_t-\gb_*v)^2+\abs{\nabla' v}^2\right)^{\frac{p-2}{2}}(v_t-\gb_*v)+e^{\gn t}\left((v_t-\gb_*v)^2+\abs{\nabla' v}^2\right)^{\frac{q}{2}}=0
\EA\ee
in $(-\infty,\ln\gd)\ti S^{N-1}_+$ where $\gn=1-(q+1-p)(\gb_*+1)=1-\frac{\gb_*+1}{\gb_q+1}>0$ and $\Gl_{\gb_*}=\gb_*(p-1)+p-N$. Notice that $\psi_*$ satisfies
\bel{eigen}-div'\left(\left(\gb_*^2\psi_*^2+\abs{\nabla' \psi_*}^2\right)^{\frac{p-2}{2}}\nabla'\psi_*\right)
-\gb_*\Gl_{\gb_*}\left(\gb_*^2\psi_*^2+\abs{\nabla' \psi_*}^2\right)^{\frac{p-2}{2}}\psi_*=0,
\ee
hence it is a supersolution for $(\ref{cons1})$. We look for a subsolution under the form
$$V(t,\gs)=\psi_*-a(t)g(\psi_*)$$
where $g$ is a continuous increasing function defined on $\BBR_+$, vanishing at $0$ and smooth on $\BBR^*_+$ and
 $a(t)=e^{\gg t}$ with $\gg>0$  to be chosen. Thus $a'=\gg a$, $a''=\gg^2 a$, $V_t=-\gg ag(\psi_*)$,  $V_t-\gb_*V=-\gb_*\psi_*+a(\gb_*-\gg)g(\psi_*)$, $\nabla' V=(1-ag'(\psi_*))\nabla'\psi_*$ and
$$\BA{ll}(V_t-\gb_*V)^2+\abs{\nabla' V}^2=\left(-\gb_*\psi_*+a(\gb_*-\gg)g(\psi_*)\right)^2+(1-ag'(\psi_*))^2\abs{\nabla'\psi_*}^2\\[4mm]
\phantom{\abs{\nabla V}}
=\left(\gb_*^2\psi_*^2+2a\gb_*(\gg-\gb_*)g(\psi_*)\psi_*\right)+\left(1-2ag'(\psi_*)\right)\abs{\nabla'\psi_*}^2
+O(a^2\norm{g(\psi)}_{C^1})\\[4mm]
\phantom{\abs{\nabla V}}
=\gb_*^2\psi_*^2+\abs{\nabla'\psi_*}^2+2a\left(\gb_*(\gg-\gb_*)\psi_*g(\psi_*)-g'(\psi_*)\abs{\nabla\psi_*}^2\right)+O(a^2\norm{g(\psi_*)}_{C^1}).
\EA$$
Therefore
$$\BA{ll}\left((V_t-\gb_*V)^2+\abs{\nabla' V}^2\right)^{\frac{p-2}{2}}\\[4mm]
\phantom{}
=\left(\gb_*^2\psi_*^2+\abs{\nabla'\psi_*}^2\right)^{\frac{p-2}{2}}\left[1+(p-2)a\myfrac{\gb_*(\gg-\gb_*)\psi_*g(\psi_*)-g'(\psi_*)\abs{\nabla\psi_*}^2
}{\gb_*^2\psi_*^2+\abs{\nabla'\psi_*}^2}\right]\\[4mm]
\phantom{---------------------------}
+O(a^2\norm{g(\psi)}_{C^1}),
\EA$$
and
$$\BA {ll}e^{\gn t}\left((V_t-\gb_*V)^2+\abs{\nabla' V}^2\right)^{\frac{q}{2}}
\\[4mm]
\phantom{-}
=e^{\gn t}\left(\gb_*^2\psi_*^2+\abs{\nabla'\psi_*}^2\right)^{\frac{q}{2}}\left[1+qa\myfrac{\gb_*(\gg-\gb_*)\psi_*g(\psi_*)-g'(\psi_*)\abs{\nabla\psi_*}^2
}{\gb_*^2\psi_*^2+\abs{\nabla'\psi_*}^2}\right]\\[4mm]
\phantom{--------------------------}
+O(e^{\gn t}a^2\norm{g(\psi_*)}_{C^1}),
\EA$$
thus
$$\BA{ll}
\left((V_t-\gb_*V)^2+\abs{\nabla' V}^2\right)^{\frac{p-2}{2}}(V_t-\gb_*V)\\[4mm]\phantom{-----}
=-\gb_*\left(\gb_*^2\psi_*^2+\abs{\nabla'\psi_*}^2\right)^{\frac{p-2}{2}}\psi_*+
a(\gb_* -\gg)\left(\gb_*^2\psi_*^2+\abs{\nabla'\psi_*}^2\right)^{\frac{p-2}{2}}g(\psi_*)\\[4mm]\phantom{--,,----}
-a\gb_*(p-2)\myfrac{\gb_*(\gg-\gb_*)\psi_*g(\psi_*)-g'(\psi_*)\abs{\nabla\psi_*}^2
}{(\gb_*^2\psi_*^2+\abs{\nabla'\psi_*}^2)^{\frac{4-p}{2}}}
\psi_*
+O(a^2\norm{g(\psi_*)}_{C^1}).
\EA$$
Finally,
\bel{cons2}\BA{ll}
-\left(\left((V_t-\gb_*V)^2+\abs{\nabla' V}^2\right)^{\frac{p-2}{2}}(V_t-\gb_*V)\right)_t\\[4mm]
\phantom{----}=a\left[
(\gg^2-\gb_* \gg)\left(\gb_*^2\psi_*^2+\abs{\nabla'\psi_*}^2\right)^{\frac{p-2}{2}}g(\psi_*)\right.
\\[4mm]\phantom{----}
\left.+\gb_*(p-2)\myfrac{\gb_*(\gg^2-\gb_*\gg)\psi_*g(\psi_*)-\gg g'(\psi_*)\abs{\nabla\psi_*}^2
}{(\gb_*^2\psi_*^2+\abs{\nabla'\psi_*}^2)^{\frac{4-p}{2}}}
\psi_*\right]
+O(a^2\norm{g(\psi_*)}_{C^2}).
\EA\ee
Since
$$\BA{lll}
\left((V_t-\gb_*V)^2+\abs{\nabla' V}^2\right)^{\frac{p-2}{2}}\nabla'V=\\[4mm]
\left(\gb_*^2\psi_*^2+\abs{\nabla'\psi_*}^2\right)^{\frac{p-2}{2}}(1-ag'(\psi_*))\left[1+a(p-2)\myfrac{\gb_*(\gg-\gb_*)\psi_*g(\psi_*)-g'(\psi_*)\abs{\nabla\psi_*}^2
}{\gb_*^2\psi_*^2+\abs{\nabla'\psi_*}^2}\right]
\nabla'\psi_*\\[4mm]
\\[2mm]\phantom{-----\left(\left((V_t-\gb_*V)^2+\abs{\nabla' V}^2\right)^{\frac{p-2}{2}}(V_t-\gb_*V)\right)_t}
+O(a^2\norm{g(\psi_*)}_{C^1})
\\[4mm]\phantom{\left((V_t-\gb_*V)^2+\abs{\nabla' V}^2\right)^{\frac{p-2}{2}}\nabla'V}
=\left(\gb_*^2\psi_*^2+\abs{\nabla'\psi_*}^2\right)^{\frac{p-2}{2}}\nabla'\psi_*
\\[4mm]\phantom{---}
+a\left(\gb_*^2\psi_*^2+\abs{\nabla'\psi_*}^2\right)^{\frac{p-2}{2}}\left[(p-2)\myfrac{\gb_*(\gg-\gb_*)\psi_*g(\psi_*)-g'(\psi_*)\abs{\nabla\psi_*}^2
}{\gb_*^2\psi_*^2+\abs{\nabla'\psi_*}^2}
-g'(\psi_*)\right]\nabla'\psi_*
\\[2mm]\phantom{-----\left(\left((V_t-\gb_*V)^2+\abs{\nabla' V}^2\right)^{\frac{p-2}{2}}(V_t-\gb_*V)\right)_t}
+O(a^2\norm{g(\psi_*)}_{C^1}),

\EA$$
we get similarly
\bel{cons32}\BA{ll}-div'\left(\left((V_t-\gb_*V)^2+\abs{\nabla' V}^2\right)^{\frac{p-2}{2}}\nabla'V\right)=
-div'\left(\left(\gb_*^2\psi_*^2+\abs{\nabla'\psi_*}^2\right)^{\frac{p-2}{2}}\nabla'\psi_*\right)
\\[4mm]
-a\,div'\left(\left(\gb_*^2\psi_*^2+\abs{\nabla'\psi_*}^2\right)^{\frac{p-2}{2}}\left[(p-2)\myfrac{\gb_*(\gg-\gb_*)\psi_*g(\psi_*)-g'(\psi_*)\abs{\nabla\psi_*}^2
}{\gb_*^2\psi_*^2+\abs{\nabla'\psi_*}^2}
-g'(\psi_*)\right]\nabla'\psi_*\right)
\\[2mm]\phantom{-----\left(\left((V_t-\gb_*V)^2+\abs{\nabla' V}^2\right)^{\frac{p-2}{2}}(V_t-\gb_*V)\right)_t}
+O(a^2\norm{g(\psi_*)}_{C^2}).
\EA\ee
Noting that
\bel{cons4}\BA{ll}-div'\left(\left(\gb_*^2\psi_*^2+\abs{\nabla'\psi_*}^2\right)^{\frac{p-2}{2}}\nabla'\psi_*\right)\psi_*
=\gb_*\Gl_{\gb_*}\left(\gb_*^2\psi_*^2+\abs{\nabla'\psi_*}^2\right)^{\frac{p-2}{2}}\psi_*,
\EA\ee
 we obtain
\bel{cons5}\BA{ll}
e^{-\gg t}\CQ[V]\\[4mm]
=\left[
(\gg^2-\gb_* \gg)\left(\gb_*^2\psi_*^2+\abs{\nabla'\psi_*}^2\right)^{\frac{p-2}{2}}g(\psi_*)
+\gb_*(p-2)\myfrac{\gb_*(\gg^2-\gb_*\gg)\psi_*g(\psi_*)-\gg g'(\psi_*)\abs{\nabla\psi_*}^2
}{(\gb_*^2\psi_*^2+\abs{\nabla'\psi_*}^2)^{\frac{4-p}{2}}}
\psi_*\right]
\\[4mm]
-div'\left(\left(\gb_*^2\psi_*^2+\abs{\nabla'\psi_*}^2\right)^{\frac{p-2}{2}}\left[(p-2)\myfrac{\gb_*(\gg-\gb_*)\psi_*g(\psi_*)-g'(\psi_*)\abs{\nabla\psi_*}^2
}{\gb_*^2\psi_*^2+\abs{\nabla'\psi_*}^2}
-g'(\psi_*)\right]\nabla'\psi_*\right)
\\[4mm]
-\Gl_{\gb_*}\left((\gg-\gb_* )\left(\gb_*^2\psi_*^2+\abs{\nabla'\psi_*}^2\right)^{\frac{p-2}{2}}g(\psi_*)
+\gb_*(p-2)\myfrac{\gb_*(\gg-\gb_*)\psi_*g(\psi_*)-g'(\psi_*)\abs{\nabla\psi_*}^2
}{(\gb_*^2\psi_*^2+\abs{\nabla'\psi_*}^2)^{\frac{4-p}{2}}}\psi_*\right)\\[4mm]
\;\;+e^{(\gn-\gg) t}\left(\gb_*^2\psi_*^2+\abs{\nabla'\psi_*}^2\right)^{\frac{q}{2}}\left[1+qa\myfrac{\gb_*(\gg-\gb_*)\psi_*g(\psi_*)-g'(\psi_*)\abs{\nabla\psi_*}^2
}{\gb_*^2\psi_*^2+\abs{\nabla'\psi_*}^2}\right]+O(a\norm{g(\psi_*)}_{C^2}).
\EA\ee
In this expression we have in particular
\bel{cons6}\BA{ll}
-div'\left(\left(\gb_*^2\psi_*^2+\abs{\nabla'\psi_*}^2\right)^{\frac{p-2}{2}}\left[(p-2)\myfrac{\gb_*(\gg-\gb_*)\psi_*g(\psi_*)-g'(\psi_*)\abs{\nabla\psi_*}^2
}{\gb_*^2\psi_*^2+\abs{\nabla'\psi_*}^2}
-g'(\psi_*)\right]\nabla'\psi_*\right)\\[4mm]\phantom{-----}
=(p-1)div'\left[g'(\psi_*)\left(\gb_*^2\psi_*^2+\abs{\nabla'\psi_*}^2\right)^{\frac{p-2}{2}}\nabla\psi_*\right]
\\[4mm]\phantom{-----}
-\gb_*div'\left(\left(\gb_*^2\psi_*^2+\abs{\nabla'\psi_*}^2\right)^{\frac{p-4}{2}}\left[(p-2)\gb_*\psi_*g'(\psi_*)+(p-2)(\gg-\gb_*)g(\psi_*)\right]\psi_*\right)
\\[4mm]\phantom{-----}
=(p-1)g''(\psi_*)\left(\gb_*^2\psi_*^2+\abs{\nabla'\psi_*}^2\right)^{\frac{p-2}{2}}|\nabla\psi_*|^2
\\[4mm]\phantom{-----}
+(p-1)g'(\psi_*)div'\left(\left(\gb_*^2\psi_*^2+\abs{\nabla'\psi_*}^2\right)^{\frac{p-2}{2}}\nabla\psi_*\right)
\\[4mm]\phantom{-----}
-(p-2)\gb_*div'\left[\myfrac{\left((\gg-\gb_*) g(\psi_*)\psi_*+\gb_*g'(\psi_*)\psi^2_*\right)}{\left(\gb_*^2\psi_*^2+\abs{\nabla'\psi_*}^2\right)^{\frac{4-p}{2}}}
\nabla'\psi_*\right].
\EA\ee
Using the equation $(\ref{eigen})$ satisfied by $\psi_*$, it infers that
\bel{cons7}\BA{ll}
-div'\left(\left(\gb_*^2\psi_*^2+\abs{\nabla'\psi_*}^2\right)^{\frac{p-2}{2}}\left[(p-2)\myfrac{\gb_*(\gg-\gb_*)\psi_*g(\psi_*)-g'(\psi_*)\abs{\nabla\psi_*}^2
}{\gb_*^2\psi_*^2+\abs{\nabla'\psi_*}^2}
-g'(\psi_*)\right]\nabla'\psi_*\right)\\[4mm]\phantom{-----}
=(p-1)\left(\gb_*^2\psi_*^2+\abs{\nabla'\psi_*}^2\right)^{\frac{p-2}{2}}\left(g''(\psi_*)|\nabla'\psi_*|^2
-\gb_*\Gl_{\gb_*}g'(\psi_*)\psi_*\right)
\\[4mm]\phantom{--------}
-(p-2)\gb_*div'\left[\myfrac{\left((\gg-\gb_*) g(\psi_*)\psi_*+\gb_*g'(\psi_*)\psi^2_*\right)}{\left(\gb_*^2\psi_*^2+\abs{\nabla'\psi_*}^2\right)^{\frac{4-p}{2}}}
\nabla'\psi_*\right].
\EA\ee
Plugging this identity into  the expression $(\ref{cons5})$, we obtain after some simplifications
\bel{cons8}\BA{ll}
e^{-\gg t}\CQ[V]=\left(\gb_*^2\psi_*^2+\abs{\nabla'\psi_*}^2\right)^{\frac{p-2}{2}}g(\psi_*)\CQ_1[V]
+e^{(\gn-\gg)t}R[V]+O(a\norm{g(\psi_*)}_{C^2}),
\EA\ee
where
\bel{cons9}\BA{ll}
R[V]=e^{\gn t}\left(\gb_*^2\psi_*^2+\abs{\nabla'\psi_*}^2\right)^{\frac{q}{2}}\left[1+q\myfrac{\gb_*(a'-\gb_*a)\psi_*g(\psi_*)-ag'(\psi_*)\abs{\nabla\psi_*}^2
}{\gb_*^2\psi_*^2+\abs{\nabla'\psi_*}^2}\right],
\EA\ee
and
\bel{cons10}\BA{ll}
\CQ_1[V]=
(\gg-\Gl_{\gb_*})(\gg-\gb_*)\left[1+(p-2)\myfrac{\gb^2_*\psi^2_*}{\gb^2_*\psi^2_*+|\nabla'\psi_*|^2}\right]
-(p-1)\gb_*\Gl_{\gb_*}\myfrac{\psi_*g'(\psi_*)}{g(\psi_*)}\\[4mm]\phantom{\CQ_1[V]=}
+\left[(p-4)\gb_*\Gl_{\gb_*}\psi_*-2\Gd'\psi_*\right]
\left(\gg-\gb_*\left(1-\myfrac{\psi_*g'(\psi_*)}{g(\psi_*)}\right)\right)\myfrac{\gb_*\psi_*}{\gb_*^2\psi_*^2+\abs{\nabla'\psi_*}^2}
\\[4mm]
-(p-2)\left[\myfrac{\psi_*g'(\psi_*)}{g(\psi_*)}\left((\gb_*+1)\gg-\gb_*\Gl_{\gb_*}+\gb_*\right)
+\gg-\gb_*+\gb_*\myfrac{\psi_*^2g''(\psi_*)}{g(\psi_*)}\right]\myfrac{\abs{\nabla'\psi_*}^2}{\gb_*^2\psi_*^2+\abs{\nabla'\psi_*}^2}
\\[4mm]\phantom{\CQ_1[V]=}
+(p-1)\myfrac{g''(\psi_*)}{g(\psi_*)}\abs{\nabla'\psi_*}^2.
\EA\ee
In this expression the difficult term to deal with  is $\left[(p-4)\gb_*\Gl_{\gb_*}\psi_*-2\Gd'\psi_*\right]$ since it has not a prescribed sign. However $\Gd'\psi_*=O(\psi_*)$ by $(\ref{Delta-1})$ in Appendix II.
\smallskip

\noindent {\it Step 2: The perturbation method and the computation with $g(\psi_*)=\psi_*$}. With such a choice
of function $g$
\bel{cons11}\BA{ll}
\CQ_1[V]=(\gg-\Gl_{\gb_*})(\gg-\gb_*)\left[1+(p-2)\myfrac{\gb^2_*\psi^2_*}{\gb^2_*\psi^2_*+|\nabla'\psi_*|^2}\right]-(p-1)\beta_*\Lambda_{\beta_*}
\\[4mm]\phantom{\CQ_1[V]----}
-(p-2)\left[(\gg-\Gl_{\gb_*})\gb_*+2\gg\right]\myfrac{\abs{\nabla'\psi_*}^2}{\gb_*^2\psi_*^2+\abs{\nabla'\psi_*}^2}+\gamma\ O(\psi_*^2).
\EA\ee
Equivalently
$$\BA{ll}\CQ_1[V]=\left[1+(p-2)\myfrac{\gb^2_*\psi^2_*}{\gb^2_*\psi^2_*+|\nabla'\psi_*|^2}\right]
\left(\gg^2-(\Gl_{\gb_*}+\gb_*)\gg\right)\\[4mm]\phantom{\CQ_1[V]----}
-\gamma\left[(p-2)(\gb_*+2)\myfrac{\abs{\nabla'\psi_*}^2}{\gb_*^2\psi_*^2+\abs{\nabla'\psi_*}^2}+O(\psi^2_*)\right]
\EA$$
and finally
\bel{cons12}\BA{ll}
\CQ_1[V]=\left[1+(p-2)\myfrac{\gb^2_*\psi^2_*}{\gb^2_*\psi^2_*+|\nabla'\psi_*|^2}\right]
\gg\left[\gg-(\Gl_{\gb_*}+\gb_*+(p-2)(\gb_*+2))+O(\psi^2_*)\right].
\EA\ee
Using the fact that $\gb_*>\frac{N-1}{p-1}$ if $1<p<2$ and $1<\gb_*<\frac{N-1}{p-1}$ if $2<p<N$ (see Theorem \ref{t:beta-est} in Appendix II), we have
\bel{cons13}\BA{ll}
\Gl_{\gb_*}+\gb_*+(p-2)(\gb_*+2)\geq\left\{\BA{ll}\Gl_{\gb_*}+\gb_*(p-1)\;&\text { if }p\geq 2\\[2mm]
N+3(p-2)>N-3&\text { if }1<p<2.
\EA\right.
\EA\ee
When $N=2$, we have explicitly $\gb_*=\frac{1+2\sqrt{p^2-3p+3}}{3(p-1)}$ (see \cite[Th 3.3]{KV}).  Therefore for all $N\geq 2$ and $p>1$, there holds
\bel{cons14}\BA{ll}
\Gl_{\gb_*}+\gb_*+(p-2)(\gb_*+2)>0.
\EA\ee
We fix $\ge_0>0$ such that, whenever $\psi_*\leq\ge_0$,  there holds
\bel{cons15}\BA{ll}
\Gl_{\gb_*}+\gb_*+(p-2)(\gb_*+2)+O(\psi^2_*)>\myfrac{1}{2}\left(\Gl_{\gb_*}+\gb_*+(p-2)(\gb_*+2)\right).
\EA\ee
If we fix $\gg_0>0$ such that
\bel{cons16}\BA{ll}
\gg_0<\min\left\{\myfrac{1}{2}\left(\Gl_{\gb_*}+\gb_*+(p-2)(\gb_*+2)\right),\gn,\gb_*\right\},
\EA\ee
we obtain
\bel{cons17}\BA{ll}
\CQ_1[V]\leq-\min\{1,p-1\}\gg m^2\qquad\forall\, 0<\gg\leq\gg_0,
\EA\ee
whenever $\psi_*\leq\ge_0$, for some $m$ depending only on $p$, $q$ and $N$ (through $\psi_*$ and $\gn$), which, in the same range of value of $\psi_*$,  yields to
\bel{cons17+}\BA{ll}
\left(\gb_*^2\psi_*^2+\abs{\nabla'\psi_*}^2\right)^{\frac{p-2}{2}}g(\psi_*)\CQ_1[V]\leq-c_{17}
\psi_*\qquad\forall\, 0<\gg\leq\gg_0,
\EA\ee
for some $c_{17}>0$ depending on $N,p,q$. This estimate is valid whatever is $p>1$, but only in a neighborhood of $\psi_*=0$. If we replace $g(\psi_*)=\psi_*$ by $g_k(\psi_*)=\psi_*e^{-k\psi_*}$ for $0<k<1$, and denote by $\CQ_{1,k}[V]$ the corresponding expression of $\CQ_1[V]$ which becomes now $\CQ_{1,0}[V]$. We define similarly $\CQ_k[V]$, and $\CQ[V]$ becomes $\CQ_0[V]$. Since $g'_k(\psi_*)=e^{-k\psi_*}-kg_k(\psi_*)$ and $g''_k=-2ke^{-k\psi_*}+k^2g_k(\psi_*)$, we obtain
\bel{cons18}\BA{ll}
\CQ_{1,k}[V]=\CQ_{1,0}[V]+k(p-1)\gb_*\Gl_{\gb_*}\psi_*+(p-1)\left(-\myfrac{2k}{\psi_*}+k^2\right)\abs{\nabla'\psi_*}^2\\[4mm]
\phantom{-----------}
+(2-p)\gb_*\left(-2k+k^2\right)\psi_*+O(\psi_*^2)
\EA\ee
Notice that $\nabla'\psi_*$ vanishes only at the North pole ${\bf e}_N$, thus there exists $k_0\in (0,1]$ such that
$$
 k(1-p)\gb_*\Gl_{\gb_*}\psi_*+(p-1)\left(\myfrac{2k}{\psi_*}-k^2\right)\abs{\nabla'\psi_*}^2\geq \myfrac{1}{2}(2-p)_+\gb_*\left(-2k+k^2\right)\psi_*\qquad\forall k\leq k_0
 $$
 whenever $\psi_*\leq\ge_0$ which yields to
\bel{cons19}\BA{ll}
\left(\gb_*^2\psi_*^2+\abs{\nabla'\psi_*}^2\right)^{\frac{p-2}{2}}g_k(\psi_*)\CQ_{1,k}[V]\leq -c_{18}k\qquad\forall k\leq k_0
\EA\ee
for some $c_{13}=c_{13}(N,p,q,\ge_0)$. There exists $c_{14}=c_{14}(N,p,q)>0$ such that
\bel{cons20}\BA{ll}\left(\gb_*^2\psi_*^2+\abs{\nabla'\psi_*}^2\right)^{\frac{q}{2}}\left[1+qe^{\gg t}\myfrac{\gb_*(\gg-\gb_*)\psi_*g_k(\psi_*)-g_k'(\psi_*)\abs{\nabla\psi_*}^2
}{\gb_*^2\psi_*^2+\abs{\nabla'\psi_*}^2}\right]\leq c_{14}
\EA\ee
 in $S^{N-1}_+\ti (-\infty,\ln\gd]$. Moreover
 \bel{cons21-}
 O(a\norm{g(\psi_*)}_{C^2})\leq e^{\gg t}\tilde c_k
 \ee
for some $\tilde c_k=\tilde c_k(N,p,q)>0$. We derive from $(\ref{cons19})$-$(\ref{cons21-})$
\bel{cons21}\BA{ll}
e^{-\gg t}\CQ_k[V]\leq -c_{13}k+c_{14}e^{(\gn-\gg) t}+e^{\gg t}\tilde c_k\qquad\forall k\leq k_0
\EA\ee
Thus there exists $T_k\leq\ln\gd$ such that $\CQ_k[V]\leq 0$, for all $t\leq T_k$ and provided $\psi_*\leq\ge_0$. This local estimate will be used in the construction of the  subsolution when $p\geq 2$. \smallskip

 \noindent {\it Step 3: The case $1<p<2$}. Since the function $\psi^*$ depends only on the azimuthal angle $\gth\in (0;\frac{\gp}{2}]$ we will write  $\psi_*(\gs)=\psi_*(\gth)$ and $\nabla'\psi_*(\gs)=\psi_{*\gth}(\gth){\bf n}$ where ${\bf n}$ is the downward unit vector tangent to $S^{N-1}$ in the hyperplane going through $\gs$ and the poles.
 From $(\ref{beta-7})$,
 \bel{cons22}\BA{ll}
(p-4)\gb_*\Gl_{\gb_*}\psi_*-2\Gd'\psi_*=(p-2)\left(\gb_*\Gl_{\gb_*}\psi_*+2\myfrac{\gb_*^2\psi_*+{\psi_{*}}_{\gth\gth}}{\gb_*^2\psi^2_{*}+\psi^2_{*\gth}}\right),
\EA\ee
since $\psi^{\,2}_{*\gth}=\abs{\nabla'\psi_*}^2$ and thus
 \bel{cons23}\BA{ll}
\left((p-4)\gb_*\Gl_{\gb_*}\psi_*-2\Gd'\psi_*\right)\myfrac{\gb_*\gg\psi_*}{\gb_*^2\psi^2_{*}+\psi^2_{*\gth}}\\
[4mm]\phantom{---------}
=(p-2)\gg\left(\Gl_{\gb_*}\myfrac{\gb_*^2\psi_*^2}{\gb_*^2\psi^2_{*}+\psi^2_{*\gth}}+2\gb_*\myfrac{\gb_*^2\psi_*^2+{\psi_{*}}_
{\gth\gth}\psi_*}{(\gb_*^2\psi^2_{*}+\psi^{\,2}_{*\gth})^2}
\right).
\EA\ee
From \rth{beta-est}-Step 4 in Appendix II, we know that $\gb_*^2\psi_*+{\psi_{*}}_{\gth\gth}\geq 0$, thus the contribution of this term to $\CQ_1[V]$ is nonpositive.  We replace this expression in $\CQ_1[V]$ with $g(\psi_*)=\psi_*$ and obtain
 \bel{cons24}\BA{ll}
\CQ_1[V]=(\gg-\Gl_{\gb_*})(\gg-\gb_*)\left(1+(p-2)\myfrac{\gb_*^2\psi^2_{*}}{\gb_*^2\psi^2_{*}+\psi^2_{*\gth}}\right)
-\Gl_{\gb_*}\gb_*(p-1)
\\[4mm]\phantom{\CQ_1--}
+(p-2)\gg\Gl_{\gb_*}\myfrac{\gb_*^2\psi^2_{*}}{\gb_*^2\psi^2_{*}+\psi^2_{*\gth}}-(p-2)\left((\gb_*+2)\gg-\Gl_{\gb_*}\gb_*\right)
\myfrac{\psi^2_{*\gth}}{\gb_*^2\psi^2_{*}+\psi^2_{*\gth}}
\\[4mm]\phantom{\CQ_1--}
+2\gb_*(p-2)\myfrac{\gb_*^2\psi_*^2+{\psi_{*}}_{\gth\gth}\psi_*}{(\gb_*^2\psi^2_{*}+\psi^2_{*\gth})^2}\gamma
\\[4mm]\phantom{}
\leq \gg\left(1+(p-2)\myfrac{\gb_*^2\psi^2_{*}}{\gb_*^2\psi^2_{*}+\psi^2_{*\gth}}\right)\left(\gg-\Gl_{\gb_*}-\gb_*\right)
-(p-2)\gg\myfrac{(\gb_*+2))\psi^2_{*\gth}-\Gl_{\gb_*}\gb_*^2\psi_*^2}{\gb_*^2\psi^2_{*}+\psi^2_{*\gth}}
\\[4mm]\phantom{}
\leq \!\gg\left(1+(p-2)\myfrac{\gb_*^2\psi^2_{*}}{\gb_*^2\psi^2_{*}+\psi^2_{*\gth}}\right)
\left(\!\gg-\!\left(\Gl_{\gb_*}\!+\gb_*+(p-2)\myfrac{(\gb_*+2)\psi^2_{*\gth}-\Gl_{\gb_*}\gb_*^2\psi_*^2}{(p-1)\gb_*^2\psi^2_{*}+\psi^2_{*\gth}}\right)\!\right).
\EA\ee
We can write
 \bel{cons25}\BA{ll}
\Gl_{\gb_*}+\gb_*+(p-2)\myfrac{(\gb_*+2)\psi^2_{*\gth}-\Gl_{\gb_*}\gb_*^2\psi_*^2}{(p-1)\gb_*^2\psi^2_{*}+\psi^2_{*\gth}}
\\[4mm]\phantom{\CQ_1----}=
\myfrac{\left(\Gl_{\gb_*}+(p-1)\gb_*\right)\gb_*^2\psi_*^2+\left(\Gl_{\gb_*}+\gb_*(p-1)+2(p-2)\right)\psi^2_{*\gth}}{(p-1)\gb_*^2\psi^2_{*}+\psi^2_{*\gth}}
\\[4mm]\phantom{\CQ_1----}
\geq c_{15}\left(\Gl_{\gb_*}+\gb_*(p-1)+2(p-2)\right)
\EA\ee	
for some positive constant $c_{15}$. This expression $\Gl_{\gb_*}+\gb_*(p-1)+2(p-2)$ is always positive: obviously if $N\geq 3$ and by using the explicit expression of $\gb_*$ if $N=2$. Thus there exists $\gg_0$ and $c_{16}>0$ such that $\CQ_1[V]<-c_{16}$ for $0<\gg\leq\gg_0$. The perturbation method of Step 2, is valid in the whole range of values of $\psi_*$ and we derive from $(\ref{cons17})$-$(\ref{cons17+})$ that
$(\ref{cons21})$ holds for all $k\leq k_0$ and $t\leq T_k$. Therefore $\CQ_k[V]\leq 0$.
\smallskip

 \noindent {\it Step 4: The case $p\geq 2$}. For $c>0$ to be fixed and $\psi_*\geq \ge_0$, $\gg\in (0,\gg_0]$, we take $g(\psi_*)=c\psi_*^{1-\frac{\gg}{\gb_*}}$. Then we derive from $(\ref{cons10})$:
  \bel{cons26}\BA{ll}
\CQ_1[V]=
(\gg-\Gl_{\gb_*})(\gg-\gb_*)\myfrac{(p-1)\gb^2_*\psi^2_*+|\nabla'\psi_*|^2}{\gb^2_*\psi^2_*+|\nabla'\psi_*|^2}
-(p-1)\gb_*\Gl_{\gb_*}\left(1-\myfrac{\gg}{\gb_*}\right)\\[4mm]\phantom{\CQ_1[V]}
-(p-1)\myfrac{\gg(\gb_*-\gg)}{\gb^2_*}\psi_*^{-1-\frac{\gg}{\gb_*}}\abs{\nabla'\psi_*}^2
-(p-2)(\gb_*-\gg)(\gg-\Gl_{\gb_*})\myfrac{\abs{\nabla'\psi_*}^2}{\gb_*^2\psi_*^2+\abs{\nabla'\psi_*}^2}
\\[4mm]\phantom{\CQ_1[V]}
=(1-p)\left[\gg(\gb_*-\gg)+\myfrac{\gg(\gb_*-\gg)}{\gb^2_*}\psi_*^{-1-\frac{\gg}{\gb_*}}\abs{\nabla'\psi_*}^2\right].
\EA\ee
For $k\leq k_0$ we fix $c$ such that $c\ge_0^{1-\frac{\gg}{\gb_*}}=\ge_0e^{-k\ge_0}\Longleftrightarrow c=\ge_0^{\frac{\gg}{\gb_*}}e^{-k\ge_0}$ and we define $g$ by
  \bel{cons27-}\BA{ll}
g(\psi_*)=\min\left\{\psi_*e^{-k\psi_*},\ge_0^{\frac{\gg}{\gb_*}}e^{-k\ge_0}\psi_*^{1-\frac{\gg}{\gb_*}}\right\}
=\left\{\BA {ll}\psi_*e^{-k\psi_*}\;&\text{if }0\leq\psi_*\leq\ge_0\\[2mm]
\ge_0^{\frac{\gg}{\gb_*}}e^{-k\ge_0}\psi_*^{1-\frac{\gg}{\gb_*}}
\;&\text{if }\ge_0\leq \psi_*\leq 1,
\EA\right.
\EA\ee
and we set $V(t,\gs)=\psi^*(\gs)-a(t)g(\psi_*(\gs))$ with $(t,\gs)\in (-\infty,T_k]\ti S^{N-1}_+$ and define
$\tilde u(r,\gs)=r^{-\gb_*}(\psi^*(\gs)-a(\ln r)g(\psi_*(\gs)))$ accordingly for $(r,\gs)\in (-\infty,e^{T_k}]\ti S^{N-1}_+$.  Since $\psi_*$ is a decreasing function the coincidence set $\{\gs\in S^{N-1}_+:\psi_*(\gs)=\ge_0\}$ is a  circular cone $\Gs_{\gth_0}$ with vertex $0$, axis ${\bf e}_N$ and angle $\gth_0$. We set $R_0=e^{T_k}$
$$\BA {ll}\Gg_1=\left\{x=(r,\gth)\in B^+_{R_0}:\gth_0< \gth<\frac{\gp}{2}\right\}=\left\{(r,\gs)\in [0,R_0)\ti S^{N-1}_+:0<\psi_*(\gs)<\ge_0\right\},\\[2mm]
\Gg_2=\left\{x=(r,\gth)\in B^+_{R_0}:0< \gth<\gth_0\right\}=\left\{(r,\gs)\in [0,R_0)\ti S^{N-1}_+:\ge_0<\psi_*(\gs)<1\right\},
\EA$$
and define
$$\BA {ll}\tilde u(r,\gs)=r^{-\gb_*}\left(\psi_*(\gs)-r^\gg g(\psi_*(\gs))\right)\\[2mm]
\phantom{\tilde u(r,\gs)}
=\left\{\BA {lll}
u_1(r,\gs)=r^{-\gb_*}(1-r^\gg e^{-k\psi_*(\gs)})\psi_*(\gs)\quad&\text{if }(r,\gth)\in \Gg_1\\[2mm]
u_2(r,\gs)=r^{-\gb_*}\left(1-r^\gg \ge_0^{\frac{\gg}{\gb_*}}e^{-k\ge_0}(\psi_*(\gs))^{1-\frac{\gg}{\gb_*}}\right)\psi_*(\gs)\quad&\text{if }(r,\gth)\in \Gg_2.
\EA\right.
\EA$$
The function $\tilde u$ is a subsolution separately on $\Gg_1$ and $\Gg_2$ and is Lipschitz continuous in $
\overline{\Gw}\setminus\{0\}$. If we denote by $g_1$ and $g_2$ the restriction of $g$ to $\Gg_1$ and $\Gg_2$ respectively, that  is to $\Gw_1$ and $\Gw_2$, then $g'_1(\ge_0)>g'_2(\ge_0)>0$. Let $\gz\in C^{1}_c(B^+_{R_0})$ which vanishes in neighborhoods of $0$ and $\prt B^+_{R_0}$, $\gz\geq 0$, then
  \bel{cons27}\BA{ll}\myint{\Gg_i}{}\abs{\nabla \tilde u}^{p-2}\nabla \tilde u.\nabla\gz dx+\myint{\Gw_i}{}\abs{\nabla \tilde u}^{q}\gz dx
\leq \myint{\Gs_{\gth_0}}{}\abs{\nabla u_i}^{p-2}\prt_{{\bf n}_i} u_i\gz dS,
\EA\ee
 where ${\bf n}_i$ is the normal unit vector on $\Gs_{\gth_0}$ outward from $\Gg_i$. Actually, ${\bf n}_2=-{\bf n}_1={\bf n}$
  thus
$$\nabla\tilde u=\tilde u_r{\bf e}+r^{-\gb_*-1}(1-r^\gg g'(\psi_*))\nabla'\psi_*=\tilde u_r{\bf e}+r^{-\gb_*-1}(1-r^\gg g'(\psi_*))\psi_{*\gth}\,{\bf n}.
$$
and on $\Gs_{\gth_0}$,
$$
\nabla\tilde u=\left\{\BA {ll}\tilde u_r{\bf e}-r^{-\gb_*-1}(1-r^\gg g_1'(\ge_0))\psi_{*\gth}\,{\bf n}\quad&\text{ in }\Gg_1\\[2mm]
\tilde u_r{\bf e}+r^{-\gb_*-1}(1-r^\gg g_2'(\ge_0))\psi_{*\gth}\,{\bf n}\quad&\text{ in }\Gg_2
\EA\right.$$
Therefore
$$\BA {ll}\abs{\nabla u_1}^{p-2}\prt_{{\bf n}_1}u_1\\[2mm]\phantom{----}
=-r^{-\gb_*-1}(1-r^\gg g_1'(\ge_0))\left(\tilde u^2_r+r^{-2\gb_*-2}(1-r^\gg g_1'(\ge_0))^2\psi^2_{*\gth}\right)^{\frac{p-2}{2}}\psi_{*\gth}\quad\text{ in }\Gg_1
\EA$$
and
$$\BA {ll}\abs{\nabla u_2}^{p-2}\prt_{{\bf n}_2}u_2\\[2mm]\phantom{----}
=r^{-\gb_*-1}(1-r^\gg g_2'(\ge_0))\left(\tilde u^2_r+r^{-2\gb_*-2}(1-r^\gg g_2'(\ge_0))^2\psi^2_{*\gth}\right)^{\frac{p-2}{2}}\psi_{*\gth}\quad\text{ in }\Gg_2.
\EA$$
By adding the two inequalities $(\ref{cons27})$
  \bel{cons28}\BA {ll}\myint{\Gw}{}\abs{\nabla \tilde u}^{p-2}\nabla \tilde u.\nabla\gz dx+\myint{\Gw}{}\abs{\nabla \tilde u}^{q}\gz dx
\leq \myint{\Gs_{\gth_0}}{}\left(\abs{\nabla u_1}^{p-2}\prt_{{\bf n_1}} u_1+\abs{\nabla u_2}^{p-2}\prt_{{\bf n_2}} u_2\right)\gz dS.
\EA\ee
 By monotonicity of the function $X\mapsto \left(\tilde u^2_r+X^2\right)^{\frac{p}{2}}$ and since
 $$r^{-\gb_*-1}(1-r^\gg g_2'(\ge_0))\geq r^{-\gb_*-1}(1-r^\gg g_1'(\ge_0))\geq 0,
 $$
 we derive
$$\BA {ll}r^{-\gb_*-1}(1-r^\gg g_2'(\ge_0))\left(\tilde u^2_r+r^{-2\gb_*-2}(1-r^\gg g_2'(\ge_0))^2\psi^2_{*\gth}\right)^{\frac{p-2}{2}}
 \\[2mm]\phantom{---}
 \geq r^{-\gb_*-1}(1-r^\gg g_1'(\ge_0))\left(\tilde u^2_r+r^{-2\gb_*-2}(1-r^\gg g_1'(\ge_0))^2\psi^2_{*\gth}\right)^{\frac{p-2}{2}}
 \EA$$
We derive that the right-hand side of $(\ref{cons28})$ is nonpositive because $\psi_{*\gth}\leq 0$, and therefore $\tilde u$ is a positive subsolution of $(\ref{A1})$ in $B^+_{R_0}$ dominated by $\Psi_*$ and satisfying $(\ref{sub1})$.\qeda
\medskip

\noindent{\it Proof of \rth{thweak1}}. Let $M=\max\{\Psi_*(x):x\in\prt B^+_{R_0}\}$, then $M=R^{-\gb_*}_0$. The function $u^{*}$ defined by
$$u^{*}(x)=\left\{\BA {ll} (\tilde u(x)-M)_+\qquad&\text{ if }\,x\in B^+_{R_0}\\
0\qquad&\text{ if }\,x\in \Gw\setminus B^+_{R_0},
\EA\right.
$$
 is indeed a subsolution of $(\ref{A1})$ in whole $\Gw$ where it satisfies $u^{*}\leq \Psi_*$ and it vanishes on $\prt\Gw\setminus\{0\}$. Since $\Phi_*$ is a positive $p$-harmonic function in $\Gw$ which vanishes on $\prt\Gw\setminus\{0\}$ and satisfies $(\ref{wl2})$, it is supersolution of $(\ref{A1})$ and therefore it dominates $u^{*}$. Therefore there exists a solution $u$ of $(\ref{A1})$
  in $\Gw$ which vanishes on $\prt\Gw\setminus\{0\}$ and satisfies $u^{*}\leq u\leq \Phi_*$. This implies that $(\ref{wl1})$ holds with $k=1$ and we conclude with \rlemma{lweak2}. This ends the proof of \rlemma{lweak3}.\qeda
\medskip

When $p=N$ the statement of \rth{thweak1} holds without the flatness assumption on $\prt\Gw$. The proof of the next theorem  is an easy adaptation to the one of \rth{thweak1}, provided \rlemma{lweak1}, \rlemma{lweak2} and \rlemma{lweak3} are modified accordingly.
\bth {thweak-N}Assume $N-1<q<N-\frac{1}{2}$ and $\Gw$ be a bounded $C^2$ domain such that $0\in\prt\Gw$. Then for any $k>0$ there exists a unique positive solution $u:=u_k$
of $(\ref{brem5})$ in $\Gw$, which belongs to $C^1(\overline\Gw\setminus\{0\})$, vanishes on $\prt\Gw\setminus\{0\}$ and satisfies
uniformly with respect to  $\gs\in S^{N-1}_+$
\bel{wl1-N}
\lim_{\scriptsize\BA{cc}x\to 0\\
x/\abs x\to\gs\EA}\abs x u_k(x)=k\psi_*(\gs).
\ee
\es\smallskip

 Since $p=N$, then $\gb_*=1$ and $\psi_*(\gs)=\frac{x_N}{\abs x}=\cos\gth_{N-1}$ with the identification of $\gs$ and $\gth_{N-1}:=\gth$. In a more intrinsic manner $(\ref{wl1-N})$ can be written under the form
 \bel{clas3N}
\lim_{\tiny\BA{l}x\to 0\\
x\in\Gw\EA}|x|^2\frac{u_k(x)}{ d(x)}=k.
\ee

 We recall that if $\gw\in \BBR^N$ and $\CI_\gw$ denotes the inversion of center $\gw$ and power $1$, i.e.
$\CI_\gw(x)=\gw+\frac{x-\gw}{|x-\gw|^2}$, then $\tilde u=u\circ\CI_\gw$ satisfies $(\ref{brem6})$. \smallskip

\blemma{weak-1N} Assume $\Gw$ be a bounded $C^2$ domain such that $0\in\prt\Gw$. Then there exists a unique  N-harmonic function $\Phi_*$ in $\Gw$, which vanishes on $\prt\Gw\setminus\{0\}$ and satisfies
\bel{weak1-N}
\lim_{\scriptsize\BA{cc}x\to 0\\
x/\abs x\to\gs\EA}\abs x \Phi_*(x)=\psi_*(\gs),
\ee
uniformly with respect to  $\gs\in S^{N-1}_+$.
\es
\noindent\Proof Uniqueness is standard. Let $\gw=-{\bf e}_N\in \overline \Gw^c$, with the notations of the proof of \rth{remov2}, $\gw'=-\gw$, $a=-\frac{1}{2}{\bf e}_N$ and $a'=-a$. We can assume that the balls $B_{\frac{1}{2}}(a)$ and $B_{\frac{1}{2}}(a')$ are tangent to $\prt\Gw$ at $0$ and respectively subset of $\Gw^c$ and $\Gw$. The function $x\mapsto \Psi(x)=-\frac{x_N}{|x|^2}$ which is $N$-harmonic in $\BBR^N_-$ and  vanishes on $\prt\BBR^{N}_-\setminus\{0\}$ is transformed by the inversion $\CI_{\gw'}$ of center $\gw'$ and power $1$ into the function $\Psi_{\gw'}=\Psi\circ\CI_\gw$ which is positive and $N$-harmonic in $B_{\frac{1}{2}}(a')$ and vanishes on $\prt B_{\frac{1}{2}}(a')\setminus\{0\}$. The function $\hat\Psi=-\Psi$ which is $N$-harmonic in $\BBR^N_+$ and  vanishes on $\prt\BBR^{N}_+\setminus\{0\}$ is transformed by the inversion $\CI_{\gw}$ of center $\gw$ and power $1$ into the function $\Psi_\gw=\hat\Psi\circ\CI_{\gw}$ which is positive and $N$-harmonic in $B^c_{\frac{1}{2}}(a)$ and vanishes on $\prt B_{\frac{1}{2}}(a)\setminus\{0\}$. For $\ge>0$ we denote by $\Phi_\ge$ the solution of
\bel{weak2-N}\BA {ll}
-\Gd_N\Phi_\ge=0\qquad&\text {in }\,\Gw\cap B^c_\ge\\
\phantom{-\Gd_N}\Phi_\ge=0\qquad&\text {in }\,(B^c_{\frac{1}{2}}(a')\cap \prt B_\ge)\cup(\prt\Gw\cap B^c_\ge)
\\
\phantom{-\Gd_N}\Phi_\ge=\Psi_{\gw'}\qquad&\text {in }\,B_{\frac{1}{2}}(a')\cap \prt B_\ge.
\EA\ee
If $0<\ge'<\ge$,  $\Phi_{\ge'}\geq \Psi_{\gw'}$ in $B_{\frac{1}{2}}(a')\cap \prt B_\ge$, thus $\Phi_{\ge'}\geq \Phi_{\ge'}$ in $\Gw\cap B^c_\ge$. We also denote by $\hat U_\ge$ the solution of
\bel{weak3-N}\BA {ll}
-\Gd_N\hat \Phi_\ge=0\qquad&\text {in }\,\Gw\cap B^c_\ge\\
\phantom{-\Gd_N}\hat \Phi_\ge=0\qquad&\text {in }\,\prt\Gw\cap B^c_\ge
\\
\phantom{-\Gd_N}\hat \Phi_\ge=\Psi_{\gw}\qquad&\text {in }\,\Gw\cap \prt B^c_\ge.
\EA\ee
In the same way as above
$$0<\ge'<\ge\Longrightarrow \hat \Phi_{\ge'}\leq\hat \Phi_\ge\quad\text {in }\,\Gw\cap \prt B^c_\ge$$
Using the explicit form of $\Psi$, $\CI_\gw:x\mapsto \gw+\frac{x-\gw}{|x-\gw|^2}$ and $\CI_{\gw'}:x\mapsto \gw'+\frac{x-\gw'}{|x-\gw'|^2}$ we see that
$$\Psi_{\gw'}\lfloor_{B_{\frac{1}{2}}(a')\cap \prt B_\ge}\leq\myfrac{1+\ge}{1-\ge}\Psi_{\gw}\lfloor_{B_{\frac{1}{2}}(a')\cap \prt B_\ge},
$$
thus
$$\Phi_\ge\leq\myfrac{1+\ge}{1-\ge}\hat \Phi_\ge\quad\text{in }\Gw\cap B^c_\ge.
$$
Letting $\ge\to 0$ we conclude that $\Phi_\ge$ converges uniformly in $\overline\Gw\setminus \{0\}$ to $\Phi_*$
which vanishes on $\prt\Gw\setminus\{0\}$ and satisfies
$(\ref{weak1-N})$.\qeda\medskip

The proof of the next statement is similar to the one of \rlemma{lweak2} up to some minor modifications, so we omit it.

\blemma {weak-2N} Let the assumptions on $q$ and $\Gw$ of \rth {thweak-N} be satisfied. If for some $k>0$ there exists a solution $u_{k}$
of $(\ref{brem5})$ in $\Gw$, which belongs to $C^1(\overline\Gw\setminus\{0\})$, vanishes on $\prt\Gw\setminus\{0\}$ and satisfies
$(\ref{wl1-N})$, then for any $k>0$ there exists such a solution.
\es

\blemma {weak-3N} Under the assumptions of \rth{thweak-N}  there exists  a Lipschitz continuous nonnegative subsolution $\tilde u$ of $(\ref{brem5})$ in $\Gw$ which vanishes on  $\prt\Gw\setminus\{0\}$, is smaller than $\Phi_*$ and satisfies
\bel{weak-6N}
\lim_{\scriptsize\BA{cc}x\to 0\\
x/\abs x\to\gs\EA}\abs x \tilde u(x)=\gs,
\ee
uniformly with respect to  $\gs\in S^{N-1}_+$.
\es
\noindent\Proof Let $\tau>0$ to be fixed and let $w$ be the solution of
\bel{weak-7N}
-\Gd_Nw+|\nabla w|^q=0\qquad\text{in }B_2^-
\ee
which vanishes on $\prt B_2^-\setminus\{0\}$ and satisfies
\bel{weakx-8N}
\lim_{\scriptsize\BA{cc}x\to 0\\
{x/|x|}\to\gs\EA}|x|w(x)=\gs
\ee
in the $C^1$-topology of $S^{N-1}_-$. Its existence follows from \rth {thweak1} and this function is dominated by the  N-harmonic  function $\Phi_*$ corresponding to this domain, obtained in \rlemma {weak-1N}. By $\CI_{\gw'}$, the half-ball  $B_2^-$ is transform into the lunule $G=B_{\frac{1}{2}}(a')\setminus B_{\frac{2}{3}}(\frac{4}{3}\gw')$ and $\tilde w=w\circ \CI_{\gw'}$ satisfies
\bel{weak-8N}
-\Gd_N\tilde w+|x-\gw'|^{2(q-N)}|\nabla \tilde w|^q=0\qquad\text{in }G.
\ee
Since $|x-\gw'|\leq 1$ in $G$, $-\Gd_N\tilde w+|\nabla \tilde w|^q\leq 0$ in $G$. We extend $\tilde w$ by $0$ in $\Gw\setminus G$ and the resulting function $\tilde u$ is a subsolution of $(\ref{brem5})$ in $\Gw$ which vanishes on $\prt\Gw\setminus\{0\})$, is smaller than the N-harmonic function $\Phi_*$ obtained in \rlemma {weak-1N}, and satisfies $(\ref{weak-6N})$.\qeda\medskip

\medskip

\section{Classification of boundary singularities}
We assume that $\Gw\subset\BBR^N$ is a $C^2$ domain and $0\in\prt\Gw$. Furthermore, in order to avoid extremely technical computations, we shall assume either that $\prt\Gw$ is flat near $0$ or $p=N$. We suppose that the tangent plane to $\prt\Gw$ at $0$ is $\prt\BBR^N_+=\{x=(x',0)\}$ and the normal inward unit vector at $0$ is ${\bf e}_N$, therefore ${\bf n}=-{\bf e}_N$ in the sequel. We denote by $\gw_{s^{{N-1}}_+}$ the unique positive solution of $(\ref{AS1})$ in $S^{N-1}_+$ and by $U_{s^{{N-1}}_+}$ the corresponding singular solution of $(\ref{A1})$ in  $\BBR^N_+$ defined by
\bel{clas1}
U_{s^{{N-1}}_+}(x)=\abs x^{-\gb_q}\gw_{s^{{N-1}}_+}(\frac{x}{\abs x}).
\ee
We recall that ${\psi_*}$ is the unique positive solution of $(\ref{AS2})$ with maximum 1 and $\Psi_*$ the corresponding $p$-harmonic function
\bel{clas2}
\Psi_*(x)=\abs x^{-\gb_*}\psi_*(\frac{x}{\abs x}).
\ee
\subsection{ The case $1<p<N$}

The first statement  points out the link between weak and strong singularities.

\bprop{strong} Under the assumptions of \rth{thweak1} there exists $\lim_{k\to\infty}u_k=u_\infty$ which is the unique element of $C(\overline\Gw\setminus\{0\})\cap C^1(\Gw)$ which vanishes on   $\prt \Omega\setminus\{0\}$, satisfies $(\ref{A1})$ in $\Omega$ and
\bel{clas8}
\lim_{x\to 0}\frac{u_{\infty}(x)}{U_{s^{{N-1}}_+}(x)}=1.
\ee
\es
\noindent\Proof Uniqueness follows from $(\ref{clas8})$ and the maximum principle. For existence, since the mapping $k\mapsto u_k$ is increasing and $u_k\leq U_{s^{{N-1}}_+}$, there exists $\lim_{k\to\infty}u_k:=u_\infty\leq U_{s^{{N-1}}_+}$ and $u_\infty\in C(\overline\Gw\setminus\{0\})\cap C^1(\Gw)$. It vanishes on   $\prt B^+_\gd\setminus\{0\}$ and satisfies $(\ref{A1})$ in $B^+_\gd$. In order to take into account the domain $B^+_\gd$ in the notations, we set $u_k=u_{k,\gd}$. Since the mapping $\gd\mapsto u_{k,\gd}$ is also increasing and $u_{k,\gd}\leq k{\Psi_*}$, there also exists
$\lim_{\gd\to\infty}u_{k,\gd}:=u_{k,\infty}\leq k{\Psi_*}$
Then, for all $\ell>0$,
\bel{clas9}T_\ell[u_{k,\gd}](x)=\ell^{\gb_q}u_{k,\gd}(\ell x)=u_{k\ell^{\gb_q},\ell^{-1}\gd}(x).
\ee
Letting $k\to\infty$, we  obtain
\bel{clas10}T_\ell[u_{\infty,\gd}](x)=\ell^{\gb_q}u_{\infty,\gd}(\ell x)=u_{\infty,\ell^{-1}\gd}(x),
\ee
and letting $\gd\to\infty$, we obtain
\bel{clas11}T_\ell[u_{\infty,\infty}](x)=\ell^{\gb_q}u_{\infty,\infty}(\ell x)=u_{\infty,\infty}(x).
\ee
This implies that
\bel{clas12}u_{\infty,\infty}(r,\gs)=r^{-\gb_q}\gw'(\gs),
\ee
and $\gw'$ is a positive solution of problem $(\ref{AS1})$. Therefore $\gw'=\gw_{s^{{N-1}}_+}$ by \rth{exist}. If we let $\ell\to 0$ in $(\ref{clas9})$ and take $\abs x=1$, $x=\gs$, we derive
\bel{clas13}
\lim_{\ell\to 0}\ell^{\gb_q}u_{\infty,\gd}(\ell, \gs)=\lim_{\ell\to 0}u_{\infty,\ell^{-1}\gd}(1,\gs)=u_{\infty,\infty}(1,\gs)=\gw_{s^{{N-1}}_+}(\gs).
\ee
This convergence holds in $C^1(\overline{S^{N-1}_+})$ because of \rlemma{lest2}. This implies $(\ref{clas8})$.\qeda\medskip

The main classification result is as follows.

\bth {clasth1}Assume $1<p<N$, $p-1<q<q^*$ and $\prt\Gw\cap B_\gd=\{x=(x',0):\abs {x'}<\gd\}$, for some $\gd>0$. If $u\in C(\overline\Gw\setminus\{0\})\cap C^1(\Gw)$ is a positive solution of $(\ref{A1})$ in $\Gw$ which vanishes on   $\prt\Gw\setminus\{0\}$, then we have the following alternative:\smallskip

\noindent (i) either there exists $k\geq 0$ such that
\bel{clast1}
\lim_{x\to 0}\frac{u(x)}{{\Psi_*}(x)}=k,
\ee

\noindent (ii) or
\bel{clast2}
\lim_{x\to 0}\frac{u(x)}{U_{s^{{N-1}}_+}(x)}=1.
\ee
\es
\noindent\Proof {\it Step 1.} Assume
\bel{clast3}
\liminf_{x\to 0}\frac{u(x)}{{\Psi_*}(x)}<\infty,
\ee
then we claim that $(\ref{clast1})$ holds. We first note that if \eqref{clast3} holds, there also holds
\bel{clast4}
\liminf_{x\to 0}\frac{u(x)}{u_1(x)}<\infty,
\ee
where $u_1$ is the solution of $(\ref{A1})$ obtained in \rth{thweak1} with $k=1$. If $\{x_n\}$ is converging to $0$ and such that for some $k>0$
$$
\liminf_{x\to 0}\frac{u(x)}{u_1(x)}=k=\lim_{n\to \infty}\frac{u(x_n)}{u_1(x_n)},
$$
there also holds by the boundary Harnack inequality $(\ref{harn3})$ applied to both $u$ and $u_1$,
$$\frac{u(x_n)}{u_1(x_n)}=\frac{u(x_n)}{d(x_n)}\frac{d(x_n)}{u_1(x_n)}\geq c_5^{-2}\frac{u(x)}{u_1(x)}\quad\forall x\text{ s.t. }\abs x=\abs {x_n}.
$$
This implies in particular
$$u(x)\leq c_5^{2}(k+\ge_n)u_1(x)\qquad\forall x\text{ s.t. }\abs x=\abs {x_n}
$$
where $\{\ge_n\}$ is converging to $0_+$, and by the comparison principle
$$
u(x)\leq Ku_1(x)\qquad\forall x\in\BBR^N_+\text{ s.t. }\abs {x_n}\leq \abs x\leq\frac{\gd}{2},
$$
for some $K>0$ and all $n\in\BBN_*$. Therefore
\bel{clast5}
\limsup_{x\to 0}\frac{u(x)}{u_1(x)}<\infty.
\ee
We can assume that $k\neq 0$, otherwise $(\ref{clast1})$ holds with $k=0$ and actually $u$ remains bounded near $0$. As a consequence of the Hopf Lemma and $C^1$ regularity, there exists $K>0$ such that
\bel{clast5'}
u(x)\leq K{\Psi_*}(x)\qquad\forall x\in B^+_{\frac{\gd}{2}}.
\ee
Let $m=\max\{u(x):\abs x=\gd\}$. For $0<\gt<\gd$ we denote by $k_\gt$ the minimum of the $\gk>0$ such that $u(x)\leq \gk{\Psi_*}(x)+m$ for $\gt\leq\abs x\leq \gd$. Then $u(x)\leq k_\gt {\Psi_*}(x)+m$, and either the graphs of the mappings
$u(.)$ and $k_\gt {\Psi_*}(.)+m$ are tangent at some $x_\gt\in  B^+_\gd\setminus\overline B^+_\gt$, or they are tangent on the boundary of the domain, and the only possibility is that they are tangent on $\abs x=\gt$. Since
$$\abs{\nabla{\Psi_*}(x)}^2=\abs x^{-2(\gb_*+1)}(\gb_*^2\psi_*^2+\abs{\nabla{\psi_*}}^2),$$
it never vanishes. If we set $w=u-(k_\gt {\Psi_*}(x)+m)$, then
\bel{clast6}
-\CL w+\abs{\nabla u}^q=0
\ee
where the operator
$$\CL=\sum_{i,j}\myfrac{\prt}{\prt x_i}\left(a_{ij}\myfrac{\prt }{\prt x_j}\right)$$
is  uniformly elliptic in a neighborhood of $x_\gt$ (see \cite[Lemma 1.3]{FV}). Furthermore $w\leq 0$ and $w(x_\gt)=0$ by the strong maximum principle
$\nabla u(x_\gt)$ must vanish, which contradicts the fact that $\nabla u(x_\gt)=\nabla w(x_\gt)$ by the tangency condition, and $\nabla w(x_\gt)\neq 0$. Therefore $\abs{x_\gt}=\gt$ and $x_\gt\notin \prt\BBR^N_+$. If $\gt'<\gt$,
$k_\gt\leq k_{\gt'}$, and we set $k=\lim_{\gt\to 0}k_{\gt}$, which is finite because of $(\ref{clast5'})$. There exists $\{\gt_n\}$ such that $\gs_n:=\gt^{-1}x_{\gt_n}\to \gs_0$. Furthermore
\bel{clast7}
r^{\gb_*}u(r,\gs)\leq k_\gt {\psi_*}(\gs)+mr^{\gb_*}\quad \text{if }\gt\leq r\leq \gd\quad \text{and }\,
\gt^{\gb_*}u(\gt,\gs_\gt)= k_\gt {\psi_*}(\gs_\gt)+m\tau^{\gb_*}.
\ee
Put
\bel{clast8} u_\gt(x)=\gt^{\gb_*}u(\gt x)
\ee
Then
$$-\Gd_p u_\gt+\gt^{p-q-\gb_*(p+1-q)}\abs{\nabla u_\gt}^q=0\qquad\text{in }B^+_{\frac{\gd}{\gt}}\setminus \{0\}
$$
and, by $(\ref{clast5'})$,
$$0\leq u_\gt(x)\leq K\abs x^{-\gb_*}\qquad\text{in }B^+_{\frac{\gd}{2\gt}}\setminus \{0\}.
$$
By \rlemma{lest2}, the set of functions $\{u_\gt(.)\}$ is relatively compact in the $C_{loc}^1$ topology of
$\overline{\BBR^N_+}\setminus \{0\}$. Therefore, as $q<q^*$, there exist a sequence $\{\gt'_n\}\subset \{\gt_n\}$ converging to $0$, and a positive $p$-harmonic function $v$ in $\BBR^N_+$, continuous in $\overline{\BBR^N_+}\setminus \{0\}$ and vanishing on
$\prt\BBR^N_+\setminus \{0\}$, such that $u_{\gt'_n}\to v$, and $v$ satisfies $(\ref{clast5'})$ in $\overline{\BBR^N_+}\setminus \{0\}$. By \rth{pph} in Appendix I,  there exists $k^*$ such that $v=k^*{\Psi_*}$. In particular,
\bel{clast9}\lim_{\gt'_n\to 0}u_{\gt'_n}(1,\gs)=k^*{\psi_*}(\gs)
\ee
in the $C^1(\overline{S^{N-1}_+})$ topology. Combining $(\ref{clast7})$, $(\ref{clast8})$and $(\ref{clast9})$ we conclude that $k^*=k$ and
\bel{clast10}\lim_{\gt'_n\to 0}\gt'^{\gb_*}_nu_{\gt'_n}(1,\gs)=k{\psi_*}(\gs)
\ee
Using \rth{thweak1}, it is equivalent to
\bel{clast11}\lim_{\gt'_n\to 0}\myfrac{u(\gt'_n,\gs)}{u_k(\gt'_n,\gs)}=1
\ee
uniformly on $S^{N-1}_+$. For any $\ge>0$, there exists $n_\ge>0$ such that $n\geq n_\ge$ implies
$$u_{k-\ge}(\gt'_n,\gs)\leq u(\gt'_n,\gs)\leq u_{k+\ge}(\gt'_n,\gs)
$$
By comparison principle,
\bel{clast12}
u_{k-\ge}\leq u\leq u_{k+\ge}+m\qquad\text{in }B^+_\gd\setminus B^+_{\gt'_n},
\ee
and finally
\bel{clast13}
u_{k-\ge}\leq u\leq u_{k+\ge}+m\qquad\text{in }B^+_\gd,
\ee
Since $\ge$ is arbitrary and using again \rth{thweak1}, it implies
\bel{clast14}
\lim_{r\to 0}\myfrac{u(r,\gs)}{{\Psi_*}(r,\gs)}=k,
\ee
locally uniformly on $S^{N-1}$. But since the convergence holds in $C^1(\overline{S^{N-1}_+})$, $(\ref{clast1})$ follows.\medskip

\noindent{\it Step 2.} Assume
\bel{clast15}
\lim_{x\to 0}\frac{u(x)}{{\Psi_*}(x)}=\infty.
\ee
For any $0<\ge<\gd$ and $k>0$, there holds
 \bel{clast16}
u_k(x)\leq u(x)\leq v_\ge(x)\qquad\text{in }B_\gd^+\setminus B_\ge^+
\ee
where $v_\ge$ has been defined in $(\ref{brem1})$ and $u_k$ is given by Theorem \ref{t:thweak1}. Letting $\ge \to 0$, $k\to\infty$, and using Proposition \ref{p:strong}, we derive
 \bel{clast17}
u_\infty(x)\leq u(x)\leq v_0(x)\qquad\text{in }B_\gd^+\setminus \{0\}.
\ee
We have seen in \rth{remov} that $v_0$ is a separable solution of $(\ref{A1})$ in $\BBR^N_+$ which vanishes on $\prt \BBR^N_+\setminus\{0\}$, therefore $v_0(x)=U_{s^{{N-1}}_+}(x)$. This implies
\bel{clast18}
u_\infty(x)\leq u(x)\leq \abs x^{-\gb_q}\gw_{s^{{N-1}}_+}(\frac{x}{\abs x})\qquad\text{in }B_\gd^+\setminus \{0\}.
\ee
We conclude using \rprop{strong}.\qeda\medskip
\subsection{ The case $p=N$}
\noindent
When $p=N$, the assumption that $\prt\Gw$ is an hyperplane near $0$ can be removed. The proof of the next results is  based upon \rth{thweak-N}.
The following result is the extension to the case $p=N$ of \rprop{strong}.


\bprop{strongN} Under the assumptions of \rth{thweak-N} there exists $\lim_{k\to\infty}u_k=u_\infty$ which is the unique element of $C(\overline\Gw\setminus\{0\})\cap C^1(\Gw)$ which satisfies $(\ref{brem5})$ in $\Gw$, vanishes on   $\prt \Omega\setminus\{0\}$ and such that
\bel{clas8*}
\lim_{x\to 0}\frac{u_{\infty}(x)}{U_{s^{{N-1}}_+}(x)}=1.
\ee
\es

\noindent\Proof We denote by $u_k^{\Gw}$ the unique positive solution of $(\ref{brem5})$ satisfying $(\ref{wl1-N})$ obtained in \rth{thweak1}. Then
\bel{clas8x}
T_\ell[u_k^\Gw]=u^{\Gw^{\ell}}_{\ell^{\gb_q-\gb_*}k},
\ee
because of uniqueness. We denote by $B:=B_{\frac{1}{2}}(a)$  and $B':=B_{\frac{1}{2}}(a')$ the two balls tangent to $\prt\Gw$ at $0$ respectively interior and exterior to $\Gw$ introduced in the proof of \rlemma{weak-1N}. Estimate $(\ref{clas3N})$ implies
\bel{clas8x1}
u_k^{B'^c}\leq u_k^\Gw\leq u_k^B
\ee
the left-hand side inequality holding in $\Gw$ and the right-hand side one in $B$. Therefore
\bel{clas8x11}
T_\ell[u_k^{B'^c}]:=u^{B'^{c\,\ell}}_{\ell^{\gb_q-\gb_*}k}\leq T_\ell[u_k^\Gw]\leq T_\ell[u_k^B]:=u^{B^{\ell}}_{\ell^{\gb_q-\gb_*}k},
\ee
the domains of validity of these inequalities being modified accordingly. Using again $(\ref{clas3N})$ we obtain
\bel{clas8x12}
T_{\ell'}[u_{k'}^{B'^c}]\leq T_\ell[u_k^{B'^c}]\qquad\text{in } B'^{c\,\ell'},
\ee
for any $0<\ell'\leq \ell$ and $\ell'^{\gb_q-\gb_*}k'\leq \ell^{\gb_q-\gb_*}k$. In the same way
\bel{clas8x13}
T_{\ell'}[u_{k'}^{B}]\geq T_\ell[u_k^{B}]\qquad\text{in } B^{\ell},
\ee
for any $0<\ell'\leq \ell$ and $\ell'^{\gb_q-\gb_*}k'\geq \ell^{\gb_q-\gb_*}k$. Since $u_k^\Gw$ $u_k^B$, $u_k^{B'^c}$ are increasing with respect to $k$, they converge respectively to $u_\infty^\Gw$ $u_\infty^B$, $u_\infty^{B'^c}$ and there holds
for any $\ell>0$
\bel{clas8x14}
T_\ell[u_\infty^{B'^c}]\leq T_\ell[u_\infty^\Gw]\leq T_\ell[u_\infty^{B}],
\ee
from $(\ref{clas8x11})$ and
\bel{clas8x15}\BA{lcc}
(i)\qquad\qquad&T_{\ell'}[u_{\infty}^{B'^c}]\leq T_\ell[u_\infty^{B'^c}]\qquad &\text{in } B'^{c\,\ell'}
\\[2mm]
(ii)\qquad\qquad& T_{\ell'}[u_{\infty}^{B}]\geq T_\ell[u_\infty^{B}]\qquad&\text{in } B^{\ell}
\EA\ee
for any $0<\ell'\leq \ell$. Notice that , replacing $\ell$ by $\ell\ell'$ we can rewrite $(\ref{clas8x14})$ as follows
\bel{clas8x16}
T_{\ell'}[T_\ell[u_\infty^{B'^c}]]\leq T_{\ell'}[T_\ell[u_\infty^\Gw]]\leq T_{\ell'}[T_\ell[u_\infty^{B}]].
\ee
Because of the monotonicity with respect to $\ell$ the following limits exist
\bel{clas8x17}U^{B'^c}=\lim_{\ell\to 0}T_\ell[u_\infty^{B'^c}]\quad\text{and }\;U^{B}=\lim_{\ell\to 0}T_\ell[u_\infty^{B}].
\ee
By \rlemma{lest2} applied with $\gf(|x|)=|x|^{-\gb_q}$ and since there holds
$u_\infty^B(x)\leq c|x|^{-\gb_q}$ and $u_\infty^{B'}(x)\leq c|x|^{-\gb_q}$, we derive

\bel{clas8xx1}\BA{lll}
(i)&|\nabla T_\ell[u_\infty^B](x)|\leq c_2|x|^{-\gb_q-1}\quad&\forall x\in B^{\ell}\\
(ii) &|\nabla T_\ell[u_\infty^B](x)-\nabla T_\ell[u_\infty^B](y)|\leq c_2|x|^{-\gb_q-1-\ga}|x-y|^{\ga}\quad&\forall x,y\in B^{\ell},\;|x|\leq|y|\\
(iii) &T_\ell[u_\infty^B](x)\leq c_2|x|^{-\gb_q-1} (\dist(x,\prt B^\ell))^{\ga}\quad&\forall x\in B^{\ell},
\EA\ee
and
\bel{clas8xx2}\BA{lll}
(i)&|\nabla T_\ell[u_\infty^{B'^c}](x)|\leq c_2|x|^{-\gb_q-1}\quad&\forall x\in B'^{c\,\ell}\\
(ii) &|\nabla T_\ell[u_\infty^{B'^c}](x)-\nabla T_\ell[u_\infty^{B'^c}](y)|\leq c_2|x|^{-\gb_q-1-\ga}|x-y|^{\ga}\quad&\forall x,y\in B'^{c\,\ell},\;|x|\leq|y|\\
(iii) &T_\ell[u_\infty^{B'^c}](x)\leq c_2|x|^{-\gb_q-1} (\dist(x,\prt B'^{c\,\ell}))^{\ga}\quad&\forall x\in B'^{c\,\ell}.
\EA\ee
Thus the sets of functions $\{T_\ell[u_\infty^B]\}$ and $\{T_\ell[u_\infty^{B'}]\}$ are equicontinuous in the $C^1$-loc topology
and by uniqueness, the limit in $(\ref{clas8x17})$ below holds in this topology. Hence $U^{B'^c}$ and $U^{B^c}$ are positive solutions of $(\ref{brem5})$ in $\BBR^N_+$ which vanish on $\prt\BBR^N_+\setminus\{0\}$. Furthermore $U^{B'^c}\leq U^{B^c}$
Since for any $\ell,\ell'>0$, $T_{\ell'}[T_\ell[u_\infty^{B'^c}]]=T_{\ell\ell'}[u_\infty^{B'^c}]$, it follows $T_{\ell'}[U^{B'^c}]=U^{B'^c}$ and in the same way $T_{\ell'}[U^{B}]=U^{B}$. This means that $U^{B}$ and $U^{B'^c}$ are self-similar solutions  of $(\ref{brem5})$ in $\BBR^N_+$ and they vanish on $\prt\BBR^N_+\setminus\{0\}$. Hence
\bel{clas8xx3}\BA{lll}
U^{B}=U^{B'^c}=U_{S^{N-1}_+}.
\EA\ee
 Applying again \rlemma{lest2} to  $u_\infty^\Gw$ with $\gf(|x|)=|x|^{-\gb_q}$ we have
\bel{clas8xx}\BA{lll}
(i)&|\nabla T_\ell[u_\infty^\Gw](x)|\leq c_2|x|^{-\gb_q-1}\quad&\forall x\in \Gw^{\ell}\\
(ii) &|\nabla T_\ell[u_\infty^\Gw](x)-\nabla T_\ell[u_k^\Gw](y)|\leq c_2|x|^{-\gb_q-1-\ga}|x-y|^{\ga}\quad&\forall x,y\in \Gw^{\ell},\;|x|\leq|y|\\
(iii) &T_\ell[u_\infty^\Gw](x)\leq c_2|x|^{-\gb_q-1} (\dist(x,\prt\Gw^\ell))^{\ga}\quad&\forall x\in \Gw^{\ell}.
\EA\ee
This implies that the set of functions $\{T_\ell[u_\infty^\Gw]\}_\ell$ is equicontinuous in the $C^1$-loc topology of $\BBR^N_+$ and there  exists a sequence $\{\ell_n\}\to 0$ and a function $U$ such that $T_{\ell_n}[u_\infty^\Gw]\to U^\Gw$ in this topology of $\BBR^N_+$, and $U$ is a positive solution
of $(\ref{brem5})$ in $\BBR^N_+$ which vanishes on $\prt\BBR^N_+\setminus\{0\}$. From $(\ref{clas8x14})$ and $(\ref{clas8xx3})$ there holds $U^\Gw= U_{S^{N-1}_+}$ and therefore
\bel{clas8x17'}
\lim_{\ell\to 0}T_\ell[u_\infty^\Gw]=U_{S^{N-1}_+}.
\ee
This implies $(\ref{clas8*})$ and
\bel{clas8x18}\lim_{r\to 0}r^{\gb_q}u_\infty^\Gw(r,\gs)=\gw_{S^{N-1}_+}(\gs)
\ee
uniformly on compact subsets of $S^{N-1}_+$.\qeda\medskip

Up to minor modifications the proof of the next classification theorem is similar to the one of \rth {clasth1}.

\bth {clasth2}Assume $N-1<q<N-\frac{1}{2}$ If $u\in C(\overline\Gw\setminus\{0\})\cap C^1(\Gw)$ is a positive solution of $(\ref{brem5})$ in $\Gw$ which vanishes on   $\prt\Gw\setminus\{0\}$, then we have the following alternative:\smallskip

\noindent (i) either there exists $k\geq 0$ such that $(\ref{clast1})$ holds,\smallskip

\noindent (ii) or $(\ref{clast2})$ holds.
\es

\mysection{Appendix I: Positive $p$-harmonic functions in a half space}
In this section we prove the following rigidity result.

\bth{pph} Assume $1<p\leq N$ and $u\in C^1(\BBR^N_+)\cap C(\overline{\BBR^N_+}\setminus\{0\})$ is a positive $p$-harmonic function which vanishes on $\prt \BBR_+^N\setminus\{0\}$ and such that $\abs x^{\gb_*} u(x)$ is bounded. Then there exists $k\geq 0$ such that
\bel{a1}
u(x)=k{\Psi_*}(x)\qquad\forall x\in\BBR^N_+.
\ee
\es
\noindent\Proof Since $\abs {x}^{\gb_*} u(x)$ is bounded, $\abs x^{\gb_*+1} \nabla u(x)$ is also bounded and there exists  $m>0$ such that $u(x)\leq m{\Psi_*}(x)$ in $B_\gd^+$. We denote by $k$ the infimum of the $c>0$ such that $u(x)\leq c\Psi_*(x)$. Then
 \bel{a2}
0\leq u(x)\leq k{\Psi_*}(x)\qquad\forall x\in \BBR^N_+\setminus\{0\}
\ee
and we assume that $k>0$ otherwise $u=0$. Assume that the graphs over $\BBR^N_+$ of the functions $x\mapsto u(x)$ and $x\mapsto k\Psi_*(x)$ are tangent at some point
$x_0\in\BBR^N_+$ or $x_0\in\prt\BBR^N_+\setminus\{0\}$. Since
$\nabla \Psi_*$ never vanishes in $\overline\BBR^N_+\setminus\{0\}$ it follows from the strong maximum principle or Hopf Lemma that $u=k{\Psi_*}$. If the two graphs are not tangent in $\overline\BBR^N_+\setminus\{0\}$, either  they are asymptotically tangent at $0$, or at $\infty$. \smallskip

\noindent (i) In the first case there exists two sequences $\{k_n\}$ increasing to $k$ and $\{x_n\}\subset\BBR^N_+$ converging to zero such that $\frac{u(x_n)}{{\Psi_*}(x_n)}=k_n$. We set
$r_n=\abs{x_n}$ and $u_{r_n}(x)=r_n^{\gb_*} u(r_nx)$. Then $u_{r_n}$ is $p$-harmonic and positive and  $0<u_{r_n}(x)\leq k\abs x^{-\gb_*}{\psi_*}(\frac{x}{\abs x})$; therefore
 \bel{a3}
\abs{\nabla u_{r_n}(x)}\leq C\abs x^{-\gb_*-1}\text { and }\abs{\nabla u_{r_n}(x)-\nabla u_{r_n}(x')}\leq C\abs x^{-\gb_*-1-\ga}
\abs{x-x'}^\ga\ee
for $0<\abs x\leq \abs {x'}$ and some constants $C>0$ and $\ga\in (0,1)$. Up to a subsequence, we can assume that
$u_{r_n}$ converges to some $ U$ in the $C^1_{loc}$ topology of $\overline\BBR^N_+\setminus\{0\}$ and $\frac{x_n}{r_n}\to \xi\in S^{N-1}_+$. The function $U$ is $p$-harmonic and positive in $\BBR^N_+$ and satisfies $0\leq U\leq k{\Psi_*}$ in $\BBR^N_+$ and
$U(\xi)= k{\Psi_*}(\xi)$ if $\xi\in S^{N-1}_+$ or $U_{x_N}(\xi)= k\Psi_{*\,x_N}(\xi)$  if $\xi\in \prt S^{N-1}_+$. It follows from the strong maximum principle or Hopf Lemma that $U= k{\Psi_*}$. Therefore $u_{r_n}\to k{\Psi_*}$ and in particular
 \bel{a4}
\lim_{r_n\to 0}\myfrac{r_n^{\gb_*} u(r_n,\gs)}{{\psi_*}(\gs)}=k\quad\text{uniformly on }\,S^{N-1}_+.
\ee
For any $\ge>0$, there exists $n_\ge\in\BBN_*$ such that for $n\geq n_\ge$, $(k-\ge){\Psi_*}(x)\leq u(x)\leq (k+\ge){\Psi_*}(x)$ if $\abs x=r_n$. This implies  $(k-\ge){\Psi_*}(x)\leq u(x)\leq (k+\ge){\Psi_*}$ for $\abs x\geq r_n$ and therefore in
$\BBR^N$.
Since $\ge$ is arbitrary, we deduce that $u=k{\Psi_*}$.\smallskip

\noindent (ii) if the two graphs are tangent at infinity, there exist two sequences $\{k_n\}$ increasing to $k$ and $\{x_n\}$ such that $r_n=\abs{x_n}\to\infty$ with $u(x_n)=k_n{\Psi_*}(x_n)$ and
 \bel{a5}
\lim_{r_n\to \infty}\myfrac{r_n^{\gb_*} u(r_n,\gs)}{{\psi_*}(\gs)}=k\quad\text{uniformly on }\,S^{N-1}_+.
\ee
Therefore we look at the supremum of the $c>0$ such that $u\geq c{\Psi_*}$. If the set of such $c$ is empty, it would mean that
$$\inf_{x\in\BBR^N_+}\frac{u(x)}{{\Psi_*}(x)}=0.
$$
Clearly, if this infimum is achieved at some point, the strong maximum principle or Hopf Lemma imply $u\equiv 0$, contradicting $(\ref{a5})$, and this relation prevents also this infimum be achieved at infinity. We are left with the case where there exists a sequence $\{z_n\}\subset \BBR^N_+$, converging to $0$, such that
 \bel{a6}
\lim_{n\to \infty}\myfrac{u(z_n)}{{\Psi_*}(z_n)}=0.
\ee
By boundary Harnack inequality \cite[th 2.11]{BVBV}, there exists $c>0$ such that
 \bel{a7}
c^{-1}\myfrac{u(z)}{{\Psi_*}(z)}\leq \myfrac{u(z_n)}{{\Psi_*}(z_n)}\leq c\myfrac{u(z)}{{\Psi_*}(z)}\quad\forall z\in\BBR^N_+\text{ s.t. }\abs z=\abs{z_n}
\ee
Combining $(\ref{a6})$ and $(\ref{a7})$, we derive that
 \bel{a8}
\lim_{n\to \infty}\sup_{\abs z=\abs{z_n}}\myfrac{u(z)}{{\Psi_*}(z)}=0,
\ee

Denoting by $\ge_n$ the supremum in the above relation, we obtain that $u\leq \ge_n{\Psi_*}$ in  $\BBR^N_+\setminus B_{\ge_n}$ and finally $u=0$, contradiction. Thus we are left with the case where there exists $k'\in (0,k]$ which is the supremum of the $c>0$ such that $u\geq c{\Psi_*}$. In particular $u\geq k'{\Psi_*}$. Remembering that  $u\leq k{\Psi_*}$ we get $k=k'$, which implies $u=k{\Psi_*}$.

Next we assume that $k'<k$. Clearly the graphs of $u$ and $k'{\Psi_*}$ cannot be tangent in $\overline\BBR^N_+$, because of strong maximum principle or Hopf Lemma. They cannot be tangent at infinity because of $(\ref{a5})$. Therefore there exist two sequences $\{k'_n\}$ increasing to $k'$ and $\{x'_n\}\subset \BBR^N_+$ converging to $0$ such that  $\frac{u(x'_n)}{{\Psi_*}(x'_n)}=k'_n$. As in case (i) we obtain that
 \bel{a9}
\lim_{r'_n\to 0}\myfrac{r_n'^{\gb_*} u(r'_n,\gs)}{{\psi_*}(\gs)}=k'\quad\text{uniformly on }\,S^{N-1}_+,
\ee
where $r'_n=\abs{x'_n}$, and finally derive that $u=k'{\Psi_*}$, a contradiction with $(\ref{a5})$. Therefore $k=k'$, which ends the proof.\qeda
\medskip

\noindent\Remark In the case $p=N$ the result holds under the weaker assumption $\displaystyle\lim_{\abs x\to\infty}u(x)=0$. This is due to the fact that this condition implies by regularity
$$\lim_{\abs x\to\infty}\myfrac{u(x)}{\gw_{s^{{N-1}}_+}(\frac{x}{\abs x})}=0
$$
and therefore
$$u(x)\leq m{\Psi_*}(x)\quad\forall x\,\text{ s.t. }\abs x\geq 1,
$$
where $m=\max_{\abs x=1}\myfrac{u(x)}{\gw_{s^{{N-1}}_+}(\frac{x}{\abs x})}$. Using the inversion $x\mapsto \frac{x}{\abs x^2}$, we obtain that the estimate $u\leq m{\Psi_*}$ holds $\BBR^N$, and we conclude by \rth{pph}. \medskip

\noindent\Remark We conjecture that the rigidity result holds under the mere condition
\bel{Z}
\lim_{\abs x\to \infty}\abs x^{- \tilde \gb}u(x)=0,
\ee
were $\tilde \gb$ is the (positive) exponent corresponding to the regular spherical $p$-harmonic function under the form
\bel{Z*}
\tilde {\Psi}=\abs x^{\tilde \gb} {\tilde\psi}(\footnotesize{\frac{x}{\abs x}}),
\ee
see \cite{To}, \cite{PV1}. Note that $\tilde \gb=1$  when $p=N$.

\mysection{Appendix II: Estimates on $\beta_*$}

When $N=2$ and $1<p\leq 2$, it is proved in \cite{KV} that
\bel{beta-0}
\gb_*=\myfrac{3-p+2\sqrt{p^2-5p+7}}{3(p-1)}.
\ee

Up to now no estimate is known when $N>2$ except in the cases $p=2$ where $\gb_*=N-1$ and $p=N$ where $\gb_*=1$, besides the classical one
\bel{beta-1}
\gb_*>\myfrac{N-p}{p-1},
\ee
valid when $p<N$. In this section we prove the following result

\bth {beta-est} Assume $1<p< N$. Then the following estimates hold:
\bel{beta-2}
1<p<2\Longrightarrow \gb_*> \frac{N-1}{p-1},
\ee
\bel{beta-3}
2<p<N\Longrightarrow \max\left\{1, \frac{N-p}{p-1}\right\}<\gb_*< \frac{N-1}{p-1}.
\ee
\es

\noindent\Remark It is worth noticing that when $p=2$ or $p=N$, there holds $\gb_*=\frac{N-1}{p-1}$.  \medskip

\noindent{\it Proof of \rth{beta-est}}. We consider the following set of spherical coordinates in $\BBR^N_+$ with $x=(x_1,...,x_N)$
\bel{beta-4}\BA{ll}
x_1&\!\!\!\!\!=r\sin\gth_{N-1}\sin\gth_{N-2}...\sin\gth_{2}\sin\gth_{1}\\
x_2&\!\!\!\!\!=r\sin\gth_{N-1}\sin\gth_{N-2}...\sin\gth_{2}\cos\gth_{1}\\
\vdots\\
x_{N-1}&\!\!\!\!\!=r\sin\gth_{N-1}\cos\gth_{N-2}\\
x_N&\!\!\!\!\!=r\cos\gth_{N-1}
\EA\ee
with $\gth_1\in [0,2\gp]$ and $\gth_k\in [0,\gp]$ for $k=2,...,N-2$ and $\gth_{N-1}\in [0,\frac{\gp}{2}]$. Under this representation, a solution $\gw$ of $(\ref{AS2})$ verifies
\bel{beta-5}\BA{ll}
-\myfrac{1}{\sin^{N-2}\gth_{N-1}}\left[\sin^{N-2}\gth_{N-1}\left(\gb^2_*\gw^2+\gw^2_{\gth_{N-1}}+\myfrac{1}{\sin^{2}\gth_{N-1}}\abs{\nabla_{\gth'}\gw}^2\right)^{\frac{p-2}{2}}\gw_{\gth_{N-1}}\right]_{\gth_{N-1}}\\[4mm]
-\myfrac{1}{\sin^{2}\gth_{N-1}}div'_{\gth'}\left[\sin^{N-2}\gth_{N-1}\left(\gb^2_*\gw^2+\gw^2_{\gth_{N-1}}+\myfrac{1}{\sin^{2}\gth_{N-1}}\abs{\nabla_{\gth'}\gw}^2\right)^{\frac{p-2}{2}}\nabla_{\gth'}\gw\right]\\[4mm]
=\gb_*\Gl_{\gb_*}\left[\sin^{N-2}\gth_{N-1}\left(\gb^2_*\gw^2+\gw^2_{\gth_{N-1}}+\myfrac{1}{\sin^{2}\gth_{N-1}}\abs{\nabla_{\gth'}\gw}^2\right)^{\frac{p-2}{2}}\gw\right]
\EA\ee
where $\nabla_{\gth'}$ and $div'_{\gth'}$ denotes respectively the spherical gradient the divergence in variables
$\gth'=(\gth_1,...,\gth_{N-2})$ parameterizing $S^{N-2}$ and $\Gl_{\gb_*}$ is defined in Introduction. If $\gw$ is the unique positive solution of $(\ref{AS2})$ (up to homothety), it depends only on $\gth_{N-1}$ and is $C^\infty$. For simplicity we set $\gth_{N-1}=\gth\in [0,\frac{\gp}{2}]$ and $\gw=\gw(\gth)$ satisfies
\bel{beta-6}\BA{ll}
-\myfrac{1}{\sin^{N-2}\gth}\left[\sin^{N-2}\gth\left(\gb^2_*\gw^2+\gw^2_{\gth}\right)^{\frac{p-2}{2}}\gw_{\gth}\right]_{\gth}
=\gb_*\Gl_{\gb_*}\left[\sin^{N-2}\gth\left(\gb^2_*\gw^2+\gw^2_{\gth}\right)^{\frac{p-2}{2}}\gw\right]
\\[4mm]\phantom{-------------------------------}\text{in }(0,\frac{\gp}{2})
\\[4mm]\phantom{-------------------------}
\gw(\frac{\gp}{2})=0\,,\;\gw_\gth(0)=0.
\EA\ee
{\it Step 1: The eigenvalue identity.}
Equation $(\ref{beta-6})$ can also be written under the form

\bel{beta-7}\BA{ll}
-\gw_{\gth\gth}-(N-2)\cot\gth\,\gw_{\gth}-(p-2)\myfrac{\gb_*^2\gw+\gw_{\gth\gth}}{\gb^2_*\gw^2+\gw^2_{\gth}}\gw^2_{\gth}
=\gb_*\Gl_{\gb_*}\gw.
\EA\ee
By multiplying \eqref{beta-7} by $\cos\theta\sin^{N-2}\theta$ and then integrating over $(0,\frac{\gp}{2})$ we obtain
 $$
- \myint{0}{\frac{\gp}{2}}\left(\gw_{\gth\gth}+(N-2)\cot\gth\,\gw_{\gth}\right)\cos\gth\sin^{N-2}\gth d\gth
 =(N-1) \myint{0}{\frac{\gp}{2}}\gw\cos\gth\sin^{N-2}\gth d\gth.
 $$
 Noticing that
$$\gb_*\Gl_{\gb_*}+1-N=(p-1)\left(\gb_*-\frac{N-1}{p-1}\right)\left(\gb_*+1\right)
$$
we derive
\bel{beta-8}\BA{ll}
(2-p)\myint{0}{\frac{\gp}{2}}\myfrac{\gb_*^2\gw+\gw_{\gth\gth}}{\gb^2_*\gw^2+\gw^2_{\gth}}\gw^2_{\gth}
\gw\cos\gth\sin^{N-2}\gth d\gth\\[4mm]\phantom{-----}
=(p-1)\left(\gb_*-\myfrac{N-1}{p-1}\right)\left(\gb_*+1\right)\myint{0}{\frac{\gp}{2}}\gw\cos\gth\sin^{N-2}\gth d\gth.
\EA\ee
{\it Step 2: Elliptic coordinates and reduction}. Writing  $\gw(\gth)=\gw(0)+a\gth^2+o(\gth^2)$, $\gw_\gth(\gth)=2a\gth+o(\gth)$ and $\gw_{\gth\gth}(\gth)=2a+o(1)$, then $-Na=\gb_*\Gl_{\gb_*}$. This implies that $\gw$ is decreasing near $0$. It is immediate that it cannot have a local minimum in $(0,\frac{\gp}{2})$, therefore it remains decreasing in the whole interval.
We parameterize the ellipse
$$E_r=\{(x,y): x>0,\,y<0,\,x^2+\gb_*^{-2}y^2=r^2\}$$
by setting
$$\gw=r\cos\gf\,\text{ and}\;-\gw_\gth=\gb_* r\sin\gf\,\text{ with }\;\gf=\gf(\gth)\,\text{ and }\;r=r(\gth).$$
The  functions $r$ and $\gf$ are $C^2$. Hence
$r_\gth\cos\gf-r\sin\gf\gf_\gth=-\gb_* r\sin\gf$, then   $r_\gth\cos\gf=(\gf_\gth-\gb_*)r\sin\gf$ and
$r_\gth=(\gf_\gth-\gb_*)r\tan\gf.$
Plugging this into $(\ref{beta-7})$, we derive
\bel{beta-9}\BA{ll}
-\left((p-1)\myfrac{r_\gth}{r}+\gf_\gth\cot\gf+(N-2)\cot\gth\right)+\Gl_{\gb_*}\cot\gf=0,
\EA\ee
and finally
\bel{beta-10}\BA{ll}
(p-1)(\gf_\gth-\gb_*)\tan\gf+(\gf_\gth-\Gl_{\gb_*})\cot\gf=(2-N)\cot\gth.
\EA\ee
{\it Step 3: Estimates on $\gf_\gth$}. We can write $(\ref{beta-10})$ under the equivalent form
\bel{beta-10'}
(p-1)(\gf_\gth-\gb_*)\tan^2\gf+\gf_\gth-\Gl_{\gb_*}=(2-N)\myfrac{\cos\gth}{\cos\gf}\myfrac{\sin\gf}{\sin\gth}.
\ee
Since
$$\lim_{\gth\to 0}\myfrac{\sin\gf}{\sin\gth}=\lim_{\gth\to 0}\myfrac{\cos\gf}{\cos\gth}\gf_\gth=\gf_\gth(0),
$$
we derive $\gf_\gth(0)-\Gl_{\gb_*}=(2-N)\gf_\gth(0)$ and thus
$\gf_\gth(0)=\myfrac{\Gl_{\gb_*}}{N-1}$.
Similarly, the expansion of $\gf(\gth)$ near $\gth=\frac{\gp}{2}$ yields to
$
\gf_\gth(\frac{\gp}{2})=\gb_*.
$
Since $p<N$, $\Gl_{\gb_*}/(N-1)<\gb_*$. We claim now that
\bel{beta-13}\BA{ll}
\gf_\gth(\gth)\leq \gb_*\qquad\forall \gth\in (0,\frac{\gp}{2}).
\EA\ee
If $\Gl_{\gb_*}\leq \gb_*$, then
$$(2-N)\cot\gth=(p-1)(\gf_\gth-\gb_*)\tan\gf+(\gf_\gth-\Gl_{\gb_*})\cot\gf\geq((p-1)\tan\gf+\cot\gf)(\gf_\gth-\gb_*)
$$
thus ($\ref{beta-13}$) holds. \smallskip

\noindent Next we assume $\gb_*< \Gl_{\gb_*}$. It means
$0<(p-2)\gb_*-(N-p)$ and thus $p>2$. We claim that
\bel{beta-14}\BA{ll}
\gb_*\leq\myfrac{N-2}{p-2}.
\EA\ee
We proceed by contradiction and assume
\bel{beta-15}
\gb_*>\myfrac{N-2}{p-2}.
\ee
Then
$$(p-2)\left(\gb_*^2-\frac{N-p}{p-2}\gb_*-\frac{N-2}{p-2}\right)=(p-2)\left(\gb_*+1\right)\left(\gb_*-\frac{N-2}{p-2}\right)>0.
$$
Equivalently
$$\gb_*(\Gl_{\gb_*}-\gb_*)>N-2.
$$
Since
$$\lim_{\gth\to\frac{\gp}{2}}\cot\gth\tan\gf=\lim_{\gth\to\frac{\gp}{2}}\frac{\cos\gth}{\cos\gf}=\lim_{\gth\to\frac{\gp}{2}}\frac{\sin\gth}{\gf_\gth\sin\gf}=\frac{1}{\gb_*}$$
and
\bel{beta-15+}\BA {ll}\displaystyle(p-1)(\gf_\gth(\gth)-\gb_*)\tan^2\gf=\Gl_{\gb_*}-\gf_\gth(\gth)+(2-N)\myfrac{\cos\gth}{\cos\gf}\myfrac{\sin\gf}{\sin\gth}\\[2mm]
\phantom{(p-1)(\gf_\gth(\gth)-\gb_*)\tan^2\gf}
\displaystyle=\myfrac{1}{\gb_*}\left(\gb_*(\Gl_{\gb_*}-\gb_*)+2-N\right)+o(1),
\EA\ee
thus, if  $(\ref{beta-15})$ holds there exists $\ge>0$ such that $\gf_\gth(\gth)>\gb_*$
for any $\gth\in[\frac{\gp}{2}-\ge,\frac{\gp}{2})$.
Since $\gf_\gth(0)<\gb_*$, there exists $\bar\gth\in (0,\frac{\gp}{2})$ such that $\gf_\gth(\bar\gth)=\gb_*$ and $\gf_{\gth\gth}(\bar\gth)\geq 0$. We compute $\gf_{\gth\gth}$ and get
$$\BA {ll}(p-1)\gf_\gth(\gf_\gth-\gb_*)\sec^2\gf+\left((p-1)\tan\gf+\cot\gf\right)\gf_{\gth\gth}
-\gf_\gth(\gf_\gth-\Gl_{\gb_*})\csc^2\gf
=(N-2)\csc^2\gth
\EA$$
Hence, at $\gth=\bar\gth$
$$\gf_{\gth\gth}(\bar\gth)\left((p-1)\tan\gf(\bar\gth)+\cot\gf(\bar\gth)\right)=
\gb_*(\gb_*-\Gl_{\gb_*})\csc^2\gf(\gth)+(N-2)\csc^2\bar\gth
$$
From $(\ref{beta-10})$,
$$\cot\gf(\bar\gth)=\frac{N-2}{\Gl_{\gb_*}-\gb_*}\cot\bar\gth
$$
Therefore
\bel{beta-16}\BA{ll}
A(\bar\gth):=\gf_{\gth\gth}(\bar\gth)\left((p-1)\tan\gf(\bar\gth)+\cot\gf(\bar\gth)\right)\\[2mm]
\phantom{A(\gth):}
=
\left(1+\left(\myfrac{N-2}{\Gl_{\gb_*}-\gb_*}\right)^2\cot^2\bar\gth\right)\gb_*(\gb_*-\Gl_{\gb_*})+(N-2)(1+\cot^2\bar\gth)\\[2mm]
\phantom{A(\gth):}
=\gb_*(\gb_*-\Gl_{\gb_*})+N-2-\left(\myfrac{(N-2)^2}{\Gl_{\gb_*}-\gb_*}+2-N\right)\cot^2\bar\gth\\[2mm]
\phantom{A(\gth):}
=-(p-2)(\gb_*+1)\left(\gb_*-\myfrac{N-2}{p-2}\right)-\myfrac{N-2}{\Gl_{\gb_*}-\gb_*}\left(\gb_*(N-1)-\Gl_{\gb_*}\right)\cot^2\bar\gth\\[2mm]
\phantom{A(\gth):}
<0,
\EA\ee
using $(\ref{beta-15})$ and the fact that $N>p$. This is a contradiction, thus $(\ref{beta-14})$ holds. \smallskip

Next, if $\gb_*< \frac{N-2}{p-2}$, it follows from   $(\ref{beta-15+})$ that there exists $\ge>0$ such that  $\gf_\gth<\gb_*$ in
$[\frac{\gp}{2}-\ge,\frac{\gp}{2})$. If $(\ref{beta-13})$ is not true, there exist $0<\gth_1<\gth_2<\frac{\gp}{2}-\ge$ such that $\gf_\gth(\gth_1)=\gf_\gth(\gth_2)=\gb_*$, $\gf_{\gth\gth}(\gth_1)\geq  0$, $\gf_{\gth\gth}(\gth_2)\leq  0$. Using the equation satisfied by $\gf_{\gth\gth}$, we obtain for $i=1,2$,
\bel{beta-17}\BA{ll}
A(\gth_i)
=(2-p)(\gb_*+1)\left(\gb_*-\myfrac{N-2}{p-2}\right)-\myfrac{N-2}{\Gl_{\gb_*}-\gb_*}(\gb_*(N-1)-\Gl_{\gb_*})\cot^2\gth_i.
\EA\ee
On one hand $A(\gth_2)\leq 0\leq A(\gth_1)$, and on the other
$$A(\gth_2)- A(\gth_1)=\myfrac{N-2}{\Gl_{\gb_*}-\gb_*}(\gb_*(N-1)-\Gl_{\gb_*})(\cot^2\gth_1-\cot^2\gth_2)>0,
$$
since $\cot$ is decreasing in $(0,\frac{\gp}{2})$, $\cot^2\gth_1>\cot^2\gth_2$, a contradiction. Therefore $\gf_\gth\leq \gb_*$ in $(0,\frac{\gp}{2})$. \smallskip

Finally, if $\gb_*= \frac{N-2}{p-2}$ and the maximum of $\phi_\gth$ on $[0,\frac{\gp}{2})$ is larger than $\gb_*$ and achieved at some $\bar\theta<\frac{\gp}{2}$ the exists $\gth_1<\bar\gth$ such that $\phi_\gth(\gth_1)=\gb_*$ and $\phi_{\gth\gth}(\gth_1)\geq 0$. In that case
$$0\leq A(\gth_1)=-\myfrac{N-2}{\Gl_{\gb_*}-\gb_*}\left(\gb_*(N-1)-\Gl_{\gb_*}\right)\cot^2\gth_1<0
$$
which is again a contradictions.
\smallskip

\noindent{\it Step 4: End of the proof.} Since $r^2=\gb^2_*\gw^2+\gw^2_\gth$, $r_\gth=r(\gf_\gth-\gb_*)\tan\gf$, we have
$$rr_\gth=\left(\gb^2_*\gw+\gw_{\gth\gth}\right)\gw_\gth= r(\gf_\gth-\gb_*)\tan\gf.
$$
Since $\gw_\gth< 0$ on $(0,\frac{\gp}{2})$, it follows from Step 3 that $\gb^2_*\gw+\gw_{\gth\gth}\geq 0$ and thus
$$\myint{0}{\frac{\gp}{2}}\myfrac{\gb_*^2\gw+\gw_{\gth\gth}}{\gb^2_*\gw^2+\gw^2_{\gth}}\gw^2_{\gth}
\gw\cos\gth\sin^{N-2}\gth d\gth>0,
$$
since the integrand cannot be identically 0. The conclusion follows from $(\ref{beta-8})$. \qeda\medskip


\noindent\Remark Since $\gw_{\gth}(\frac{\gp}{2})=-c^2<0$, it follows $\gw(\gth)=-\gw_\gth(\gth)\cot\gth+O(\frac{\gp}{2}-\gth)$ as $\gth\to\frac{\gp}{2}$, and  from the eigenfunction equation $(\ref{beta-7})$
$$\myfrac{\gb_*^2\gw+\gw_{\gth\gth}}{\gb^2_*\gw^2+\gw^2_{\gth}}\gw^2_{\gth}=(\gb_*^2\gw+\gw_{\gth\gth})(1+o(1)).
$$
Therefore
$$-(p-1)\gw_{\gth\gth}=(\gb_*\Gl_{\gb_*}+(p-2)\gb_*^2+2-N)\gw (1+o(1))\quad\text{as }\gth\to\frac{\gp}{2}$$
and since $\Gd'\gw:=\gw_{\gth\gth}+(N-2)\cot\gth\,\gw_{\gth}$
$$-\Gd'\gw=\frac{\gb_*(\gb_*(2p-3)+p-N)+(p-2)(N-2)}{p-1}\gw (1+o(1))\quad\text{as }\gth\to\frac{\gp}{2}.
$$
Because $\gw$ is $C^\infty$ we obtain finally
\bel{Delta-1}
\abs{\Gd'\gw}\leq c\gw,
\ee
for some $c>0$.
\medskip


\end{document}